\documentclass[a4paper,12pt]{article}     

\usepackage{graphicx}		  
\usepackage{amsmath,amssymb}	

\usepackage{makeidx}

\newtheorem{lemma}{Lemma}[section]
\newtheorem{theorem}[lemma]{Theorem}
\newtheorem{proposition}[lemma]{Proposition}

\newtheorem{assumption}{Assumption}
\newtheorem{remark}[lemma]{Remark}

\newcommand\LABEL[1]{\label{#1}}

\def\authorfont{\footnotesize}

\def\ccode#1{\par		
	\vspace*{8pt}
	{\authorfont{\leftskip18pt\rightskip\leftskip
	\noindent #1\par}}\par}

\newenvironment{proof}{
\hspace*{-9mm}
{ \it Proof.}}
{\hfill {$\square$}\vspace{1.5em}}

\begin{document}

\begin{center}{
{\Large 
The structure of a minimal $n$-chart with\\
 two crossings II:\\
Neighbourhoods of $\Gamma_1\cup\Gamma_{n-1}$}
\vspace{10pt}
\\ 
Teruo NAGASE and Akiko SHIMA\footnote{The second author is supported by JSPS KAKENHI Grant Number 18K03309.}
}
\end{center}


\begin{abstract}
Given a 2-crossing minimal chart $\Gamma$, 
a minimal chart with two crossings, set\\
$\alpha=\min\{~i~|~$there exists an edge 
of label $i$ containing a white vertex$\}$, and\\
$\beta=\max\{~i~|~$there exists an edge 
of label $i$ containing a white vertex$\}$.\\
In this paper we study the structure of 
a neighbourhood of 
$\Gamma_\alpha\cup\Gamma_\beta$,
and propose a normal form 
for 2-crossing minimal $n$-charts,
here $\Gamma_\alpha$ and $\Gamma_\beta$ mean 
the union of all the edges of 
label $\alpha$ and $\beta$ respectively.
\end{abstract}
%
%
%
%
%

\ccode{2010 Mathematics Subject Classification. Primary 57Q45; Secondary 57Q35.}
\ccode{ {\it Key Words and Phrases}. surface link, chart, crossing. }

\setcounter{section}{0}
\section{Introduction}


Charts are oriented labeled graphs
in a disk 
with three kinds of vertices
called black vertices, crossings,
and white vertices (see Section~\ref{s:Prel} for the precise definition of charts, 
black vertices, crossings, and white vertices).
From a chart, we can construct an oriented closed surface 
embedded in 4-space ${\Bbb R}^4$ 
 (see \cite[chapter 14, chapter 18 and chapter 23]{BraidBook}). 
A C-move 
is a local modification between two charts
in a disk (see Section~\ref{s:Prel}).
A C-move between two charts induces 
an ambient isotopy between oriented closed surfaces 
corresponding to the two charts.
Two charts are said to be {\it C-move equivalent}\index{C-move~equivalent} 
if there exists
a finite sequence of C-moves 
which modifies one of the two charts 
to the other.

We will work in the PL or smooth category. 
All submanifolds are assumed to be locally flat.
A {\it surface link} is a closed surface embedded in 4-space ${\Bbb R}^4$. 
A {\it $2$-link} is a surface link each of whose connected component is a $2$-sphere.
A {\it $2$-knot}
is a surface link which is a $2$-sphere.
An orientable surface link is called a 
{\it ribbon surface link}
if there exists an immersion of a 3-manifold $M$
into ${\Bbb R}^4$ sending the boundary of $M$ onto the surface link
such that each connected component of $M$ is a handlebody
and its singularity
consists of ribbon singularities,
here a ribbon singularity
is a disk in the image of $M$
whose pre-image consists of 
two disks;
one of the two disks is a proper disk of $M$ 
and
the other is a disk in the interior of $M$.
In the words of charts,
a ribbon surface link is
a surface link corresponding to a {\it ribbon chart}, 
a chart C-move equivalent to 
a chart
without white vertices \cite{BraidThree}.
A chart is called a {\it $2$-link chart}
if a surface link corresponding to the chart is a $2$-link.

In this paper, 
we denote the closure, the interior, 
the boundary, and the complement of $(...)$ by 
$Cl(...)$, Int$(...)$, 
$\partial(...)$, $(...)^c$ 
respectively.
Also for a finite set $X$, 
the notation $|X|$ denotes 
the number of elements in $X$.

At the end of this paper 
there are lists of terminologies and notations 
which are used in this paper.

Kamada showed that 
any $3$-chart is a ribbon chart 
\cite{BraidThree}.
Kamada's result was extended by Nagase and Hirota:
Any $4$-chart with at most one crossing
is a ribbon chart \cite{NH}.
We showed that any $n$-chart with at most one crossing is a ribbon chart
\cite{OneCrossing}.
We also showed that any $2$-link chart 
with at most two crossings
 is a ribbon chart  \cite{TwoCrossingI},  
 \cite{TwoCrossingII}.

The purpose of our research
is a classification of surface links by using charts.
However we do not classify ribbon surface links.
It is known that any minimal chart $\Gamma$ decomposes a ribbon chart and a main part of a chart, ${\rm Main}(\Gamma)$
by C-I-M1 moves and C-I-M2 moves 
(see Section~\ref{s:Prel} for the definitions of 
${\rm Main}(\Gamma)$, C-I-M1 moves, and 
C-I-M2 moves).
We want to classify charts (modulo ribbon charts)
up to C-move equivalent
and we want to find a new minimal chart.

Let $\Gamma$ be a chart in a disk $D^2$, and 
$D$ a disk in $D^2$.
The pair $(\Gamma\cap D,D)$ is called 
a {\it tangle} provided that\index{tangle}
\begin{enumerate}
\item[(i)]
$\partial D$ does not contain
any white vertices, 
black vertices
nor crossings of $\Gamma$, 
\item[(ii)]
if an edge of $\Gamma$ intersects $\partial D$, 
then the edge intersects $\partial D$ transversely, 
\item[(iii)] $\Gamma\cap D\not=\emptyset$.
\end{enumerate}

For each label $m$ of a chart $\Gamma$, 
\index{$\Gamma_m$}we define

$$\Gamma_m=\text{ 
the union of 
all the edges of label $m$ and 
their vertices in }\Gamma.$$

For a chart $\Gamma$ 
containing a white vertex,  
we define

$\alpha(\Gamma)=\min\{~i~|~\Gamma_i$ 
contains a 
white vertex$\}$,\index{$\alpha(\Gamma)$}

$\beta(\Gamma)=\max\{~i~|~\Gamma_i$ 
contains a 
white vertex$\}$.\index{$\beta(\Gamma)$}

Let $\Gamma$ be a chart in a disk $D^2$. 
A closed edge in $\Gamma$ 
without vertices is called
a {\it hoop}.\index{hoop}
A hoop is said to be 
{\it simple}\index{simple hoop} 
if one of the complementary domains
of the hoop in $D^2$
does not contain any white vertices. 
An edge with two black vertices is called 
a {\it free edge}.\index{free edge}

In this paper and \cite{StI},
we investigate 
the structure of minimal charts with two crossings
(see Section~\ref{s:Prel}
for the precise definition of a minimal chart),
and give an enumeration of the charts 
with two crossings 
(see Section~\ref{s:NormalForm}). 
First, we split a minimal chart with two crossings 
into two kinds of tangles; 
one is called a net-tangle, and 
the other is called an IO-tangle. 

We investigate net-tangles in \cite{StI},
and
IO-tangles in this paper.
In short, for any minimal $n$-chart $\Gamma$ 
with two crossings in a disk $D^2$,
setting 
$\alpha=\alpha(\Gamma),\beta=\beta(\Gamma)$,
there exist 
two cycles 
$C_\alpha\subset\Gamma_\alpha$ and 
$C_\beta\subset\Gamma_\beta$
with $C_\alpha\cap C_\beta$ the two crossings 
(see Lemma~\ref{StIILemma1}).
If $\Gamma_\alpha$ or $\Gamma_\beta$
contains at least three white vertices,
then
after shifting all the free edges and 
simple hoops 
into a regular neighbourhood of $\partial D^2$
by applying C-I-M1 moves and C-I-M2 moves, 
we can find 
an annulus $A$ 
containing all the white vertices of $\Gamma$ 
but not intersecting 
any hoops nor free edges
such that (see Fig.~\ref{fig01}(a))
\begin{enumerate}
\item[(1)] 
each connected component of $Cl(D^2-A)$ 
contains a crossing,
\item[(2)]
$\Gamma\cap \partial A=
(C_\alpha\cup C_\beta)\cap \partial A$, and 
$\Gamma\cap \partial A$ 
consists of eight points.
\end{enumerate}
\begin{figure}
\begin{center}
\includegraphics{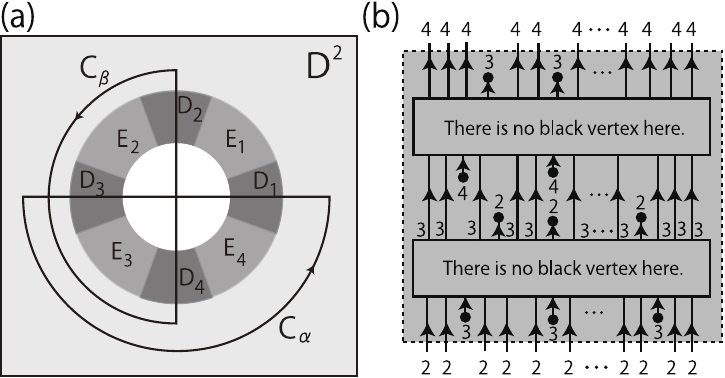}
\end{center}
\caption{\LABEL{fig01}
(b) a tangle $(\Gamma\cap E_i,E_i)$ with $\Gamma\cap E_i\subset \Gamma_2\cup \Gamma_3\cup \Gamma_4$
for the case $\alpha=1$ and $\beta=5$,
here all the free edges and simple hoops are in a regular neighbourhood of $\partial D^2$,
and the numbers are labels of the chart.
}
\end{figure}
In Section~\ref{s:NormalForm}, we show 
the annulus $A$ can be split 
into 
mutually disjoint four disks 
$D_1,D_2,D_3,D_4$ and 
mutually disjoint four disks 
$E_1,E_2,E_3,E_4$ such that
\begin{enumerate}
\item[(3)] 
for each $i=1,2,3,4$, 
if $\alpha+1=\beta-1$, 
then $\Gamma\cap E_i$ consists of 
parallel proper arcs of $E_i$ of label $\alpha+1$, 
otherwise 
the tangle $(\Gamma\cap E_i,E_i)$ is 
a net-tangle with 
$\Gamma\cap E_i\subset
\cup^{\beta-1}_{j=\alpha+1}\Gamma_j$
as shown in Fig.~\ref{fig01}(b) 
(see Section~\ref{s:Indices} for the definition of 
a net-tangle), and
\item[(4)] 
for each $i=1,3$ (resp. $i=2,4$)
the tangle $(\Gamma\cap D_i,D_i)$ 
is an IO-tangle of label $\alpha$ 
(resp. label $\beta$) 
(see two paragraphs before 
Theorem~\ref{StIITheorem2} 
for the definition of an IO-tangle).
\end{enumerate}
We count the number of edges between edges 
with black vertices
in Fig.~\ref{fig01}(b) 
to enumerate charts with two crossings.
As important results, 
from the enumeration 
we can calculate the fundamental group of 
the exterior of the surface link 
represented by $\Gamma$, 
and the braid monodromy of the surface braid 
represented by $\Gamma$.

A chart with exactly two crossings is called\index{2-crossing chart}
a {\it $2$-crossing chart}.

The following lemma is 
the first step to investigate 
a 2-crossing minimal chart 
and shown in Section~\ref{s:ProofLemma1}.

\begin{lemma}
\LABEL{StIILemma1} 
{\rm$($cf. \cite[Lemma~6.3]{TwoCrossingII}$)$}
Let $\Gamma$ be a $2$-crossing minimal chart in a disk $D^2$.
Set $\alpha=\alpha(\Gamma),
\beta=\beta(\Gamma)$.
Then there exists a minimal chart $\Gamma'$ obtained from $\Gamma$ by applying C-I-M1 moves and C-I-M2 moves
satisfying the following conditions.
\begin{enumerate}
\item[{\rm (a)}]
There exist two cycles 
$C_\alpha,C_\beta$ with $C_\alpha\subset\Gamma_\alpha',C_\beta\subset\Gamma_\beta'$ such that 
$C_\alpha\cap C_\beta$ consists of the two crossings.
\item[{\rm (b)}]
There exists an annulus $A$ 
with $A\cap \partial D^2=\emptyset$ 
such that 
$A$ contains all the white vertices of $\Gamma'$ 
but does not intersect hoops nor free edges.
\item[{\rm (c)}]
Each connected component of $Cl(D^2-A)$ contains a crossing.
\item[{\rm (d)}]
 $\Gamma'\cap \partial A=(C_\alpha\cup C_\beta)\cap \partial A$,
$\Gamma'\cap \partial A$ consists of eight points.
\item[{\rm (e)}] 
$\Gamma_\alpha'\cap A$ 
consists of 
two connected components 
$X_1,X_3$ 
separated by $C_\beta$. 
\item[{\rm (f)}]
$\Gamma_\beta'\cap A$ 
consists of 
two connected components 
$X_2,X_4$ 
separated by $C_\alpha$. 
\end{enumerate}
\end{lemma}

Let $\Gamma$ be a chart in a disk $D^2$, 
and $E$ a disk. 
Suppose that an edge $e$ of $\Gamma$ 
transversely intersects $\partial E$.
Let $p$ be a point in 
$e\cap \partial E$, 
and 
$N$ a regular neighbourhood of $p$ 
in $D^2$.
Then the orientation of $e$
induces the one of the arc $e\cap N$.
The edge $e$ is said to be 
{\it locally inward}\index{locally inward} 
(resp. {\it locally outward})\index{locally outward}
at $p$ 
with respect to $E$
if the oriented arc $e\cap N$ is 
oriented
from a point outside 
(resp. inside) $E$ 
to a point inside 
(resp. outside) $E$.
We often say that 
$e$ is locally inward (resp. outward) at $p$ 
instead of saying that
$e$ is locally inward (resp. outward) at $p$ 
{\it with respect to $E$}, 
if there is no confusion.

An edge of a chart $\Gamma$ is called 
a {\it terminal edge}\index{terminal edge}
if it contains a white vertex and 
a black vertex.

For a simple arc $\ell$, we set\\
\ \ \ $
\begin{array}{rl}
\partial \ell&= \text{ the set of its two endpoints, 
and}\\
{\rm Int}~\ell&=~\ell-\partial \ell.
\end{array}
$

Let $\Gamma$ be a chart, and 
$m$ a label of the chart.
A tangle $(\Gamma\cap D,D)$
is called an 
{\it IO-tangle of label $m$}\index{IO-tangle}
provided that (see Fig.~\ref{fig02})
\begin{enumerate}
\item[(i)]
no terminal edge nor free edge intersects $\partial D$, 
\item[(ii)]
there exists a label $k$ with $|m-k|=1$ and
$\Gamma\cap D\subset\Gamma_m\cup\Gamma_k$,
\item[(iii)]
there exist two arcs $L_I,L_O$ on $\partial D$
with $L_I\cap L_O=\partial L_I
=\partial L_O=\Gamma_m\cap \partial D$,
\item[(iv)]
for any point 
$p\in\Gamma\cap$~Int~$L_I$,
there exists an edge 
of label $k$ 
locally inward at $p$, 
and\\
for any point 
$p\in\Gamma\cap$~Int~$L_O$,
there exists an edge 
of label $k$ 
locally outward 
at $p$.
\end{enumerate}
The pair $(L_I,L_O)$ is called a 
{\it boundary IO-arc pair}\index{boundary IO-arc pair}
of the tangle.
An IO-tangle of label $m$ is 
said to be {\it simple}\index{simple IO-tangle}
if all the terminal edge in $D$ is of label $m$.

\begin{figure}
\begin{center}
\includegraphics{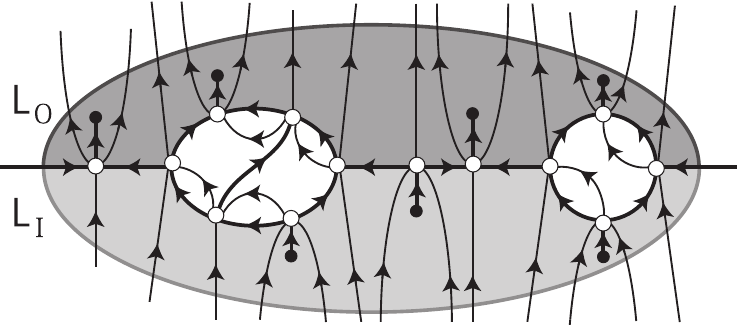}
\end{center}
\caption{\LABEL{fig02}
An IO-tangle of label $m$. The thick edges are of label $m$. The light gray arc is $L_I$, and the dark gray arc is $L_O$}
\end{figure}

Let $X_1,X_2,X_3,X_4$ be 
the connected components in 
Lemma~\ref{StIILemma1} above.
For each $i=1,2,3,4$,
let $D_i^*$ be 
a regular neighbourhood of 
the SC-closure $SC(X_i)$ in 
the annulus $A$ 
(see Section~\ref{s:ProofLemma1} 
for the definition of 
SC-closures).
Then the tangle $(\Gamma'\cap D_i^*,D_i^*)$
is called a {\it fundamental tangle}\index{fundamental tangle}
 of the $2$-crossing chart.

We obtain the following theorem for fundamental tangles.

\begin{theorem}
\LABEL{StIITheorem2} 
Let $\Gamma$ be 
a $2$-crossing minimal chart, 
and $(\Gamma\cap D,D)$ 
a fundamental tangle. 
If $D$ contains 
at least two white vertices,
then the tangle
is a simple IO-tangle.
\end{theorem}

A simple arc $\ell$ in a compact surface $E$ 
is called a {\it proper} arc\index{proper arc} 
provided that 
$\ell\cap \partial E=\partial\ell$.
We decompose a simple IO-tangle.
Let $\Gamma$ be a chart, and 
$m$ a label of $\Gamma$.
An IO-tangle 
$(\Gamma\cap D,D)$ of label $m$ 
is called 
a {\it Type-I elementary IO-tangle 
of label $m$}\index{Type-I}
provided that
\begin{enumerate}
\item[(I)] 
there exists a disk $E$ 
containing 
all the white vertices in $D$ 
with $\partial E\subset\Gamma_m\cap D$.
\end{enumerate}
An IO-tangle 
$(\Gamma\cap D,D)$ of label $m$ 
is called 
a {\it Type-II elementary IO-tangle 
of label $m$}\index{Type-II}
provided that
\begin{enumerate}
\item[(II)]
$\Gamma_m\cap D$ contains no cycle,
but there exists 
a proper arc $L_D$ of $D$ 
containing 
all the white vertices in $D$ 
with $L_D\subset\Gamma_m\cap D$
(hence 
the arc $L_D$ connects
the two points 
$\Gamma_m\cap\partial D$).
\end{enumerate}
Let $D^\dagger$ be a disk. 
A tangle 
$(\Gamma\cap D^\dagger,D^\dagger)$
is called a {\it trivial tangle}\index{trivial tangle} 
of label $m$
if $\Gamma\cap D^\dagger$ 
is a proper arc of $D^\dagger$ 
contained in the interior of 
an edge of label $m$. 
We consider the trivial tangle 
as a Type-II elementary IO-tangle.
Now Type-I elementary IO-tangles and
Type-II elementary IO-tangles\index{elementary IO-tangles} 
are called {\it elementary IO-tangles}.

Let $\Gamma$ be a chart, and 
$m$ a label of the chart.
Let $(\Gamma\cap D,D)$ be 
an IO-tangle of label $m$.
The tangle is said to have an 
{\it elementary IO-tangle decomposition}\index{decomposition}
$(\Gamma\cap D_1,D_1)\#(\Gamma\cap D_2,D_2)\#
\cdots\#(\Gamma\cap D_{2s+1},D_{2s+1})$
provided that 
(see Fig.~\ref{fig15})
\begin{enumerate}
\item[(i)]
for each $i=1,2,\cdots,2s+1$,\\
if $i$ is even, then
$(\Gamma\cap D_i,D_i)$ is 
a Type-I elementary IO-tangle of label $m$,\\
otherwise 
$(\Gamma\cap D_i,D_i)$ is 
a Type-II elementary IO-tangle of label $m$,
\item[(ii)]
$D=\cup_{i=1}^{2s+1}D_i$,
\item[(iii)]
there exist $2s$ proper arcs 
$\ell_1,\ell_2,\cdots,\ell_{2s}$ 
of $D$
such that 
for each $i=1,2,\cdots,2s$,
$\ell_i\cap \Gamma
=\ell_i\cap \Gamma_m
=\text{a point}$, 
\item[(iv)] 
for each $1\le i<j\le 2s+1$,\\
$
D_i\cap D_j=
\left\{
\begin{array}{lll}
\ell_i&\ \ \ &\text{if $j=i+1$},\\ 
\emptyset& &\text{otherwise}.
\end{array}
\right.
$
\end{enumerate}

\begin{theorem}
\LABEL{StIITheorem3} 
Let $\Gamma$ be 
a minimal chart, and
$m$ a label of the chart. 
Then any simple IO-tangle 
$(\Gamma\cap D,D)$
of label $m$ 
has an elementary IO-tangle decomposition.
\end{theorem}

Let $\Gamma$ be a chart.
If an edge is oriented 
from a vertex $v_1$ 
to a vertex $v_2$,
then 
edge is said to be 
{\it inward}\index{inward} at the vertex $v_2$, 
also 
{\it outward}\index{outward} at the vertex $v_1$.

Let $\Gamma$ be a chart, 
and $m$ a label of $\Gamma$. 
Let $p$ be a positive integer.
A Type-I elementary IO-tangle 
$(\Gamma\cap D,D)$ of label $m$ 
with a boundary IO-arc pair $(L_I,L_O)$ 
is said to be of {\it Type-I$_p$}\index{Type-I$_p$}
provided that (see Fig.~\ref{fig03})
\begin{enumerate}
\item[(i)] 
there exists a disk $E$ 
containing 
all the white vertices in $D$ 
with 
$\partial E\subset\Gamma_m\cap D$, 
but 
Int~$E$ does not contain a white vertex,
\item[(ii)] 
the closure of $\Gamma_m\cap (D-E)$ 
consists of 
two terminal edges $\tau_I,\tau_O$ 
and two arcs such that 
one of the two arcs is contained in 
an edge $e_I$ inward at 
a white vertex in $\partial E$, and 
the other is contained in 
an edge $e_O$ outward at 
a white vertex in $\partial E$,

\item[(iii)] 
$\partial E-(e_I\cup e_O)$ 
consists of
two components $J_I',J_O'$ with 
$J_I'\cap \tau_I\neq\emptyset,
J_O'\cap \tau_O\neq\emptyset$
such that
each of $J_I'$ and $J_O'$ contains 
exactly $p$ white vertices,
\item[(iv)] 
one of the two components of 
$J_I'-\tau_I$ 
does not contain a white vertex, 
and\\
one of the two components of 
$J_O'-\tau_O$ 
does not contain a white vertex.
\end{enumerate}

\begin{figure}
\begin{center}
\includegraphics{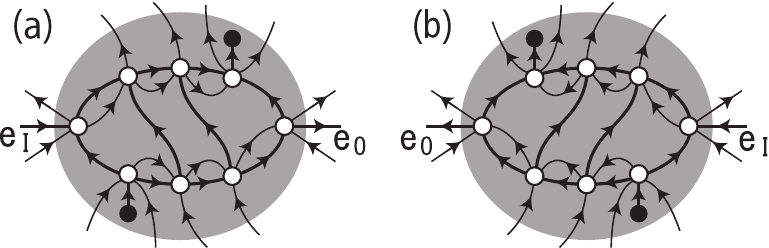}
\end{center}
\caption{\LABEL{fig03}
Type-I$_3$ elementary IO-tangles 
of label $m$.
The thick arcs are of label $m$.}
\end{figure}

The following theorem simplifies 
the structure of 
fundamental tangles 
for a 2-crossing minimal chart.

\begin{theorem}\LABEL{StIITheorem4} 
Let $\Gamma$ be a minimal chart, and 
$m$ a label of $\Gamma$. 
Then 
any Type-I elementary IO-tangle 
of label $m$ is of Type-I$_p$ 
for some integer $p$.
\end{theorem}

Our paper is organized as follows:
In Section~\ref{s:Prel}, 
we introduce 
the definition of charts and its related words.
In Section~\ref{s:ProofLemma1},
we shall prove Lemma~\ref{StIILemma1}.
In Section~\ref{s:Th2Th3},
we shall prove Theorem~\ref{StIITheorem2}
and Theorem~\ref{StIITheorem3}.
In Section~\ref{s:OneWayPath} and 
Section~\ref{s:Principal},
we investigate a directed path.
In Section~\ref{s:Th4},
we shall prove Theorem~\ref{StIITheorem4}.
In Section~\ref{s:Indices}, 
we define indices for simple IO-tangles and 
N-tangles.
In Section~\ref{s:NormalForm}, 
we define a normal form for 
$2$-crossing minimal charts.



\section{Preliminaries}
\LABEL{s:Prel}

In this section, 
we introduce 
the definition of charts and its related words.

Let $n$ be a positive integer.
An $n$-{\it chart}  
(a braid chart of degree $n$ \cite{KS}
or a surface braid chart of degree $n$ \cite{BraidBook}) 
is 
an oriented labeled graph in the interior of a disk,
which may be empty 
or
have closed edges without vertices
satisfying the following four conditions
(see Fig.~\ref{fig04}):
\begin{enumerate}
\item[(i)] 
Every vertex has degree $1$, $4$, or $6$.
\item[(ii)] 
The labels of edges are 
in $\{1,2,\dots,n-1\}$.
\item[(iii)]
In a small neighborhood of
each vertex of degree $6$,
there are six short arcs,
three consecutive arcs are
oriented inward 
and
the other three are outward,
and
these six are labeled $i$ and $i+1$
alternately for some $i$,
where the orientation and label of
each arc are inherited from
the edge containing the arc.
\item[(iv)]
For each vertex of degree $4$,
diagonal edges have the same label
and
are oriented coherently,
and the labels $i$ and $j$ of
the diagonals satisfy $|i-j|>1$.
\end{enumerate}
We call a vertex of degree $1$ a {\it black vertex},
a vertex of degree $4$ a {\it crossing}, and 
a vertex of degree $6$ a {\it white vertex}
respectively.
Among six short arcs
in a small neighborhood of
a white vertex,
a central arc of each three consecutive arcs
oriented inward (resp. outward) 
is called a   
{\it middle arc}\index{middle arc} at the white vertex
(see Fig.~\ref{fig04}(c)).
For each white vertex $v$, 
there are two middle arcs at $v$ 
in a small neighborhood of $v$.
An edge $e$ is said to be 
{\it middle at} a white vertex $v$\index{middle at $v$} if it 
contains a middle arc at $v$. 



\begin{figure}[htb]
\begin{center}
\includegraphics{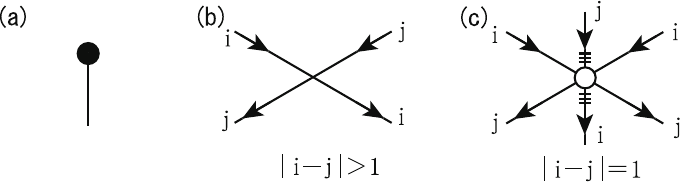}
\end{center}
\caption{ \LABEL{fig04} (a) A black vertex. (b) A crossing. (c) A white vertex. 
Each arc with three transversal short arcs is a middle arc at the white vertex. }
\end{figure}

Now {\it C-moves} are local modifications 
of charts as shown in Fig.~\ref{fig05}
(cf. \cite{KS}, 
\cite{BraidBook} and \cite{Tanaka}).
We often use C-I-M2 moves, C-I-M3 moves, C-II moves
and C-III moves.

Kamada originally defined CI-moves
as follows: 
A chart $\Gamma$ is obtained from
a chart $\Gamma'$ in a disk $D^2$
by a {\it CI-move},
if there exists a disk $E$ 
in $D^2$ such that 
\begin{enumerate}
\item[(i)] 
the two charts $\Gamma$ and $\Gamma'$
intersect the boundary of $E$ transversely
or
do not intersect the boundary of $E$, 
\item[(ii)] 
$\Gamma\cap E^c=\Gamma'\cap E^c$, and
\item[(iii)]
neither $\Gamma\cap E$ nor 
$\Gamma'\cap E$ 
contains a black vertex,
\end{enumerate}
where $E^c$ is 
the complement of $E$ in the disk $D^2$.

\begin{remark}
{\rm Any CI-move is realized by a finite sequence of 
seven types: C-I-R2, C-I-R3, C-I-R4, 
C-I-M1, C-I-M2, C-I-M3, C-I-M4.}
\end{remark}

\begin{figure}[htb]
\begin{center}
\includegraphics{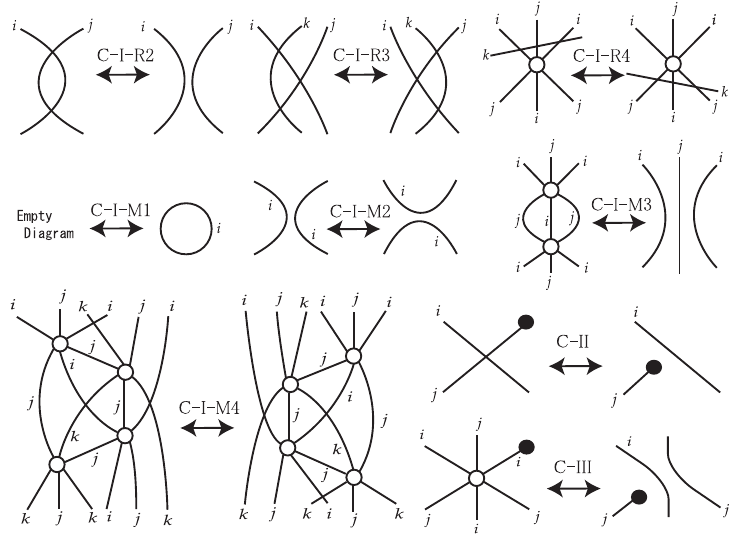}
\end{center}
\caption{ \LABEL{fig05} For the C-III move, the edge containing the black vertex does not contain a middle arc at
a white vertex in the left figure. }
\end{figure}

Let $\Gamma$ be a chart. 
Let $e_1$ and $e_2$ be edges of $\Gamma$
which connect two white vertices $w_1$ and $w_2$
where possibly $w_1=w_2$.
Suppose that 
the union $e_1\cup e_2$ bounds 
an open disk $U$.
Then $Cl(U)$ 
is called 
a {\it bigon} of $\Gamma$\index{bigon}
provided that
any edge containing $w_1$ or $w_2$ 
does not intersect the open disk $U$
(see Fig.~\ref{fig06}).
Note that neither $e_1$ nor $e_2$ contains a crossing.

\begin{figure}
\begin{center}
\includegraphics{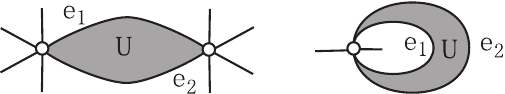}
\end{center}
\caption{ \LABEL{fig06} Bigons.}
\end{figure}
%

Let $\Gamma$ be a chart.
Let 
$w(\Gamma),~f(\Gamma),~c(\Gamma)$, 
and 
$b(\Gamma)$ be 
the number of white vertices, 
the number of free edges, 
the number of crossings, 
and 
the number of bigons of $\Gamma$
respectively.
The 4-tuple $(c(\Gamma),w(\Gamma),-f(\Gamma),-b(\Gamma))$ is called a 
{\it $c$-complexity} of the chart $\Gamma$.
The 4-tuple $(w(\Gamma),c(\Gamma),-f(\Gamma),-b(\Gamma))$ is called a 
{\it $w$-complexity} of the chart $\Gamma$.
The 3-tuple $(c(\Gamma)+w(\Gamma),-f(\Gamma),-b(\Gamma))$ is called a 
{\it $cw$-complexity} of the chart $\Gamma$
(see \cite{BraidThree} 
for complexities of charts).

A chart $\Gamma$ is said to be 
{\it $c$-minimal $($resp. $w$-minimal or $cw$-minimal$)$} if
its $c$-complexity (resp. $w$-complexity or $cw$-complexity) is minimal among the charts 
which are C-move equivalent to 
the chart $\Gamma$
with respect to 
the lexicographical order of the 
4-tuple (or 3-tuple) of the integers.
If a chart is $c$-minimal, $w$-minimal or $cw$-minimal, 
then we say that the chart is {\it minimal}\index{minimal chart}
in this paper.

An {\it oval nest} is a free edge 
together\index{oval nest} 
with some concentric simple hoops. 

\begin{proposition}
\LABEL{MoveOutFreeEdge}
{\rm (\cite[Proposition~2.2]{StI})}
Let $\Gamma$ be a chart in a disk $D^2$.
For any regular neighbourhood $N$ 
of $\partial D^2$ in $D^2$,
there exists a chart $\Gamma'$ obtained from $\Gamma$ 
by C-I-M2 moves and ambient isotopies of $D^2$ 
without changing  
the complexity such that $($see Fig.~\ref{fig07}$)$
\begin{enumerate}
\item[{\rm (a)}]
$\Gamma'\cap (D^2-N)$ contains no free edge,
\item[{\rm (b)}]
$\Gamma'\cap N$ consists of oval nests, 
simple hoops and free edges.
\end{enumerate}
\end{proposition}

\begin{figure}
\begin{center}
\includegraphics{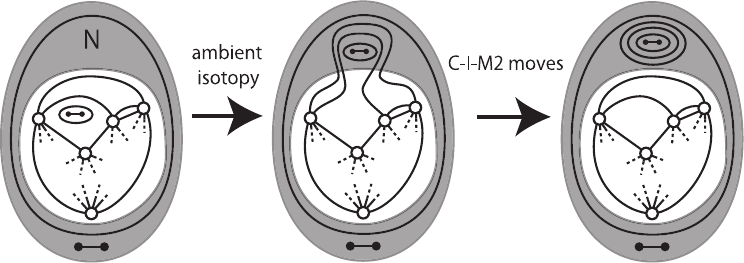}
\end{center}
\caption{ \LABEL{fig07} Moving free edges into 
a gray annulus $N$, a regular neighbourhood of $\partial D^2$ in $D^2$.}
\end{figure}
%


Let $m$ be a label of a chart $\Gamma$.
A simple closed curve in $\Gamma_m$
is called a {\it ring},\index{ring}
if it contains at least one crossing 
but does not contain a white vertex
nor black vertex.

\begin{proposition}
\LABEL{Assumption0}
{\rm (\cite[Remark~2.3]{MinimalChart}, 
\cite[Proposition~2.3]{StI})}
Let $\Gamma$ be a minimal chart in $D^2$. 
Then we have 
the following:
\begin{enumerate}
\item[{\rm (a)}]
If an edge of $\Gamma$ contains a black vertex, 
then the edge is a terminal edge or a free edge.
\item[{\rm (b)}]
Any terminal edge of $\Gamma$ contains a middle arc 
at its white vertex.
\item[{\rm (c)}] Each complementary domain 
of any ring in $D^2$ 
must contain at least one white vertex.
\end{enumerate}
\end{proposition}


\begin{proposition}
\LABEL{MoveOutSimpleHoop}
{\rm (\cite[Proposition~2.4]{StI})}
Let $\Gamma$ be a minimal chart in a disk $D^2$.
For any regular neighbourhood $N$ 
of $\partial D^2$ in $D^2$,
there exists a minimal chart $\Gamma'$ 
obtained from $\Gamma$ 
by C-I-M2 moves 
and ambient isotopies of $D^2$ such that
\begin{enumerate}
\item[{\rm (a)}]
$\Gamma'\cap (D^2-N)$ contains neither free edge
nor simple hoop,
\item[{\rm (b)}]
$\Gamma'\cap N$ consists of oval nests, 
simple hoops and free edges.
\end{enumerate}
\end{proposition}

For any minimal chart in a disk $D^2$
we can move free edges and simple hoops into 
a regular neighbourhood of $\partial D^2$ 
in $D^2$ 
by C-I-M2 moves and ambient isotopies of $D^2$
by Proposition~\ref{MoveOutFreeEdge}
and Proposition~\ref{MoveOutSimpleHoop}.
Even during argument,
if free edges or simple hoops appear, 
we immediately move them 
into a regular neighbourhood of $\partial D^2$ 
in $D^2$.
Thus we assume the following 
(cf. \cite{OneCrossing}, \cite[Assumption 1]{MinimalChart}):

\begin{assumption}
\LABEL{AssumptionFreeEdge}
{\it For any minimal chart in a disk $D^2$, 
all the free edges and 
simple hoops  
are in a regular neighbourhood of 
$\partial D^2$ in $D^2$.}
\end{assumption}

Let $\Gamma$ be a minimal chart in 
a disk $D^2$, and 
$X$ the union of all the free edges 
and simple hoops.
Now $X$ is in a regular neighbourhood $N$ of 
$\partial D^2$ in $D^2$ 
by Assumption~\ref{AssumptionFreeEdge}.
We say that $\Gamma$
is a chart with a {\it brim} 
$N$ in $D^2$.\index{brim}
We define \index{${\rm Main}(\Gamma)$}
$${\rm Main}(\Gamma)=\Gamma-X.$$
Let $\widehat D=Cl(D^2-N)$.  
Then $\Gamma\cap\widehat D=$Main$(\Gamma)$.
Hence $(\Gamma\cap\widehat D,\widehat D)$ is 
a tangle without free edges 
and simple hoops.

In this paper we always assume that
\begin{enumerate}
\item[] 
{\it for any tangle $(\Gamma\cap D,D)$, 
the disk $D$ does not contain 
any free edge nor a simple hoop}.
\end{enumerate}


Let $E$ be a disk, and
$\ell_1,\ell_2,\ell_3$ three arcs on $\partial E$
such that each of $\ell_1\cap \ell_2$ and $\ell_2\cap \ell_3$ is one point and $\ell_1\cap \ell_3=\emptyset$
(see Fig.~\ref{fig08}(a)),
say $p=\ell_1\cap \ell_2$,
$q=\ell_2\cap \ell_3$.
Let $\Gamma$ be a chart in a disk $D^2$.
Let $e_1$ be a terminal edge of 
 $\Gamma$. 
A triplet $(e_1,e_2,e_3)$ of 
mutually different edges of $\Gamma$
is called 
a {\it consecutive triplet}\index{consecutive triplet}
if there exists
a continuous map $f$ from the disk $E$ 
to the disk $D^2$ such that (see Fig.~\ref{fig08}(b) and (c))
\begin{enumerate}
\item[(i)] the map $f$ is injective on $E-\{p,q\}$,
\item[(ii)] 
$f(\ell_3)$ is an arc in $e_3$, and $f({\rm Int}~E)\cap\Gamma=\emptyset$,
$f(\ell_1)=e_1$,
$f(\ell_2)=e_2$,
\item[(iii)]
each of $f(p)$ and $f(q)$ is a white vertex.
\end{enumerate}
If the label of $e_3$ is different
from the one of $e_1$ 
then the consecutive triplet is said to be
{\it admissible}.


\begin{remark}
{\rm Let $(e_1,e_2,e_3)$
be a consecutive triplet. 
Since $e_2$ is an edge of $\Gamma$, 
the edge $e_2$ MUST NOT contain a crossing.}
\end{remark}

\begin{figure}
\begin{center}
\includegraphics{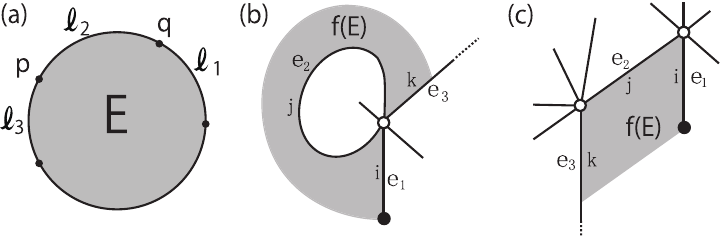}
\end{center}
\caption{ \LABEL{fig08} (b),(c) Consecutive triplets.}
\end{figure}


\begin{lemma}
\LABEL{ConsecutiveTripletLemma} 
{\rm [Consecutive Triplet Lemma]}
{\rm $($\cite[Lemma 1.1]{OneCrossing}, \cite[Lemma 3.1]{MinimalChart}$)$}
{\it Any consecutive triplet 
in a minimal chart is admissible.}
\end{lemma}

Let $\Gamma$ be a chart.
A tangle $(\Gamma\cap D,D)$ is 
called an {\it NS-tangle of label $m$} 
(new significant tangle) 
provided that\index{NS-tangle}
\begin{enumerate}
\item[(i)] if $i\neq m$, 
then $\Gamma_i\cap \partial D$ is 
at most one point,
\item[(ii)] 
$\Gamma\cap D$ contains at least one white vertex, and 
\item[(iii)]
for each label $i$, 
the intersection $\Gamma_i\cap D$ contains 
at most one crossing.
\end{enumerate}

\begin{lemma} 
{\rm $($\cite[Theorem 1.2]{MinimalChart}$)$ }
\LABEL{LemNS-Tangle}
In a minimal chart, 
there does not exist 
an NS-tangle of any label.
\end{lemma}

\begin{lemma}\LABEL{BoundaryConditionLemma}
{\rm [Boundary Condition Lemma]
$($\cite[Lemma 4.1]{TwoCrossingI}, 
\cite[Lemma 11.1]{MinimalChart}$)$}
Let $(\Gamma\cap D,D)$ be a tangle 
in a minimal chart $\Gamma$ 
such that 
$D$ does not contain any crossing,
free edge nor simple hoop.
Let $a=\min\{~i~|~\Gamma_i\cap\partial D\not=\emptyset\}$ and
$b=\max\{~i~|~\Gamma_i\cap\partial D\not=\emptyset\}$.
Then
 $\Gamma_i\cap D=\emptyset$ 
 except for $a\le i \le b$.
\end{lemma}

Let $\Gamma$ be a chart, and 
$m$ a label of $\Gamma$. 
A simple closed curve in $\Gamma_m$ is\index{cycle} 
called a {\it cycle of label $m$}. 

Let $\Gamma$ be a chart,
and $m$ a label of $\Gamma$.
Let $C$ be a cycle of label $m$ in $\Gamma$ 
bounding a disk $E$.
Then an edge $e$ of 
\underline{label $m$} 
is called
an {\it outside edge for $C$} provided that\index{outside edge}
\begin{enumerate}
\item[(i)]
$e\cap C$ consists of one white vertex or two white vertices, and
\item[(ii)]
$e\not\subset E$.
\end{enumerate}
For a cycle $C$ of label $m$,\index{${\mathcal{W}}_O^{{\rm Mid}}(C,m)$} 
we define
$$\begin{array}{ll}
{\mathcal{W}}(C)&= \{ w \ | \text{ $w$ is a white vertex in $C$} \},\\
{\mathcal{W}}_O^{{\rm Mid}}(C,m)&= \{ w\in \mathcal{W}(C)  \ | \text{ there exists an outside edge for $C$ {\it middle} at $w$} \}.
\end{array}$$

The following lemma will be used in
the proof of Lemma~\ref{Lemma4}.

\begin{lemma} 
{\rm $($\cite[Lemma 2.11]{StI}$)$} 
\LABEL{LemTwoColorTangle}
Let $\Gamma$ be a minimal chart, and 
$m,k$ labels of $\Gamma$ with $|m-k|=1$. 
Let $(\Gamma\cap D,D)$ be 
a tangle
with $\Gamma\cap D\subset\Gamma_m\cup\Gamma_k$
but without free edges nor simple hoops. 
Then for any cycle $C$ of label $m$ in $D$, 
we have 
$|{\mathcal W}_O^{{\rm Mid}}(C,m)|\ge 2$.
\end{lemma}



\section{Proof of Lemma~1.1}
\LABEL{s:ProofLemma1}

Let $\Gamma$ be a minimal chart 
with a brim $N$ 
in a disk $D^2$.
By Assumption~\ref{AssumptionFreeEdge},
all the simple hoops and free edges of $\Gamma$
are in the brim $N$.
Set $C^*=\partial N-\partial D^2$.
Let $e$ be an edge in Main($\Gamma$) 
of label $m$ 
such that 
there exists a simple arc $\ell$ in $Cl(D^2-N)$ 
connecting 
a point $p$ in Int~$e$ and 
a point $q$ in $C^*$ 
with $\ell\cap{\rm Main}(\Gamma)=p$ and 
$\ell\cap N=q$
(see Fig.~\ref{fig09}(a)).
We construct a chart from $\Gamma$
by a C-I-M1 move and C-I-M2 moves as follows.

First we create a simple hoop $H$ of label $m$ 
surrounding the point $q$
by a C-I-M1 move (see Fig.~\ref{fig09}(b)),
where $H$ is oriented in such a way that 
we can apply a C-I-M2 move 
between $e$ and $H$.
Next apply a C-I-M2 move to 
the hoop $H$ along $C^*$
(see Fig.~\ref{fig09}(c) and (d)).
Then we obtain two simple hoops 
parallel to $C^*$;
one hoop $H_1$ is in the brim $N$ and 
the other hoop $H_2$ is in $D^2-N$.
Finally apply a C-I-M2 move 
between the edge $e$ and $H_2$
along the arc $\ell$ 
to get a new edge $e^*$ of label $m$ 
(see Fig.~\ref{fig09}(e)). 
Let
$$\Gamma^*=(\Gamma-e)\cup e^*\cup H_1.$$
Then $\Gamma^*$ is a chart 
C-move equivalent to $\Gamma$.
We say that
the chart $\Gamma^*$ is obtained from 
$\Gamma$ by\index{DH-trick}
{\it a double hoops trick $($DH-trick$)$
along the arc $\ell$}.

Let $\Gamma$ be a minimal chart 
with a brim $N$ 
in a disk $D^2$.
Let $\ell$ be a simple arc in $Cl(D^2-N)$
connecting 
a point $p$ in $D^2-N$ and 
a point $q$ in $\partial N-\partial D^2$
with $\ell \cap N=q$ 
such that 
the arc $\ell$ transversely intersects 
an edge of $\Gamma$ at each point in 
Main$(\Gamma)\cap$~Int~$\ell$. 
Set Main$(\Gamma)\cap\ell=\{v_1,v_2,\cdots,v_s\}$ 
here we assume that 
$p,v_s,v_{s-1},\cdots, v_2,v_1,q$ 
are situated on $\ell$ in this order, 
here possibly $p=v_s$. 
For each $i=1,2,\cdots,s$ 
let $\ell[v_i,q]$ be the subarc of $\ell$ 
with $\partial \ell[v_i,q]=\{v_i,q\}$.
Let $\Gamma^0=\Gamma$. 
For each $i=1,2,\cdots,s$ 
let $\Gamma^i$ be a chart 
obtained from $\Gamma^{i-1}$ 
by performing a DH-trick along $\ell[v_i,q]$ 
with $\ell[v_{i+1},q]\cap \Gamma^i=v_{i+1}$ 
for $i<s$. 
Then we say that {\it $\Gamma^s$ is 
obtained from $\Gamma$\index{DH-trick} 
by DH-tricks along $\ell$}.  

\begin{remark}
{\rm  In the definition of 
DH-tricks along the arc $\ell$ above,  
if the chart $\Gamma$ is a minimal chart, 
we have the following.\\
(1) If the point $p$ is not in a bigon, 
then 
the number of bigons does not change 
by the DH-tricks along $\ell$ 
(see Fig.~\ref{fig09}(f)), 
and neither does the complexity.\\
(2) If two bigons intersect by an edge, 
then we can eliminate the two white vertices of 
the bigons by C-I-M2 moves and a C-I-M3 move. 
Thus applying C-I-M2 moves from outer bigons 
(cf. Fig.~\ref{fig07}),
we can assume that the interior of each bigon
does not intersect $\Gamma$.
}
\end{remark}

\begin{figure}
\begin{center}
\includegraphics{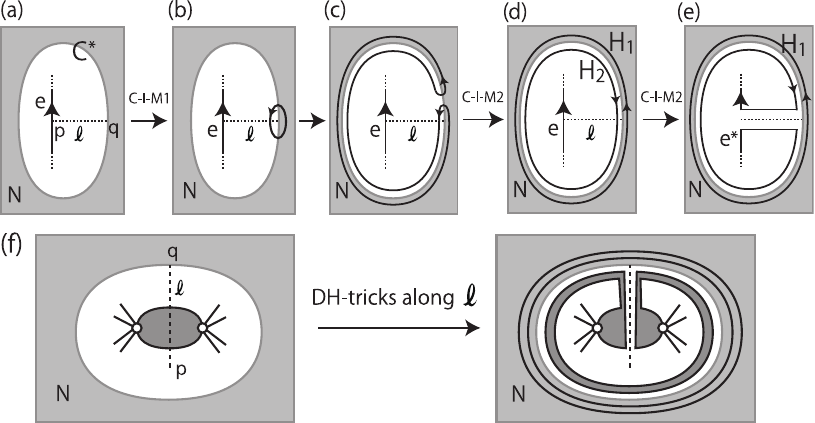}
\end{center}
\caption{ \LABEL{fig09} Light gray is a brim $N$ in $D^2$. (f) Dark gray is a bigon.}
\end{figure}

%

\begin{lemma}\LABEL{Lemma3-1}
Let $\Gamma$ be a minimal chart 
in a disk $D^2$.
If $2\leq c(\Gamma)\le3$,
then any hoop is simple.
\end{lemma}

\begin{proof}
Let $N$ be a brim of $D^2$. 
Suppose that there exists 
a non-simple hoop $C$. 
Then each complementary domains 
of the hoop in $D^2$ 
contains a white vertex.
Since there exist at most three crossings, 
one of the two complementary domains contains 
at most one crossing,
say $U$.
The closure $Cl(U)$ is a disk 
with $\partial U=C$ 
or 
an annulus containing $\partial D^2$ 
with 
$\partial U=C\cup\partial D^2$.
There are two cases.

{\bf Case 1.} 
If $Cl(U)$ is a disk, then 
let $D$ be 
a regular neighbourhood of $Cl(U)$ in $D^2$.
Since $\Gamma\cap \partial D=\emptyset$,
the tangle $(\Gamma\cap D,D)$ 
is an NS-tangle. 
This contradicts Lemma~\ref{LemNS-Tangle}.

{\bf Case 2}. If $Cl(U)$ is an annulus,  
then by applying DH-tricks along 
a simple arc $\ell$ connecting 
a point $p$ in $C$ and 
a point in $\partial N-\partial D^2$ 
with $\ell \cap C=p$,
we can assume that 
the hoop $C$ bounds a disk 
with at least one white vertex and 
with at most one crossing.
Thus we can find an NS-tangle 
by the same way as the one of Case 1.
This contradicts Lemma~\ref{LemNS-Tangle}.
\end{proof}

\begin{lemma}\LABEL{Lemma3-2}
Let $\Gamma$ be a minimal chart 
in a disk $D^2$.
If $2\le c(\Gamma)\le3$,
then there exists no ring.
\end{lemma}

\begin{proof}
Let $N$ be a brim of $D^2$.
Suppose that there exists a ring $C$.
Let $A$ be 
a regular neighbourhood of $C$ 
in $D^2$.
Since there exist 
at most three crossings, 
one of a complementary domain of $A$ 
in $D^2$ 
contains at most one crossing.
Let $D$ be the closure of 
the complementary domain of $A$.
If necessary, 
applying DH-tricks along 
a simple arc $\ell$ 
connecting a point on $C$ and 
a point in $\partial N-\partial D^2$ 
with $\ell\cap C$ being a point,
we can assume that $D$ is a disk.
Since each complementary domain of $C$ 
in $D^2$ 
contains a white vertex 
by Proposition~\ref{Assumption0}(c), 
and 
since there are 
at most three crossings on $C$,
the tangle $(\Gamma\cap D,D)$ is an NS-tangle.
This contradicts Lemma~\ref{LemNS-Tangle}.
\end{proof}

Let $\Gamma$ be a minimal chart.
For a subset $X$ of $\Gamma$,
let
\begin{enumerate}
\item[]
$B(X)=$ the union of all the disk 
bounded by a cycle in $X$, and
\item[]
$T(X)=$ the union of all the terminal edge 
intersecting $X\cup B(X)$.
\end{enumerate}
The set $X\cup B(X)\cup T(X)$ is called
the SC-{\it closure} of $X$ and\index{SC-closure $SC(X)$}
denoted by $SC(X)$. 
Each connected component of  $SC(X)$ 
is simply connected.

Let $\Gamma$ be a chart, and 
$m$ a label of the chart. 
Let $\mathcal W$ be the set of 
all the white vertices of $\Gamma$. 
The closure of a connected component 
of $\Gamma_m-\mathcal W$ 
is called an {\it internal} edge of label $m$\index{internal edge} 
if it contains a white vertex 
but does not contain any black vertex, 
here we consider $\Gamma_m$ as a topological set.

Let $G$ be a subgraph of a chart $\Gamma$. 
An internal edge $e$ in $G$ is called 
a {\it cut-edge} for $G$\index{cut-edge}
if $Cl(G-e)$ is not connected.

\begin{lemma} 
\LABEL{NoCutEdge}
Let $\Gamma$ be a minimal chart 
in a disk $D^2$. 
Let $\alpha=\alpha(\Gamma)$, and 
$\beta=\beta(\Gamma)$. 
If $2\le c(\Gamma)\le 3$, 
then there does not exist 
an internal cut-edge for
$\Gamma_\alpha$ 
nor $\Gamma_\beta$.
\end{lemma}

\begin{proof}
Let $N$ be a brim of $D^2$.
Suppose that there exists 
an internal cut-edge $\overline{e}^*$
for $\Gamma_\alpha$. 
If necessary, 
applying DH-tricks along a simple arc 
connecting a point near 
$\overline{e}^*$ 
but not in a bigon, 
and 
a point in $\partial N-\partial D^2$, 
we can assume that 
\begin{enumerate}
\item[(1)] 
there exists a simple arc $L$ connecting 
a point $p$ in Int~${\overline e}^*$ and 
a point in the brim 
with $L\cap\Gamma_\alpha=p$ 
(see Fig.~\ref{fig10}(a)).
\end{enumerate}
Let $X_1,X_2$ be the 
connected components of 
$Cl(\Gamma_\alpha-{\overline e}^*)$ 
such that
\begin{enumerate}
\item[(2)] 
each of $X_1\cap {\overline e}^*$ and 
$X_2\cap {\overline e}^*$ 
consists of exactly one point.
\end{enumerate}
Further, the existence of 
the arc $L$ of Statement (1) 
assures that  
the SC-closures of $X_1,X_2$ 
do not intersect each other, 
i.e. 
$SC(X_1)\cap SC(X_2)=\emptyset$.
Furthermore, 
$c(\Gamma)\le 3$ implies that 
\begin{enumerate}
\item[(3)] 
one of $SC(X_1)$ and $SC(X_2)$ contains 
at most one crossing, say $SC(X_1)$.
\end{enumerate}
If $X_1$ does not contain a crossing, 
let $E$ be a regular neighbourhood of $SC(X_1)$ 
in $D^2$. 
Since $SC(X_1)$ contains at most one crossing,
the disk $E$ contains at most one crossing.
Hence $(\Gamma\cap E,E)$ 
is an NS-tangle of label $\alpha+1$.
This contradicts Lemma~\ref{LemNS-Tangle}.

Suppose that $X_1$ contains a crossing $v$ 
in $\Gamma_\alpha\cap\Gamma_k$ 
for some label $k$ with $1<|\alpha-k|$. 
Then there exists 
an internal edge ${\overline e}$ of 
label $\alpha$ containing the crossing $v$. 

If ${\overline e}$ is not an internal cut-edge 
for $X_1$, 
let $E$ be a regular neighbourhood of $SC(X_1)$ 
in $D^2$ 
(see Fig.~\ref{fig10}(b)). 
Then $E$ contains 
a white vertex of ${\overline e}$.
Hence by Statement $(3)$, 
the tangle $(\Gamma\cap E,E)$ 
is an NS-tangle of label $\alpha+1$.
This contradicts Lemma~\ref{LemNS-Tangle}.

If ${\overline e}$ is an internal cut-edge 
for $X_1$, 
then  $Cl(X_1-\overline e)$ consists of 
two connected components.
By Statement (2), 
one of the connected components
does not intersect the edge ${\overline e}^*$, 
say $X$ 
(see Fig.~\ref{fig10}(c)). 
Then $SC(X)$ does not contain $\overline e$. 
Let $E$ be a regular neighbourhood of 
$SC(X)$ in $D^2$. 
Now $X\subset X_1$ implies 
$SC(X)\subset SC(X_1)$. 
Thus $E$ does not contain a crossing 
by Statement $(3)$.
Further $\Gamma_\alpha\cap\partial E$ 
consists of one point.
Hence the tangle $(\Gamma\cap E,E)$ 
is an NS-tangle of label $\alpha+1$.
This contradicts Lemma~\ref{LemNS-Tangle}.
Thus there does not exist 
an internal cut-edge for  
$\Gamma_\alpha$.

Similarly we can show that 
there does not exist 
an internal cut-edge for 
$\Gamma_\beta$.
Thus Lemma~\ref{NoCutEdge} holds.
\end{proof}

\begin{figure}
\begin{center}
\includegraphics{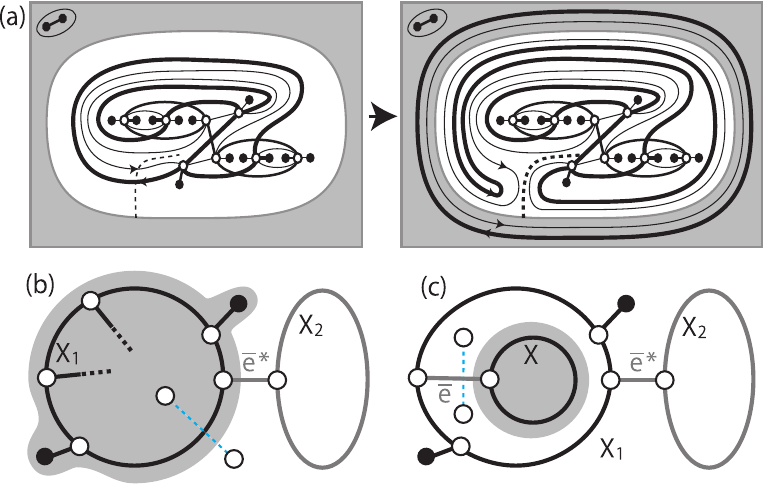}
\caption{\LABEL{fig10} (a) Thick edges are of label $\alpha$. (b) $X_1$ is thick black, $\overline e^*$ is a cut-edge. (c) $X$ is thick black, edges $\overline e^*,\overline e$ are cut-edges. }
\end{center}
\end{figure}

Let $\Gamma$ be a chart.
For each label $m$,
we define 
${\rm Main}(\Gamma_m)=
\Gamma_m\cap {\rm Main}(\Gamma)$.

\begin{lemma} 
\LABEL{Gamma1Connected}
Let $\Gamma$ be a minimal chart 
in a disk $D^2$. 
Let $\alpha=\alpha(\Gamma)$, and 
$\beta=\beta(\Gamma)$. 
If $2\le c(\Gamma)\le 3$, 
then ${\rm Main}(\Gamma_\alpha)$ and 
${\rm Main}(\Gamma_\beta)$ 
are connected.
\end{lemma}

\begin{proof}
If Main$(\Gamma_\alpha)$ 
is not connected, 
by assuming there exists 
an 'imaginary' cut-edge 
connecting two components of $\Gamma_\alpha$, 
we get a similar contradiction as the one of 
Lemma~\ref{NoCutEdge}. 
Thus ${\rm Main}(\Gamma_\alpha)$ is connected. 
Similarly we can show that 
${\rm Main}(\Gamma_\beta)$ is connected.
Thus Lemma~\ref{Gamma1Connected} holds.
\end{proof}

\begin{lemma} 
\LABEL{ContainCrossing}
Let $\Gamma$ be a minimal chart 
in a disk $D^2$. 
Let $\alpha=\alpha(\Gamma)$, and 
$\beta=\beta(\Gamma)$. 
If $2\le c(\Gamma)\le 3$, 
then each of $\Gamma_\alpha$ 
and $\Gamma_\beta$ 
contains a crossing.
\end{lemma}

\begin{proof}
Let $N$ be a brim of $D^2$.
Suppose that $\Gamma_\alpha$ does not 
contain a crossing.
There are two cases.\\
Case 1. There exists a complementary domain $U$ 
of ${\rm Main}(\Gamma_\alpha)$ 
in $D^2$ 
containing exactly one crossing. \\
Case 2. There exists a complementary domain $U$ of 
${\rm Main}(\Gamma_\alpha)$ 
in $D^2$ containing all the crossings.

{\bf Case 1}. 
Suppose that 
there exists a complementary domain $U$ 
of ${\rm Main}(\Gamma_\alpha)$ 
in $D^2$ 
containing exactly one crossing.  
If necessary, 
applying DH-tricks along a simple arc  
connecting a point in 
Int $SC({\rm Main}(\Gamma_\alpha))$  
but not in a bigon, 
and a point in $\partial N-\partial D^2$, 
we can assume that 
\begin{enumerate}
\item[(1)] 
$U$ does not intersect the brim.
\end{enumerate}
Let $A$ be a regular neighbourhood of 
$\partial U$ in $D^2$, and $E=Cl(U-A)$. 
Then Statement $(1)$ and 
Lemma~\ref{Gamma1Connected} 
imply that 
$E$ is a disk. 
Since $U$ contains 
exactly one crossing, so does $E$.
Hence $(\Gamma\cap E,E)$ 
is an NS-tangle of label $\alpha+1$. 
This contradicts Lemma~\ref{LemNS-Tangle}.

{\bf Case 2}. 
Suppose that 
there exists a complementary domain $U$ of 
${\rm Main}(\Gamma_\alpha)$ 
in $D^2$ containing 
all the crossings. 
Then 
\begin{enumerate}
\item[(2)] 
$\Gamma_\alpha$ does not contain any crossing.
\end{enumerate} 
If necessary, 
applying DH-tricks along a simple arc  
connecting 
a point in $U$ but not in a bigon, 
and a point in $\partial N-\partial D^2$, 
we can assume that 
$U$ intersects the brim. 
Then 
$U\cap B({\rm Main}(\Gamma_\alpha))
=\emptyset$, 
here 
$B({\rm Main}(\Gamma_\alpha))$ 
is the union of all the disk bounded by 
a cycle in ${\rm Main}(\Gamma_\alpha)$. 
Since any terminal edge does not 
contain a crossing,
the SC-closure 
$SC({\rm Main}(\Gamma_\alpha))
={\rm Main}(\Gamma_\alpha)\cup
B({\rm Main}(\Gamma_\alpha))\cup
T({\rm Main}(\Gamma_\alpha))$ 
does not contain any crossing. 
Let $E$ be a regular neighbourhood of 
$SC({\rm Main}(\Gamma_\alpha))$ in $D^2$. 
Since $\Gamma_\alpha$ contains a white vertex, 
so does $E$.
Hence $(\Gamma\cap E,E)$
is an NS-tangle of label $\alpha+1$ 
without crossing. 
This contradicts Lemma~\ref{LemNS-Tangle}.
Hence $\Gamma_\alpha$ contains a crossing.

Similarly we can show that 
$\Gamma_\beta$ contains a crossing.
Thus Lemma~\ref{ContainCrossing} holds.
\end{proof}

\begin{lemma} 
\LABEL{ContainTwoCrossing}
Let $\Gamma$ be a minimal chart 
in a disk $D^2$. 
Let $\alpha=\alpha(\Gamma)$, and 
$\beta=\beta(\Gamma)$. 
If $2\le c(\Gamma)\le 3$, 
then there exist cycles $C_\alpha$ 
and $C_\beta$ with 
$C_\alpha\subset\Gamma_\alpha,
 C_\beta\subset\Gamma_\beta$ 
each of which contains two crossings.
\end{lemma}

\begin{proof}
Let $N$ be a brim of $D^2$.
By Lemma~\ref{ContainCrossing}, 
$\Gamma_\alpha$ contains at least one crossing $v$.
By Lemma~\ref{NoCutEdge},
there exist two complementary domains 
$U_1,U_2$ of ${\rm Main}(\Gamma_\alpha)$
in $D^2$ 
with 
\begin{enumerate}
\item[(1)]
$v\in\partial U_1\cap\partial U_2$.
\end{enumerate}
Then $c(\Gamma)\le 3$ implies that 
one of $U_1,U_2$ contains at most one crossing,
say $U_1$. 
If necessary, 
applying DH-tricks along a simple arc $\ell$ 
connecting a point $p$ in 
$U_2$ but not in a bigon,  
and a point in $\partial N-\partial D^2$  
with $\ell\cap \partial U_2$ a point, 
we can assume that
$U_1$ does not intersect the brim. 
Now Lemma~\ref{NoCutEdge} and 
Lemma~\ref{Gamma1Connected}
imply that  
$Cl(U_1)$ is a disk. 
Let $C_\alpha=\partial(Cl(U_1))$.
Then 
$C_\alpha$ is a simple closed curve.
Now $C_\alpha$ contains the crossing $v$ by (1).

Suppose that $C_\alpha$ 
contains exactly one crossing $v$.
Let 
$A$ be a regular neighbourhood of 
$\partial U$ in $D^2$, and 
$E=Cl(U_1-A)$.
Then $E$ is a disk. 
Now $v\in C_\alpha=\partial(Cl(U_1))$ 
implies that
the disk $E$ contains a white vertex. 
Thus $(\Gamma\cap E,E)$ is 
an NS-tangle of label $\alpha+1$.
This contradicts Lemma~\ref{LemNS-Tangle}.
Hence $C_\alpha$ contains 
at least two crossings.

Similarly we can show that 
there exists a cycle 
$C_\beta\subset\Gamma_\beta$ 
contains at least two crossings.
Thus Lemma~\ref{ContainTwoCrossing} holds.
\end{proof}

{\bf Proof of Lemma~\ref{StIILemma1}}.
Let $\Gamma$ be a 
$2$-crossing minimal chart 
with two crossings $x_1,x_2$.
By Lemma~\ref{ContainTwoCrossing},
there exist cycles $C_\alpha$ 
and $C_\beta$ with 
$C_\alpha\subset\Gamma_\alpha,
 C_\beta\subset\Gamma_\beta$ 
and $C_\alpha\cap C_\beta=\{ x_1,x_2 \}$.
Thus Statement (a) holds.

Let $N_1,N_2$ be regular neighbourhoods of 
$x_1,x_2$ in $D^2$ respectively.
Let $N$ be a brim containing 
all the free edges and simple hoops.
If necessary, 
applying DH-tricks along a simple arc $\ell$  
connecting 
a point $p$ in $\partial N_1-{\rm Main}(\Gamma)$ 
and a point in $\partial N-\partial D^2$ 
with $\ell\cap N_1=p$,
we can assume 
that there exists a simple arc $\gamma$ 
connecting 
a point $p^*$ in $\partial N_1$ and 
a point $q^*$ in $\partial N$
with 
$\gamma\cap({\rm Main}(\Gamma)\cup N_2)
=\emptyset$, 
$\gamma\cap N_1=p^*$, and 
$\gamma\cap N=q^*$. 
Let $N(\gamma)$ be 
a regular neighbourhood 
of $\gamma$ in $D^2$, and 
$A=Cl(D^2-(N_1\cup N_2\cup N(\gamma)\cup N))$ 
(see Fig.~\ref{fig11}).
Then $A$ is an annulus
and 
$\Gamma\cap \partial A
=\Gamma\cap(\partial N_1\cup\partial N_2)
=(C_\alpha\cup C_\beta)\cap \partial A$. 
Thus 
$\Gamma\cap \partial A$ consists of eight points.
Hence Statement (d) holds.

Since $N_1\cup N_2\cup N(\gamma)\cup N$
does not contain any white vertex,
the annulus $A$ contains all the white vertices 
of $\Gamma$.
By Lemma~\ref{Lemma3-1}, 
any hoop is simple. 
Since any simple hoops and free edges are in the brim $N$, 
the annulus $A$ does not intersect 
hoops nor free edges.
Hence Statement (b) holds.

Now $N_2$ and $N_1\cup N(\gamma)\cup N$
are the connected components of $Cl(D^2-A)$ 
each of which contains a crossing.
Thus Statement (c) holds.

Since $C_\beta\cap \Gamma_\alpha$
consists of the two crossings,
the set $\Gamma_\alpha\cap A$ 
consists of 
two connected components 
separated by $C_\beta$ 
(see Fig.~\ref{fig11}).
Similarly
$\Gamma_\beta\cap A $ consists of 
two connected components 
separated by $C_\alpha$.
Thus Statement (e) and Statement (f) hold.
Hence Lemma~\ref{StIILemma1} holds. 
{\hfill {$\square$}\vspace{1.5em}}

\begin{figure}
\begin{center}
\includegraphics{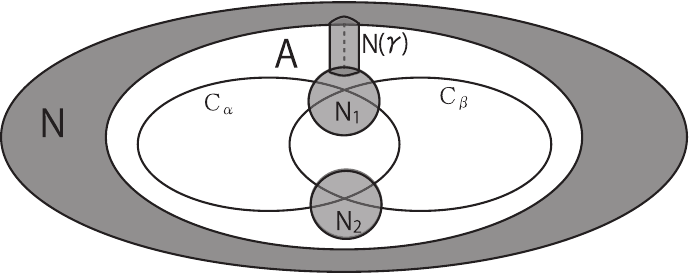}
\caption{\LABEL{fig11} The annulus $A$ is the white area. }
\end{center}
\end{figure}



\section{Proofs of Theorem~1.2 and Theorem~1.3}
\LABEL{s:Th2Th3}

Let $\Gamma$ be a chart.
A tangle $(\Gamma\cap D,D)$ is said to be
{\it admissible}\index{admissible tangle}
provided that
\begin{enumerate}
\item[(i)] 
any edge intersecting $\partial D$
is contained in an internal edge, 
\item[(ii)] 
if an internal edge $\overline e$  
intersects $\partial D$,
then each connected component of 
$\overline e\cap D$ contains a white vertex.
\end{enumerate}

\begin{lemma}
{\rm(\cite[Theorem 1.3]{MinimalChart})}
\LABEL{TwoColorGateTangle} 
If $(\Gamma\cap D,D)$ is 
an admissible tangle 
in a minimal chart $\Gamma$ such that 
\begin{enumerate}
\item[{\rm $($a$)$}] 
$\Gamma\cap D\subset \Gamma_m\cup\Gamma_{m-1}$ 
or 
$\Gamma\cap D\subset \Gamma_m\cup\Gamma_{m+1}$ 
for some label $m$,
\item[{\rm $($b$)$}] 
$\Gamma_m\cap \partial D$ 
consists of exactly two points, 
and 
\item[{\rm $($c$)$}] 
$\Gamma_m\cap D$ contains a cycle, 
\end{enumerate}
then
the tangle $(\Gamma\cap D,D)$ is 
a simple IO-tangle of label $m$.
\end{lemma}

Let $m$ be a label 
of a chart $\Gamma$, 
and $(\Gamma\cap D,D)$ 
a Type-II elementary IO-tangle 
of label $m$ 
with a boundary IO-arc pair $(L_I,L_O)$. 
Let $L_D$ be 
a simple proper arc of $D$ 
contained in $\Gamma_m\cap D$ 
and containing 
all the white vertices 
$v_1,v_2,\cdots,v_t$ 
of $\Gamma\cap D$, 
situated on $L_D$ in this order.
Let $\Delta_I,\Delta_O$ 
be the closures of 
connected components 
of $D-L_D$ 
with $L_I\subset \Delta_I,L_O\subset \Delta_O$. 
Then the tangle 
$(\Gamma\cap D,D)$ is
said to be 
of {\it Type-II$_t$ } provided that\index{Type-II$_t$}
(see Fig.~\ref{fig12})
\begin{enumerate}
\item[(i)] 
for each $i=1,2,\cdots,t$, 
there exists a terminal edge $\tau_i$ 
of label $m$ 
containing $v_i$ and 
\item[(ii)] 
for each $i=1,2,\cdots,t$, 
the terminal edge $\tau_i$ 
is contained in 
$\Delta_I$ or $\Delta_O$ alternately.
\end{enumerate}

\begin{figure}
\centerline{\includegraphics{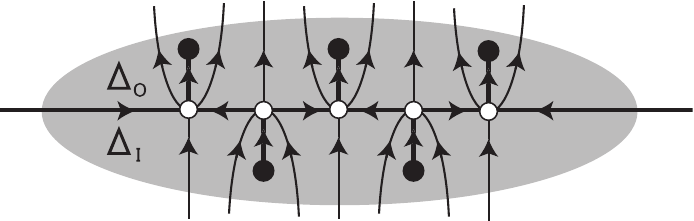}}
\caption{\LABEL{fig12}
A Type-II$_5$ elementary IO-tangle of label $m$. The thick lines are of label $m$.}
\end{figure}

\begin{lemma}\LABEL{lemEQ}
Let $\Gamma$ be a minimal chart, and 
$(\Gamma\cap D,D)$ a non-trivial admissible tangle
such that
$\Gamma\cap D\subset\Gamma_m\cup\Gamma_k$ 
for some labels $m,k$ with $|m-k|=1$.
Then the following conditions are equivalent.
\begin{enumerate}
\item[{\rm (a)}] 
The tangle 
is a Type-II elementary IO-tangle 
of label $m$.
\item[{\rm (b)}]  
$\Gamma_m\cap D$ contains no cycle, 
$\Gamma_m\cap \partial D$ 
consists of two points, 
and 
there exists 
a proper arc $L_D$ of $D$ in 
$\Gamma_m\cap D$ 
containing all the white vertices 
in $D$.
\item[{\rm (c)}] 
The tangle is 
a Type-II$_t$ elementary IO-tangle 
of label $m$ 
for some integer $t$.
\end{enumerate}
\end{lemma}

\begin{proof}
Clearly 
Statement (a) implies Statement (b).

Suppose that 
Statement (b) holds.
Let $v_1,v_2,\cdots,v_t$ be 
all the white vertices in $L_D$
situated on $L_D$ in this order.

Now for each $i=1,2,\cdots,t$, 
the white vertex $v_i$ is contained in 
a terminal edge $\tau_i$
of label $m$.
For, if not, let $e$ be 
the edge of label $m$ 
with $e\cap L_D\ni v_i$ and $e\not\subset L_D$ 
but not a terminal edge. 
Since $\Gamma_m\cap \partial D$ 
consists of two points
$\partial L_D$, 
we have $e\cap\partial D=\emptyset$. 
Since the arc $L_D$ contains 
all the white vertices in $D$, 
the union $e\cup L_D$ contains 
a cycle in $\Gamma_m\cap D$. 
This is a contradiction.
Thus each white vertex $v_i$ is contained in 
a terminal edge $\tau_i$.

Now $D-L_D$ consists of 
two connected components.
If $\tau_1$ is inward 
(resp. outward) at $v_1$,
then let $\Delta_I$ (resp. $\Delta_O$) 
be 
the closure of the one of 
the two connected components 
containing $\tau_1$, and
$\Delta_O$ (resp. $\Delta_I$) 
the closure of 
the other connected component.

If $t=1$, then 
the tangle is clearly a Type-II$_1$ 
elementary IO-tangle. 
Hence we assume $t\ge 2$. 

Suppose that 
for some $i\in\{1,2,\cdots,t-1\}$,
two terminal edges 
$\tau_i$ and $\tau_{i+1}$ are 
contained in $\Delta_I$ or $\Delta_O$
simultaneously.
Then we can eliminate 
the two white vertices $v_i,v_{i+1}$
by two C-I-M2 moves and a C-I-M3 move
(see Fig.~\ref{fig13}).
This contradicts that $\Gamma$ is 
a minimal chart.
Hence the terminal edges 
$\tau_1,\tau_2,\cdots,\tau_t$
alternately belong to 
$\Delta_I$ and $\Delta_O$.

Thus each terminal edge of label $m$ in $\Delta_I$ 
is inward at a white vertex of $L_D$, and
each terminal edge of label $m$ in $\Delta_O$
is outward at a white vertex of $L_D$.
Hence if an edge intersects Int~$\Delta_I$ 
(resp. Int~$\Delta_O$),
then the edge contains an inward (resp. outward)
arc in $\Delta_I$ (resp. $\Delta_O$) 
at a white vertex in $L_D$.
Let $L_I=Cl(\partial \Delta_I - L_D)$
and $L_O=Cl(\partial \Delta_O - L_D)$.
Then $L_I,L_O$ satisfy 
Condition (iii) and (iv) for
an IO-tangle.
Thus the tangle $(\Gamma\cap D,D)$ 
is 
a Type-II$_t$ elementary IO-tangle of
label $m$. 
Hence Statement (c) holds.

It is also clear that 
Statement (c) implies Statement (a).
This proves Lemma~\ref{lemEQ}.
\end{proof}

\begin{figure}
\centerline{\includegraphics{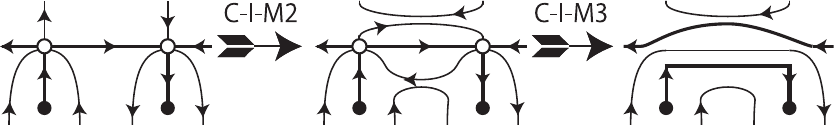}}
\caption{The thick lines are of label $m$.
\LABEL{fig13}}
\end{figure}

\begin{lemma}\LABEL{ArcConnectTwo}
Let $\Gamma$ be a minimal chart, and 
$m,k$ labels of $\Gamma$ with 
$|m-k|=1$. 
Let $(\Gamma\cap D,D)$ be a tangle 
with 
$\Gamma\cap D\subset \Gamma_m\cup\Gamma_k$
and 
$|\Gamma_m\cap\partial D|=2$. 
If no terminal edge nor free edge 
intersects $\partial D$, 
then $\Gamma_m\cap D$ contains 
a proper arc connecting 
the two points in 
$\Gamma_m\cap\partial D$.
\end{lemma}

\begin{proof}
If not, let $X$ be 
a connected component
of $\Gamma_m\cap D$ containing 
exactly one point in 
$\Gamma_m\cap\partial D$.
Then $X$ contains a white vertex, 
because neither terminal edge 
nor free edge 
intersects $\partial D$. 
Let $D^*$ be 
a regular neighbourhood of 
the SC-closure $SC(X)$ in $D$.
Then $(\Gamma\cap D^*,D^*)$
is an NS-tangle of label $k$. 
This contradicts 
Lemma~\ref{LemNS-Tangle}.
Thus
$\Gamma_m\cap D$ contains 
a proper arc 
connecting the two points 
$\Gamma_m\cap \partial D$.
\end{proof}

Let $\Gamma$ be a chart, 
$v$ a white vertex, 
and 
$e$ a terminal edge inward 
(resp. outward) at $v$.
Then the two edges inward 
(resp. outward) at $v$ 
different from $e$\index{sibling edge} 
are called the {\it sibling edges} of $e$. 

\begin{lemma}
\LABEL{LemSimpleIO-Tangle1}
If $(\Gamma\cap D,D)$ 
is an admissible tangle 
in a minimal chart $\Gamma$
such that
\begin{enumerate}
\item[{\rm $($a$)$}] 
$\Gamma\cap D\subset \Gamma_m\cup\Gamma_k$ 
for some labels $m,k$ with $|m-k|=1$,
\item[{\rm $($b$)$}] 
$\Gamma_m\cap \partial D$ consists of 
exactly two points, 
and
\item[{\rm (c)}]
the disk 
$D$ contains at least two white vertices, 
\end{enumerate}
then 
$(\Gamma\cap D,D)$ 
is a simple IO-tangle 
of label $m$. 
Further, if $D$ does not contain 
a cycle of label $m$, then
the tangle is a Type-II elementary IO-tangle 
of label $m$.
\end{lemma}

\begin{proof}
Neither terminal edge nor free edge is 
an internal edge.
By Condition (i) for an admissible tangle,
the boundary $\partial D$ intersects 
neither terminal edge nor free edge.
Thus the tangle $(\Gamma\cap D, D)$ 
satisfies Condition (i) for an IO-tangle.
Condition (a) is the same one of
Condition (ii) for an IO-tangle.

If $\Gamma_m\cap D$ 
contains a cycle,
then the tangle $(\Gamma\cap D,D)$ 
is a simple IO-tangle 
of label $m$ by 
Lemma~\ref{TwoColorGateTangle}.

Suppose $\Gamma_m\cap D$ contains no cycle.
By Lemma~\ref{ArcConnectTwo}
there exists  
a proper arc $L$ in $\Gamma_m\cap D$ 
connecting 
the two points $\Gamma_m\cap\partial D$.

{\bf Claim.}
The arc $L$ contains 
all of the white vertices in $D$.

{\bf Proof of Claim}. 
Suppose that there exists 
a white vertex in $D-L$. 
Let $X$ be the closure of 
the connected component of 
$\Gamma_m\cap (D-L)$ containing the vertex. 

If $X\cap L=\emptyset$, 
then $X$ does not intersect $\partial D$, 
because $L$ contains the two points 
$\Gamma_m\cap \partial D$. 
Let $D'$ be a regular neighbourhood of
the SC-closure $SC(X)$ in $D$. 
Then the tangle $(\Gamma\cap D',D')$ 
is an NS-tangle of label $k$. 
This contradicts Lemma~\ref{LemNS-Tangle}.

If $X\cap L\neq\emptyset$, 
then $X\cap L$ consists of 
one white vertex $v$, 
because there does not exist 
a cycle of label $m$. 
Let $e$ be the edge in $X$ containing 
the vertex $v$. 
Let $X'=Cl(X-e)$, and 
$D'$ a regular neighbourhood of
the SC-closure $SC(X')$ in $D$. 
Then the tangle $(\Gamma\cap D',D')$ 
is an NS-tangle of label $k$. 
This contradicts Lemma~\ref{LemNS-Tangle}.
Hence Claim holds.

Hence the tangle $(\Gamma\cap D,D)$ satisfies 
(b) of Lemma~\ref{lemEQ}.
Thus by Lemma~\ref{lemEQ}, 
the tangle $(\Gamma\cap D,D)$ 
is 
a Type-II elementary IO-tangle of
label $m$.

To show the tangle is simple, 
suppose that there exists a terminal edge $e$ of label $k$ in $D$.
Let $e^*,e^{**}$ be the sibling edges of $e$.
If one of  $e^*,e^{**}$ is contained 
in the disk $D$,
then we can find 
a non-admissible consecutive triplet 
which contradicts
 Consecutive Triplet Lemma
(Lemma~\ref{ConsecutiveTripletLemma}).
Thus both of $e^*$ and $e^{**}$ 
intersect $\partial D$. 
Hence $L=(e^*\cup e^{**})\cap D$ 
contains all the white vertices in $D$ 
by Claim. 
But $L$ contains exactly one white vertex.
This contradicts Condition (c).
Thus there does not exist 
any terminal edge of label $k$ in $D$.
Therefore the tangle is a simple IO-tangle.
\end{proof}

\begin{remark}\LABEL{remWhite1}
{\rm 
If $(\Gamma\cap D,D)$ 
is an admissible tangle 
in a minimal chart $\Gamma$
such that
\begin{enumerate}
\item[{\rm $($a$)$}] 
$\Gamma\cap D\subset \Gamma_m\cup\Gamma_k$ 
for some labels $m,k$ with $|m-k|=1$,
\item[{\rm $($b$)$}] 
$\Gamma_m\cap \partial D$ consists of 
exactly two points, 
and
\item[(c)]
the disk 
$D$ contains 
exactly one white vertex,
\end{enumerate}
then 
$(\Gamma\cap D,D)$ 
is one of the tangles shown in 
Fig.~\ref{fig14}. 
The tangles shown in 
Fig.~\ref{fig14}(a),(b) 
are simple IO-tangles. 
The tangles shown in 
Fig.~\ref{fig14}(c),(d) 
are not simple IO-tangles, but 
IO-tangles and N-tangles 
(see Section~\ref{s:Indices} 
for the precise definition of N-tangles).
}
\end{remark}

\begin{figure}
\centerline{\includegraphics{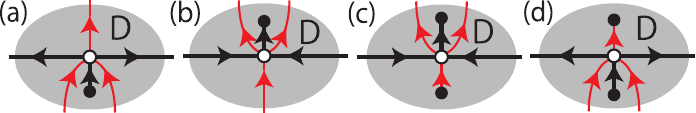}}
\caption{Admissible tangles with exactly one white vertex.\LABEL{fig14}}
\end{figure}

{\bf Proof of Theorem~\ref{StIITheorem2}.} 
We use all the notation in the definition of 
fundamental tangles of 
a 2-crossing chart 
mentioned just before 
Theorem~\ref{StIITheorem2}.
Let $(\Gamma\cap D^*_i,D^*_i)$ be 
a fundamental tangle. 
If $i=1,3$ then 
let 
$m=\alpha(\Gamma),k=m+1$ 
otherwise 
$m=\beta(\Gamma),k=m-1$.

Since $D^*_i$ is 
a regular neighbourhood 
of $SC(X_i)$ in the annulus $A$ 
with $X_i\cap \partial A$ two points, 
the intersection 
$\Gamma_m\cap\partial D^*_i$ 
consists of two points.
Further by Boundary Condition Lemma
(Lemma~\ref{BoundaryConditionLemma}) 
we have
$\Gamma\cap D^*_i\subset\Gamma_m\cup\Gamma_k$.
Since $D^*_i$ contains 
at least two white vertices,
the fundamental tangle is 
a simple IO-tangle
by Lemma~\ref{LemSimpleIO-Tangle1}. 
{\hfill {$\square$}\vspace{1.5em}}

\begin{lemma}
\LABEL{2ColorBeads}
{\rm(\cite[Lemma 10.1]{MinimalChart})}
Let $\Gamma$ be a minimal chart, and
$m,k$ labels of $\Gamma$ 
with $|m-k|=1$.
Let $(\Gamma\cap D,D)$ be
an admissible tangle with 
$\Gamma\cap D\subset\Gamma_m\cup\Gamma_k$
and 
$|\Gamma_m\cap\partial D|=2$.
If $\Gamma_m\cap D$ contains a cycle, 
then 
there exist disks $E_1,E_2,\cdots, E_d$ in {\rm Int}~$D$ and simple arcs $L_0,L_1,\cdots,L_d$ in 
$Cl(D-\cup^{d}_{i=1}E_i)$
 such that
\begin{enumerate}
\item[{\rm (a)}] 
$\partial E_i\subset \Gamma_m$ for each $i=1,2,\cdots,d$ and $L_j\subset \Gamma_m$ for each $j=0,1,\cdots,d$,
\item[{\rm (b)}] 
for each $j=1,\cdots,d-1$, $L_j$ connects 
a white vertex in $\partial E_j$ and 
a white vertex in $\partial E_{j+1}$, \\
the arc $L_0$ connects a point in $\partial D$ 
and a white vertex in $\partial E_{1}$, and\\
the arc $L_d$ connects a white vertex 
in $\partial E_{d}$ 
and a point in $\partial D$,
\item[{\rm (c)}]
if an edge of label $m$ intersects 
$D-((\cup_{i=1}^{d} E_i)\cup (\cup_{j=0}^{d} L_j))$, 
then it is a terminal edge.
\end{enumerate}
\end{lemma}

{\bf Proof of Theorem~\ref{StIITheorem3}.}
Let $k$ be a label of $\Gamma$ 
with $\Gamma\cap D\subset\Gamma_m\cup\Gamma_k$.
By Condition (iii) for an IO-tangle,
we have $|\Gamma_m\cap \partial D|=2$. 
Now any non-trivial simple IO-tangle 
is admissible. 

Suppose that $\Gamma_m\cap D$ 
contains a cycle.
By Lemma~\ref{2ColorBeads},
there exist disks $E_1,E_2,\cdots, E_d$ in {\rm Int}~$D$ and simple arcs $L_0,L_1,\cdots,L_d$ in 
$Cl(D-\cup^{d}_{i=1}E_i)$
satisfying (a), (b) and (c) in 
Lemma~\ref{2ColorBeads}.
Let $X=(\cup^{d}_{i=1}E_i)\cup(\cup^{d}_{j=0}L_j)$. Since $X$ 
is connected and since 
$|\partial D\cap X|=2$, 
for the SC-closure $SC(X)$, 
$D-SC(X)$ consists of two components. 
Let $D_I,D_O$ be 
the closures of the two components. 

{\bf Claim.} Any edge of label $k$ 
intersecting $D_I$ or $D_O$ 
intersects $\partial D$.

{\it Proof of Claim}.
Suppose there exists an edge $e$ 
of label $k$ in $D_I$. 
Since there is no terminal edge of label $k$ 
by the condition for simple IO-tangle, 
the edge $e$ is oriented 
from a white vertex $w_1$ 
to a white vertex $w_2$.
Let $D'$ be a regular neighbourhood of $SC(X)$ 
in $D$. Then $(\Gamma\cap D',D')$ is 
an IO-tangle by 
Lemma~\ref{TwoColorGateTangle}.
Let $D'_I,D'_O$ be the closures of 
the connected components of $D'-SC(X)$ 
with 
$D'_I\subset D_I,D'_O\subset D_O$. 
Then $Cl(e\cap D'_I)$ consists of two arcs 
$e'_1,e'_2$ here $w_1\in e'_1,w_2\in e'_2$.  

Now $e'_1,e'_2\subset D'_I$ and 
$e'_1,e'_2$ intersect $\partial D'$. 
But $e'_1$ is outward at $w_1$ 
and $e'_2$ is inward at $w_2$. 
This contradicts the fact that 
$(\Gamma\cap D',D')$ is 
an IO-tangle. 
Thus Claim holds.

Now $D-X$ does not contain any white vertices 
by (c) of Lemma~\ref{2ColorBeads}.
Since there is no terminal edge 
of label $k$
by the condition for a simple IO-tangle,
for each $i=1,2,\cdots,d$ 
Claim implies that 
there exist 
two proper arcs $\ell_{2i-1}$ and $\ell_{2i}$ 
of $D$ of label $k$ 
(see Fig.~\ref{fig15}(a)) such that

$\ell_{2i-1}=$ the proper arc of $D$ 
contained in $\Gamma_k$ with 
$\ell_{2i-1}\cap\Gamma_m=
L_{i-1}\cap E_i$, and

$\ell_{2i}=$ the proper arc of $D$ 
contained in $\Gamma_k$ with 
$\ell_{2i}\cap\Gamma_m=E_i\cap L_i$.\\
Then proper arcs 
$\ell_{2i-1}$ and $\ell_{2i}$
split the disk $D$ into three disks.
Let $D_{2i}^*$ be the one of the three disks
containing $E_i$
(see Fig.~\ref{fig15}(a)).
Let $D_{2i}$ be 
a regular neighbourhood of $D_{2i}^*$ 
in $D$ 
(see Fig.~\ref{fig15}(b)). 
Now $Cl(D-\cup_{i=1}^dD_{2i})$
consists of $d+1$ disks.
For each $i=0,1,2,\cdots,d$,
let $D_{2i+1}$ be the one of the disks
intersecting the arc $L_i$ 
(see Fig.~\ref{fig15}(b)).
Then for each $i=2,4,\cdots,2d$,
the tangle $(\Gamma\cap D_i,D_i)$ is 
a Type-I elementary IO-tangle of label $m$.
Further for each $i=1,3,\cdots,2d+1$,
the tangle $(\Gamma\cap D_i,D_i)$ is 
a Type-II elementary IO-tangle of label $m$.
Hence the tangle has
an elementary IO-tangle decomposition
$(\Gamma\cap D_1,D_1)\#(\Gamma\cap D_2,D_2)\#
\cdots (\Gamma\cap D_{2d+1},D_{2d+1})$.

Suppose that $\Gamma_m\cap D$ 
does not contain a cycle.
If $\Gamma\cap D$ contains 
more than one white vertex, 
then
by Lemma~\ref{LemSimpleIO-Tangle1} 
the simple IO-tangle is 
an elementary IO-tangle.
If $\Gamma\cap D$ contains 
exactly one white vertex, 
then by Remark~\ref{remWhite1} 
the simple IO-tangle is 
an elementary IO-tangle.
If $\Gamma\cap D$ contains 
no white vertex, 
then by the definition of elementary IO-tangles,
the simple IO-tangle is 
an elementary IO-tangle.
This proves Theorem~\ref{StIITheorem3}.
{\hfill {$\square$}\vspace{1.5em}}

\begin{figure}
\centerline{\includegraphics{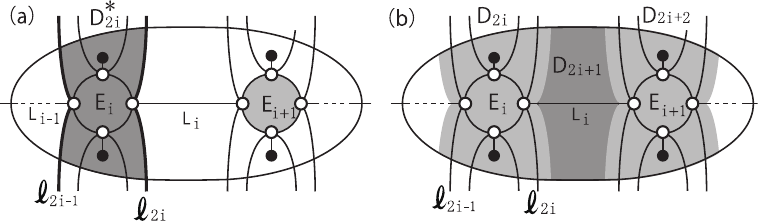}}
\caption{(a) The gray area is the disk $D^*_{2i}$. 
(b) The gray areas are the disks $D_{2i},D_{2i+2}$.
The dark gray area is the disk $D_{2i+1}$.
\LABEL{fig15}}
\end{figure}



\section{Directed paths}
\LABEL{s:OneWayPath}

In this section
we investigate a directed path in $\Gamma_m$.

Let $\Gamma$ be a chart, 
and $m$ a label of $\Gamma$.
A simple arc $P$ in $\Gamma$
is called a {\it path}
provided that 
the endpoints of $P$ are vertices of $\Gamma$.
In particular, 
if the path $P$ is in $\Gamma_m$,
then $P$ 
is called a {\it path of label $m$}.\index{path of label $m$}
Suppose that 
$v_0,v_1,\cdots,v_{p}$ 
are  
all the vertices in a path $P$ 
situated in this order on $P$.
For each $i=1,\cdots,p$, 
let $e_i$ be the edge in $P$ 
with $\partial e_i=\{v_{i-1},v_{i}\}$.
Then the $(p+1)$-tuple 
$(v_0,v_1,\cdots,v_{p})$
is called 
a {\it vertex sequence}\index{vertex sequence}
of $P$,
and the $p$-tuple 
$(e_1,e_2,\cdots,e_{p})$
is called 
an {\it edge sequence}\index{edge sequence}
of $P$.
For two integers $i,j$ 
with $0\le i<j\le p$, 
we denote the path 
$e_{i+1}\cup e_{i+2}\cup\cdots \cup e_{j}$\index{$P[v_i,v_j]$} 
by $P[v_i,v_j]$.

Let $m$ be a label of a chart $\Gamma$, and 
$P$ a path of label $m$
in $\Gamma$ with 
a vertex sequence $(v_0,v_1,\cdots,v_p)$ and
an edge sequence $(e_1,e_2,\dots,e_p)$. 
If for each $i=1,2,\cdots,p$,
the edge $e_i$ is oriented 
from $v_{i-1}$ to $v_i$, 
then $P$ 
is called a {\it directed path}.\index{directed path}
If the edge $e_1$ is middle at $v_0$ 
and if 
the edge $e_p$ is middle at $v_p$, 
then $P$  
is called an {\it M$\&$M path}.\index{M$\&$M path}
A path of label $m$ in a chart $\Gamma$ 
is called 
a {\it dichromatic path}\index{dichromatic path}
if there exists 
a label $k$ with $|m-k|=1$ such that
any vertex of the path 
is contained in $\Gamma_m\cap\Gamma_k$.\\

{\bf Warning.}
If $P$ is a directed path with 
a vertex sequence $(v_0,v_1,\cdots,v_p)$ and
an edge sequence $(e_1,e_2,\dots,e_p)$, 
then {\it we always assume that
\begin{center}
each edge $e_i~(i=1,2,\cdots,p)$ is 
oriented from $v_{i-1}$ to $v_i$.
\end{center}
}

\begin{lemma}
{\rm (\cite[Lemma~3.1]{StI})}
 \LABEL{LemNoM-M}
In a a minimal chart,  
for any label $m$ there does not 
exist
any dichromatic M$\&$M directed path 
of label $m$. {\hfill {$\square$}\vspace{1.5em}}
\end{lemma}

Let $\Gamma$ be a chart, 
and $m$ a label of $\Gamma$.
A simple closed curve in $\Gamma_m$ 
is called a {\it loop}\index{loop} 
if it contains 
exactly one white vertex. 
Let $C$ be 
a cycle of label $m$.
Let  
$v_0,v_1,\cdots,v_{p-1}$ 
be all the vertices in $C$, 
and 
$e_1,e_2,\cdots,e_p$
all the edges in $C$.
Then the cycle $C$  
is called 
a {\it directed cycle} provided that\index{directed cycle}
for each $i=1,2,\cdots,p$,
the edge $e_i$ is oriented from $v_{i-1}$ to $v_i$,
where $v_p=v_0$.
We consider a loop as a directed cycle.

\begin{lemma}
{\rm (\cite[Lemma~3.2]{StI})} 
\LABEL{LemOneWayCycle}
Let $\Gamma$ be a minimal chart, 
$m,k$ labels of $\Gamma$ with $|m-k|=1$, and 
$(\Gamma\cap D,D)$ a tangle
with
$\Gamma\cap D\subset\Gamma_m\cup\Gamma_k$
but without free edges nor simple hoops. 
Then $D$ does not contain any directed cycle 
of label $m$. {\hfill {$\square$}\vspace{1.5em}}
\end{lemma} 

Let $\Gamma$ be a chart, 
and $m$ a label of the chart.
Let $e$ be an edge of label $m$, and
$P$ a directed path of label $m$ 
with an edge sequence 
$(e_1,e_2,\dots,e_p)$.
If $e=e_1$,
then the path $P$ is called 
a directed path {\it starting from} $e$.\index{starting from $e$}
Also if $e=e_p$,
then the path $P$ is called\index{leading to $e$} 
a directed path {\it leading to} $e$. 

Let $P$ be a directed path of label $m$
in a chart, and  
$e$ an edge containing a vertex $v$ 
in Int~$P$ but 
$e\not\subset P$. 
Suppose that $e$ is not a loop.
The edge $e$ is said to be 
{\it locally right-side }
(resp. {\it locally left-side}) at $v$
provided that\index{locally right-side}\index{locally left-side}
for a regular neighbourhood $N$ 
of $v$,
the arc $e\cap N$ 
is situated right (resp. left) side
of $P$ with respect to the direction of $P$.
If the edge $e$ is locally 
right-side at 
a vertex $v$ 
and inward (resp. outward) at $v$, 
then
the edge 
is called a 
locally right-side edge {\it inward}\index{inward} 
(resp. {\it outward}) at $v$.\index{outward}
Similarly 
if the edge $e$ is locally 
left-side at 
a vertex $v$ 
and inward (resp. outward) at $v$, 
then
the edge 
is called a 
locally left-side edge {\it inward} 
(resp. {\it outward}) at $v$. 
In Fig.~\ref{fig16},
the edge $e_2$ is a locally right-side edge 
inward at $v_1$,
the edge $e_3$ is a locally right-side edge 
outward at $v_1$,
the edge $e_4$ is a locally left-side edge 
inward at $v_2$,
the edge $e_5$ is a locally left-side edge 
outward at $v_2$.
But $e_1$ is not 
locally left-side at $v_0$
nor $e_6$ is not 
locally right-side at $v_3$.

\begin{figure}[htb]
\begin{center}
\includegraphics{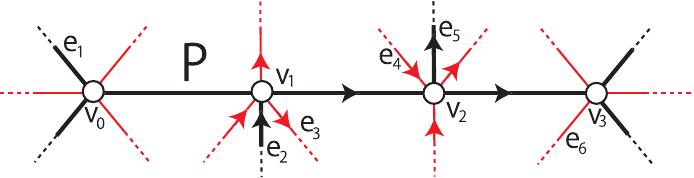}
\end{center}
\caption{ \LABEL{fig16}
The path $P$ is a directed path 
with a vertex sequence
$(v_0,v_1,v_2,v_3)$.}
\end{figure}

Let $\Gamma$ be a chart, 
and $m$ a label of the chart.
A directed path $P$ of label $m$ 
is said to be
{\it upward-right-selective} 
(resp. {\it upward-left-selective})\index{upward-right-selective}\index{upward-left-selective}
if any edge 
of label $m$
locally right-side 
(resp. left-side) 
at a vertex in Int~$P$
is inward at the vertex
(see Fig.~\ref{fig17}(a) and (b)).
A directed path $P$ of label $m$
is said to be
{\it downward-right-selective} 
(resp. {\it downward-left-selective})\index{downward-right-selective}\index{downward-left-selective}
if any edge 
of label $m$ 
locally right-side (resp. left-side) 
at a vertex in Int~$P$
is outward at the vertex 
(see Fig.~\ref{fig17}(c) and (d)).

\begin{figure}[bht]
\begin{center}
\includegraphics{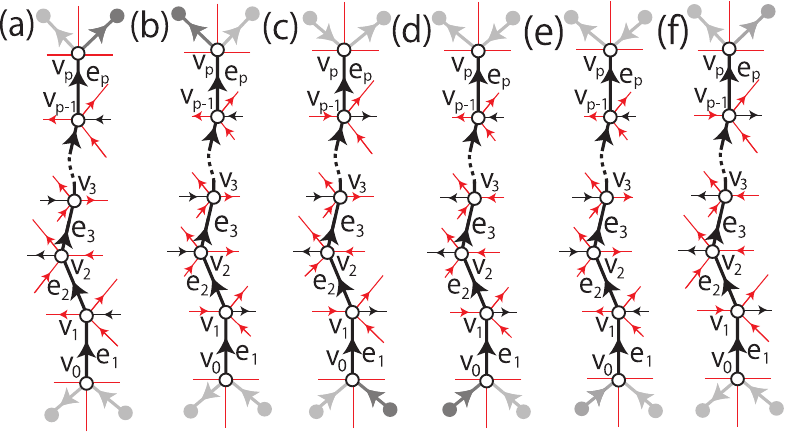}
\end{center}
\caption{\LABEL{fig17} 
(a) The black thick path is an upward-right-selective path. 
(b) The black thick path is an upward-left-selective path.
(c) The black thick path is a downward-right-selective path.
(d) The black thick path is a downward-left-selective path.
(e) The black thick path is an upward principal path.
(f) The black thick path is a downward principal path.}
\end{figure}


Let $\Gamma$ be a chart, and 
$E$ a disk.
An edge $e$ of the chart $\Gamma$
is called 
an {\it I-edge $($resp. O-edge$)$}\index{I-edge}\index{O-edge}
for $E$ provided that (see Fig.~\ref{fig18}(a))
\begin{enumerate}
\item[(i)] 
the edge $e$ possesses two white vertices,
one is in Int~$E$  
and 
the other in $E^c$,
\item[(ii)] 
the edge $e$ intersects $\partial E$ 
by exactly one point, and
\item[(iii)] 
the edge $e$ is inward (resp. outward) at
the vertex in Int~$E$.
\end{enumerate}
We often say just {\it an I-edge} 
instead of {\it an I-edge for $E$}
if there is no confusion.
Similarly we often say just {\it an O-edge} 
instead of {\it an O-edge for $E$}.

\begin{figure}[hbt]
\begin{center}
\includegraphics{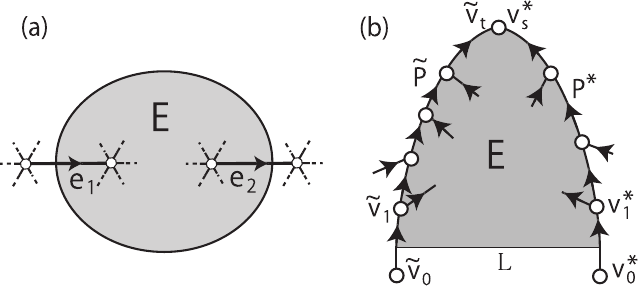}
\end{center}
\caption{ \LABEL{fig18}
(a) The gray area is a disk $E$.
The edge $e_1$ is an I-edge for $E$, and
the edge $e_2$ is an O-edge for $E$. (b) A half spindle.
}
\end{figure}

Let $\Gamma$ be a chart, 
and $m$ a label of $\Gamma$.
Let $P^*$ be 
an upward-right-selective directed path 
of label $m$
with a vertex sequence
$(v^*_0,v^*_1,\cdots,v^*_s)$,
and $\widetilde P$ 
an upward-left-selective directed path 
of label $m$
with a vertex sequence 
$(\widetilde v_0,\widetilde v_1,\cdots,
\widetilde v_t)$
equipped with 
$P^*\cap \widetilde P=v^*_s=\widetilde v_t$.
Let $E$ be a disk 
and $L$ an arc on $\partial E$.
The disk $E$ is called a 
{\it half spindle} for $\Gamma$\index{half spindle}
with an arc triplet 
$(\widetilde P,L,P^*)$
provided that 
(see Fig.~\ref{fig18}(b))
\begin{enumerate}
\item[(i)]
$\Gamma\cap E\subset\Gamma_m\cup\Gamma_k$
for some label $k$ with $|m-k|=1$,
\item[(ii)] 
$(P^*\cup\widetilde P)\cap$~Int~$E=\emptyset$,
and
$E\not\ni v^*_0,\widetilde v_0$, 
\item[(iii)]
$P^*\cap E$ 
is an arc containing 
$P^*[v^*_1,v^*_s]$,
$\widetilde P\cap E$
is an arc containing 
$\widetilde P[\widetilde v_1,\widetilde v_t]$,
and
$L=
Cl(\partial E-(P^*\cup\widetilde P))$,
\item[(iv)] 
the three arcs 
$\widetilde P\cap E,L,P^*\cap E$ 
are counterclockwise situated 
on $\partial E$ in this order,
\item[(v)]
if an edge 
intersects Int~$L$,
then it is an I-edge of label $m$ for $E$, 
\item[(vi)]
$s\ge 2$ and $t\ge 2$.
\end{enumerate}

The following lemma will be used in
the proof of Lemma~\ref{lemTwoPrincipal} 
and Lemma~\ref{lemNoPrincipalQuad}.

\begin{lemma}
{\rm (\cite[Lemma~6.3]{StI})} 
\LABEL{lemHalfSpindle}
For any minimal chart 
there does not exist a half spindle.{\hfill {$\square$}\vspace{1.5em}}
\end{lemma}

\begin{lemma} \LABEL{lemNoTerminal}
Let $\Gamma$ be a minimal chart, and 
$m,k$ labels of $\Gamma$ with 
$|m-k|=1$.
Let $E$ be a disk with 
$\partial E\subset\Gamma$ and
$\Gamma\cap E\subset\Gamma_m\cup\Gamma_k$. 
Then $E$ does not contain a terminal edge.
\end{lemma}

\begin{proof} 
Since $\Gamma\cap E\subset\Gamma_m\cup\Gamma_k$, 
the disk $E$ does not contain a crossing.
Suppose that $E$ contains a terminal edge $e$. 
Without loss of generality 
we can assume that $e$ is of label $m$.
Let $e^*$ be a sibling edge of $e$. 
Then $e^*\subset E$. 
Hence there exists an edge $e^{**}$ of label $m$ 
so that the triplet $(e,e^*,e^{**})$ is 
a non-admissible consecutive triplet. 
This contradicts Consecutive Triplet Lemma 
(Lemma~\ref{ConsecutiveTripletLemma}).
\end{proof}

Let $\Gamma$ be a chart, 
$P$ a path of $\Gamma$, and 
$E$ a disk. 
If each edge in $P$ 
intersects Int~$E$, and 
if $P\cap E$ is connected, 
then we say that the path $P$ 
{\it is dominated by} the disk $E$ or 
the disk $E$ {\it dominates} the path $P$\index{dominate} 
(see Fig.~\ref{fig19}). 

\begin{figure}
\begin{center}
\includegraphics{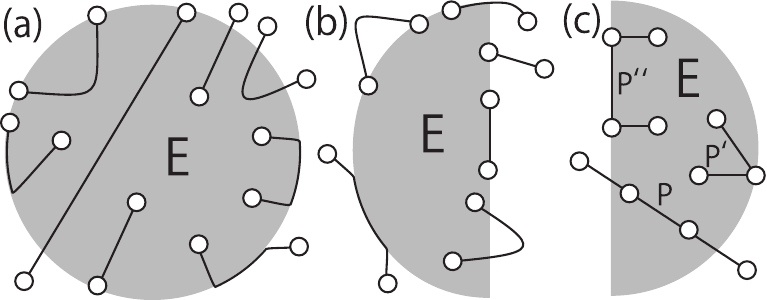}
\end{center}
\caption{\LABEL{fig19} 
(a) The disk $E$ dominates all the edges. 
(b) The disk $E$ dominates none of the edges. 
(c) The disk $E$ dominates the path $P$ and $P'$ but does not the path $P''$.}
\end{figure}

Let $\Gamma$ be a chart, 
and $m$ a label of $\Gamma$. 
Let $E$ be a disk, and
$P$ a directed path 
of label $m$ dominated by $E$
with an edge sequence 
$(e_1,e_2,\cdots,e_p)$. 
If $P$ is not contained in another 
directed path 
of label $m$ 
starting from $e_1$ 
dominated by $E$,
then the path $P$ is said to be 
{\it upward maximal}\index{upward maximal} 
with respect to $E$.
Similarly if $P$ 
is not contained in 
another directed path of label $m$ 
leading to $e_p$ 
dominated by $E$, 
then the path $P$ is said to be 
{\it downward maximal}\index{downward maximal} 
with respect to $E$.

\begin{lemma}\LABEL{lemMaxAtBoundaryUp} 
Let $\Gamma$ be a minimal chart, 
and $m,k$ labels of $\Gamma$  
with $|m-k|=1$.
Let $E$ be a disk with 
$\partial E\subset\Gamma$ and 
$\Gamma\cap E\subset\Gamma_m\cup\Gamma_k$. 
Suppose that $E$ dominates 
a directed path $P$
of label $k$ 
with a vertex sequence 
$(v_0,v_1,v_2,\cdots,v_p)$.
Then we have the following:
\begin{enumerate}
\item[{\rm (a)}] 
If $P$ is 
upward maximal with respect to $E$, 
then $v_p\in\partial E$.
\item[{\rm (b)}] 
If $P$ is downward maximal 
with respect to $E$, 
then $v_0\in\partial E$. 
\end{enumerate}
\end{lemma}

\begin{proof}
{\bf Statement (a)}. 
Since $\Gamma\cap E\subset\Gamma_m\cup\Gamma_k$,
the disk $E$ does not contain 
a terminal edge 
by Lemma~\ref{lemNoTerminal}.
Thus if $v_p\in$~Int~$E$, 
then there exists 
an edge $e$ of label $k$
outward at $v_p$. 
Since there does not exist 
a directed cycle of label $k$
by Lemma~\ref{LemOneWayCycle},
we have $P\cap e=v_p$. 
Hence 
$P\cup e$ is 
a directed path of label $k$
dominated by $E$. 
This contradicts that
the path $P$ is upward maximal with respect to $E$.
Thus Statement (a) holds. 
Similarly we can show Statement (b).
\end{proof}



\section{Principal paths}
\LABEL{s:Principal}

Let $P$ be a directed path of label $m$ 
in a chart with 
a vertex sequence $(v_0,v_1,\cdots,v_p)$ and 
an edge sequence $(e_1,e_2,\cdots,e_p)$. 
The path $P$ is {\it upward principal}\index{upward principal} 
provided that for each  $i=1,2,\cdots,p$ 
the edge $e_i$ is middle at $v_{i-1}$. 
The path $P$ is {\it downward principal}\index{downward principal} 
provided that for each  $i=1,2,\cdots,p$ 
the edge $e_i$ is middle at $v_{i}$. 

\begin{remark}{\rm (\cite[Remark 7.2]{StI})}
\LABEL{remarkPrincipal}
{\rm 
Let $\Gamma$ be a chart, and 
$m$ a label of $\Gamma$.
Let $P$ be a directed path 
of label $m$ in $\Gamma$ 
with a vertex sequence 
$(v_0,v_1,\cdots,v_p)$. 
\begin{enumerate}
\item[(1)]
If $P$ is upward principal,
then  
for any edge $e$ of label $m$ with 
$e\cap P=v_i$ for some $i~(0<i<p)$, 
the edge $e$ is inward at $v_i$ 
(see Fig.~\ref{fig17}(e)).
\item[(2)] 
If $P$ is upward principal, 
then 
$P$ is upward-right-selective and  
upward-left-selective.
\item[(3)] 
If $P$ is downward principal 
for any edge $e$ of label $m$ with 
$e\cap P=v_i$ for some $i~(0<i<p)$, 
then 
the edge $e$ is outward at $v_i$ 
(see Fig.~\ref{fig17}(f)).
\item[(4)] 
If $P$ is downward principal, 
then 
$P$ is downward-right-selective and  
downward-left-selective.
\end{enumerate}
} 
\end{remark}

\begin{lemma}{\rm (\cite[Lemma 7.3]{StI})}
\LABEL{lemPrincipalPath}
Let $\Gamma$ be a minimal chart, and 
$m$ a label of $\Gamma$.
Let $P$ be a dichromatic directed path 
of label $m$ in $\Gamma$ 
with a vertex sequence 
$(v_0,v_1,\cdots,v_p)$ and 
an edge sequence 
$(e_1,e_2,\cdots,e_p)$. 
\begin{enumerate}
\item[{\rm (a)}] 
If $e_1$ is middle at $v_0$, 
then $P$ is upward principal.
\item[{\rm (b)}] 
If $e_p$ is middle at $v_p$, 
then $P$ is downward principal.
\end{enumerate}
\end{lemma}

Let $\Gamma$ be a chart, 
and $E$ a disk. 
Let $e$ be an edge of $\Gamma$ 
such that 
$e\cap \partial E$ consists of 
one white vertex or two white vertices. 
If $e\subset Cl(E^c)$, then 
we call $e$  
an {\it outside edge} for $E$\index{outside edge} 
\index{outside edge}. 
If $e\subset E$, then 
we call $e$  
an {\it inside edge} for $E$ 
\index{inside edge}.\index{inside edge}
Let $X$ be a subset of $\partial E$. 
If an outside edge (resp. inside edge) 
for $E$ intersects $X$, then 
we call the edge an {\it outside edge} 
(resp. {\it inside edge}) 
{\it for} $(E,X)$.

Let $v$ be a white vertex 
of a chart $\Gamma$, and 
$X_1,X_2,X_3$ edges or paths in $\Gamma$ 
with $v\in\partial X_i~~(i=1,2,3)$. 
In a regular neighbourhood of $v$ in $\Gamma$, 
there are six short arcs at $v$. 
For each $i=1,2,3$, let 
$\gamma_i$ be the one 
of the six short arcs in $X_i$. 
If the three short arcs 
$\gamma_1,\gamma_2,\gamma_3$ are 
consecutive among the six short arcs 
in this order, 
then $X_2$ is said to be 
{\it situated between 
$X_1$ and $X_3$ around $v$}\index{situated between $X_1$ and $X_3$ around $v$} 
(see Fig.~\ref{fig20}).


\begin{figure}
\begin{center}
\includegraphics{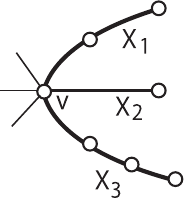}
\end{center}
\caption{\LABEL{fig20} 
The edge $X_2$ is situated between $X_1$ and $X_3$ around the vertex $v$.}
\end{figure}


\begin{lemma}
\LABEL{lemTwoPrincipal}
Let $\Gamma$ be a minimal chart 
in a disk $D^2$, and 
$m,k$ labels of $\Gamma$ 
with $|m-k|=1$. 
Let $E$ be a disk 
with $\partial E\subset \Gamma$ and
$\Gamma\cap E\subset\Gamma_m\cup\Gamma_k$. 
Let $P^\dagger,P^*$ be 
directed paths in $E$  
of label $k$ 
with vertex sequences 
$(v^\dagger_0,v^\dagger_1,\cdots,v^\dagger_p),
(v^*_0,v^*_1,\cdots,v^*_q)$ respectively. 
Then we have the following.
\begin{enumerate}
\item[{\rm (a)}] 
If $P^\dagger,P^*$ are 
upward principal, 
and if 
there exists a path $L$ 
on $\partial E$ of label $m$ 
with 
$\partial L= \{v^\dagger_0,v^*_0\}$ 
such that
\begin{enumerate}
\item[{\rm (i)}]
any outside edge for $(E,{\rm Int}~L)$ 
is of label $k$ inward at 
a white vertex in {\rm Int}~$L$,
and 
\item[{\rm(ii)}]
for each white vertex $v$ in $\partial L$, 
there exists an outside edge for $E$ 
of label $k$ at $v$,
\end{enumerate}
then $P^\dagger\cap P^*=\emptyset$.
\item[{\rm (b)}]
If $P^\dagger,P^*$ 
are downward principal, 
and if 
there exists a path $L$ 
on $\partial E$ of label $m$ 
with 
$\partial L= \{v^\dagger_p,v^*_q\}$
such that
\begin{enumerate}
\item[{\rm (i)}]
any outside edge for $(E,{\rm Int}~L)$ 
is of label $k$ outward at 
a white vertex in {\rm Int}~$L$,
and 
\item[{\rm(ii)}]
for each white vertex $v$ in $\partial L$, 
there exists an outside edge for $E$ 
of label $k$ at $v$,
\end{enumerate}
then $P^\dagger\cap P^*=\emptyset$.
\end{enumerate}
\end{lemma}

\begin{proof}
Let $(e^\dagger_1,e^\dagger_2,\cdots,e^\dagger_p), 
(e^*_1,e^*_2,\dots,e^*_q)$ 
be edge sequences of 
$P^\dagger,P^*$ respectively.
Suppose that 
$P^\dagger\cap P^*\neq\emptyset$. 

{\bf Statement (a)}. 
Let 
$s=\min\{i~|~v^\dagger_i\in P^*\}$. 
Then $v^\dagger_s=v^*_t$ 
for some $t~(1\le t\le q)$. 
By Condition (ii), 
there exists an outside edge $e^\dagger$ 
for $E$ 
of label $k$ 
at $v^\dagger_0$ such that 
no outside edge  
for $E$ 
is situated between 
$e^\dagger$ and $L$ 
around $v^\dagger_0$.
Since $e^\dagger_1$ is 
middle at $v^\dagger_0$, 
the edge $e^\dagger$ is 
inward at $v^\dagger_0$. 
Hence $e^\dagger\cup P^\dagger$ 
is a directed path of label $k$.
Let $e^*$ be 
an outside edge for $E$ 
of label $k$ 
at $v^*_0$ such that 
no outside edge  
for $E$ 
is situated between 
$e^*$ and $L$ around $v^*_0$.
Similarly $e^*\cup P^*$ 
is a directed path of label $k$.
Then we can get a half spindle 
by splitting 
a regular neighbourhood of $E$ in $D^2$
by $e^\dagger\cup 
P^\dagger[v^\dagger_0,v^\dagger_s]
\cup P^*[v^*_0,v^*_t]\cup e^*$. 
This contradicts Lemma~\ref{lemHalfSpindle}.

{\bf Statement (b)}. 
Changing orientations of 
all the edges of $\Gamma$, 
the minimal chart $\Gamma$ changes into 
a new minimal chart and 
the paths change into upward principal 
paths in the new minimal chart. 
Thus we can get a contradiction 
by Statement (a). 
\end{proof}

Let $\Gamma$ be a chart, and 
$m,k$ labels of $\Gamma$ with $|m-k|=1$. 
Let $L$ be a dichromatic path of label $m$ 
with two end points $v^\dagger,v^*$, 
and 
$e^\dagger,e^*$ 
edges of label $k$
middle at the white vertices 
$v^\dagger,v^*$ respectively. 
Then the triplet\index{staple}
$(e^\dagger,L,e^*)$ 
is called a {\it staple 
with label pair} $(m,k)$ 
(see Fig.~\ref{fig21}(a))
provided that
\begin{enumerate}
\item[(i)] 
$e^\dagger,~e^*$ 
are outward (resp. inward) 
at $v^\dagger,~v^*$ respectively, and
\item[(ii)] 
there exists a disk $D$ such that 
\begin{enumerate}
\item[(a)] 
set $\ell=
\partial D\cap(e^\dagger\cup L\cup e^*)$, 
then $\ell$ is an arc with 
$L\subset~{\rm Int}~\ell$,
\item[(b)]
each edge dominated by $D$ is 
outward (resp. inward) at a vertex on $L$.
\end{enumerate} 
\end{enumerate} 
The disk $D$ is called an {\it associated disk 
for the staple}.

\begin{figure}
\begin{center}
\includegraphics{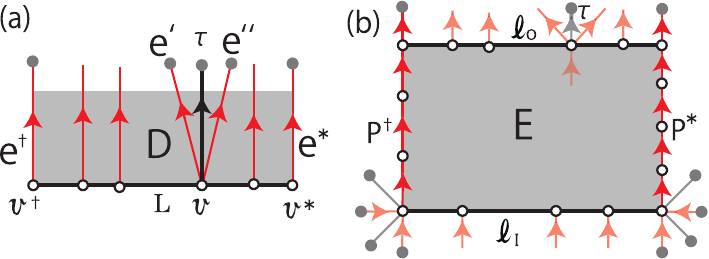}
\end{center}
\caption{\LABEL{fig21}
The gray vertices are 
black vertices, white vertices, 
or crossings.
(a) a staple $(e^\dagger,L,e^*)$ 
with an associated disk $D$.
(b) a principal quad 
$(P^\dagger,P^*,\ell_I,\ell_O)$ 
bounding a disk $E$. 
}
\end{figure}

\begin{lemma} \LABEL{lemStaple}
Let $\Gamma$ be a minimal chart, 
and $m,k$ labels of $\Gamma$
with $|m-k|=1$.
If $(e^\dagger,L,e^*)$ 
is a staple with label pair $(m,k)$, 
then any associated disk for the staple
dominates 
an edge of label $m$ 
middle at a white vertex in {\rm Int}~$L$.
\end{lemma}

\begin{proof}
We use all the notations in the definition 
of a staple. 
Let $D$ be an associated disk 
for the staple $(e^\dagger,L,e^*)$.
We only show the case 
that each of $e^\dagger,e^*$ is outward 
at $v^\dagger,v^*$ respectively.
We can show for the other case similarly.
It suffices to show 
for the case that 
\begin{enumerate}
\item[(1)] 
the disk $D$ does not dominate any edge 
of label $k$ middle at a white vertex in Int~$L$.
\end{enumerate}

{\bf Claim 1}.
The disk $D$ does not dominate 
an edge of label $m$ 
at $v^\dagger$ nor $v^*$.

For, if the disk $D$ dominates 
an edge $e'$ of label $m$ at $v^\dagger$, 
then $D$ also dominates an edge $e''$ 
of label $k$ at $v^\dagger$. 
Since $e'$ and $e''$ are dominated by $D$, 
the edges $e',e''$ are 
outward at $v^\dagger$. 
Since $e^\dagger$ is 
middle at $v^\dagger$,
there exists an edge $e'''$ of label $m$ 
outward at $v^\dagger$ 
not dominated by $D$. 
Hence the four edges 
$e^\dagger,e',e'',e'''$ are 
outward at $v^\dagger$. 
This contradicts 
the definition of charts. 
Similarly $D$ does not dominate 
an edge of label $m$ at $v^*$. 
Thus Claim1 holds.

{\bf Claim 2}.
The disk $D$ dominates 
an edge $\tau$ of label $m$ 
at a vertex $v$ in Int~$L$.

For, if the disk  
$D$ does not dominate 
an edge of label $m$ 
at a vertex in Int~$L$. 
Then $L$ is a directed path of label $m$ 
by (1).
Since $e^\dagger$ is outward 
at $v^\dagger$, 
the path $L$ is oriented 
from $v^\dagger$ to $v^*$ 
by Claim 1. 
Since $e^*$ is 
middle at $v^*$, 
the edge $e^*$ is 
inward at $v^*$ by Claim 1. 
This contradicts the assumption that
the edge $e^*$ is 
outward at $v^*$. 
Thus Claim 2 holds. 

Hence $D$ dominates an edge $\tau$ 
of label $m$ 
at a white vertex $v$ in Int~$L$. 
Then $D$ dominate 
two edges of label $k$; 
an edge $e'$ situated between 
$\tau$ and $L[v^\dagger,v]$ 
and 
an edge $e''$ situated between 
$\tau$ and $L[v,v^*]$ 
around the vertex $v$. 
Then the edges $e',e''$ are outward at $v$. 
Since the edge $\tau$ is situated between 
$e'$ and $e''$ around $v$, 
the edge $\tau$ is middle at $v$.
This proves Lemma~\ref{lemStaple}.
\end{proof}

Let $\Gamma$ be a chart 
in a disk $D^2$, 
and $k,m$ labels of $\Gamma$ with $|k-m|=1$.  
Let $P^\dagger,P^*$ be disjoint 
upward principal paths of label $k$  
with vertex sequences 
$(v^\dagger_0,v^\dagger_1,\cdots,v^\dagger_p),
(v^*_0,v^*_1,\cdots,v^*_q)$ and 
edge sequences 
$(e^\dagger_1,e^\dagger_2,\cdots,e^\dagger_p),
(e^*_1,v^*_2,\cdots,e^*_q)$ respectively, 
here $p\geq 2,q\geq 2$. 
Let $\ell_I,\ell_O$ be 
disjoint paths of label $m$ 
with 
$\partial\ell_I
=\{v^\dagger_0,v^*_0\}$ and
$\partial\ell_O
=\{v^\dagger_{p-1},v^*_{q-1}\}$. 
The quadruplet $(P^\dagger,P^*,\ell_I,\ell_O)$ 
is called a {\it principal quad 
with label pair $(k,m)$}\index{principal quad}
provided that (see Fig.~\ref{fig21}(b))
\begin{enumerate}
\item[(i)] 
$\ell_I\cap(P^\dagger\cup P^*)
=\partial\ell_I$ and 
$\ell_O\cap(P^\dagger\cup P^*)
=\partial\ell_O$,
\item[(ii)] for the disk $E$ bounded by 
$P^\dagger[v^\dagger_0,v^\dagger_{p-1}]\cup 
\ell_O\cup P^*[v^*_0,v^*_{q-1}]\cup \ell_I$, 
\begin{enumerate}
\item[(a)] $e^\dagger_p,e^*_q\not\subset E$, 
and
$\Gamma\cap E\subset \Gamma_k\cup\Gamma_m$,
\item[(b)] 
each outside edge for $(E,{\rm Int}~\ell_I)$ 
is of label $k$ and 
inward at a vertex in Int $\ell_I$,
\item[(c)] 
each outside edge 
for $(E,{\rm Int}~\ell_O)$ 
is 
outward at a vertex in Int $\ell_O$,
\item[(d)] 
for each vertex $v\in\partial\ell_I$, 
there does not exist any inside edge 
for $E$ at $v$.
\end{enumerate}
\end{enumerate}
The disk $E$ is called an 
{\it associated disk} for the principal quad. 

\begin{lemma}
\LABEL{lemPrincipalQuad}
Let $\Gamma$ be a minimal chart, and 
$k,m$ labels of $\Gamma$ with 
$|k-m|=1$. 
Let $(P^\dagger,P^*,\ell_I,\ell_O)$ 
be a principal quad in $\Gamma$ with 
an associated disk $E$ and  
a label pair $(k,m)$. 
Then we have the following.
\begin{enumerate}
\item[{\rm (a)}] 
Among the outside edges 
for $(E,{\rm Int}~\ell_I)$ of label $k$, 
at most one edge 
is middle at a vertex in {\rm Int}~$\ell_I$.
\item[{\rm (b)}]
Suppose that 
$E$ dominates a directed path $P$ 
of label $k$ 
with a vertex sequence $(v_0,v_1,\cdots,v_r)$,
an edge sequence $(e_1,e_2,\cdots,e_r)$, 
and $v_r\not\in P^\dagger\cup P^*$. Then
\begin{enumerate}
\item[{\rm (i)}] 
$P\cap(P^\dagger\cup P^*)=\emptyset$, and
\item[{\rm (ii)}]
if $P$ is downward maximal 
with respect to $E$,
then $v_0$ is in {\rm Int}~$\ell_I$ and 
there exist an outside edge for 
$(E,{\rm Int}~\ell_I)$ 
of label $k$ middle at $v_0$ and 
an inside edge for $(E,{\rm Int}~\ell_I)$ 
of label $k$ 
outward at $v_0$
different from $e_1$.
\end{enumerate}
\item[{\rm (c)}]
Among the outside edges 
for $(E,{\rm Int}~\ell_O)$ of label $m$, 
at most one edge  
is at a vertex in {\rm Int}~$\ell_O$.
\item[{\rm (d)}]
No outside edge 
for $(E,{\rm Int}~\ell_O)$ of label $k$ 
is middle at a vertex in {\rm Int}~$\ell_O$.
\end{enumerate}
\end{lemma}

\begin{proof}
Let 
$(v^\dagger_0,v^\dagger_1,\cdots,v^\dagger_p)$
be a vertex sequence of $P^\dagger$, and 
$(v^*_0,v^*_1,\cdots,v^*_q)$ 
a vertex sequence of $P^*$.

{\bf Statement (a)}.
Suppose there exist 
two outside edges $e',e''$ for 
$(E,{\rm Int}~\ell_I)$ 
of label $k$ middle at white vertices $v',v''$ 
in Int~$\ell_I$ respectively. 
Since $(e',\ell_I[v',v''],e'')$ 
is a staple with label pair $(m,k)$, 
there exists an outside edge for $E$ 
of label $m$ middle at a vertex 
in Int~$\ell_I[v',v'']\subset$~Int~$\ell_I$. 
This contradicts Condition (ii)(b) of 
a principal quad.
Thus Statement (a) holds.

{\bf Statement (b)(i)}. 
If $P\cap P^\dagger\neq\emptyset$, 
let $s=\max\{i~|~v_i\in P^\dagger\}$.
Then $v_r\not\in P^\dagger\cup P^*$ 
implies $0\leq s<r$. 
Now 
$v^\dagger_t=v_s$ 
for some integer $t~(0<t<p)$. 
Then $e_{s+1}$ is outward at 
$v_s=v^\dagger_t$. 
This contradicts 
Remark~\ref{remarkPrincipal}(1) 
because $P^\dagger$ is upward principal.
Thus $P\cap P^\dagger=\emptyset$. 
Similarly $P\cap P^*=\emptyset$.
Thus Statement (b)(i) hold.

{\bf Statement (b)(ii)}.
Now $v_0\not\in P^\dagger\cup P^*$ 
by Statement (b)(i), and 
$v_0\in\partial E$ 
by Lemma~\ref{lemMaxAtBoundaryUp}(b).

If $v_0\in{\rm ~Int~}\ell_O$, then 
there exists an outside edge for 
$(E,{\rm Int}~\ell_O)$ 
of label $k$ outward at $v_0$. 
Further  
the edge $e_1$ is an inside edge for 
$(E,{\rm Int}~\ell_O)$ 
of label $k$ outward at $v_0$, 
because $P$ is a directed path. 
Hence 
there exists an edge $e'$ 
of label $k$ inward at $v_0$. 
Then $e'$ is an inside edge for $E$ 
by Condition (ii)(c) of 
a principal quad. 
Thus $P\cup e'$ is a directed path 
of label $k$ 
leading to $e_r$ dominated by $E$. 
This contradicts the fact that 
$P$ is downward maximal with respect to $E$.
Thus $v_0\not\in{\rm ~Int~}\ell_O$.

Hence $v_0\in {\rm ~Int~}\ell_I$. 
There exists an outside edge $e$ for 
$(E,{\rm Int}~\ell_I)$ 
of label $k$ inward at $v_0$. 
Let $e'$ be the third edge of label $k$ 
at $v_0$. 
If $e'$ is an outside edge for $E$, 
then there exists an edge of label $m$ 
between $e$ and $e'$ around $v_0$. 
This contradicts Condition (ii)(b) 
of a principal quad. 
Thus $e'$ is an inside edge for 
$(E,{\rm Int}~\ell_I)$. 
If $e'$ is inward at $v_0$, 
again $P\cup e'$ is a directed path  
leading to $e_r$ dominated by $E$. 
This contradicts the fact that 
$P$ is downward maximal with respect to $E$. 
Thus $e'$ is outward at $v_0$.
Hence the outside edge $e$ 
is middle at $v_0$.
This proves Statement (b)(ii).
 
{\bf Statement (c)}.
Suppose there exist
two outside edges $\tau',\tau''$ 
for $(E,{\rm Int}~\ell_O)$ 
of label $m$.
Let $v'=\tau'\cap\ell_O, 
v''=\tau''\cap\ell_O$.
Then $\tau',\tau''$ are middle at 
$v',v''$ respectively 
by Condition (ii)(c) of 
a principal quad.
Let $e',e''$ be the inside edges for 
$(E,{\rm Int}~\ell_O)$ 
of label $k$ 
at $v',v''$ respectively. 
Then $e',e''$ are inward and
middle at $v',v''$ respectively. 
Without loss of generality 
we can assume that 
\begin{enumerate}
\item[(1)] there does not exist
an outside edge for 
$(E,{\rm Int}~\ell_O[v',v''])$ 
of label $m$.
\end{enumerate}
Let $P',P''$ be directed paths 
of label $k$ 
leading to $e',e''$ 
downward maximal with respect to $E$ 
respectively. 

Since $e',e''$ are 
middle at $v',v''$ respectively, 
the paths $P',P''$ are downward principal paths 
dominated by $E$ 
by Lemma~\ref{lemPrincipalPath}(b). 
Hence  
$P'\cap P''=\emptyset$ 
by (1) and Lemma~\ref{lemTwoPrincipal}(b). 
Thus by Statement (b)(ii)  
there exist two outside edges $e^*,e^{**}$ for 
$(E,{\rm Int}~\ell_I)$ 
of label $k$ middle at 
vertices $v^*,v^{**}$ 
in Int~$\ell_I$ respectively.
This contradicts Statement (a).
Thus Statement (c) holds.

{\bf Statement (d)}.
Let 
$(e^\dagger_1,e^\dagger_2,\cdots,e^\dagger_p)$
be an edge sequence of $P^\dagger$, and 
$(e^*_1,e^*_2,\cdots,e^*_q)$ 
an edge sequence of $P^*$.
Suppose that 
there exists 
an outside edge $e'$ for 
$(E,{\rm Int}~\ell_O)$ 
of label $k$ 
middle at a white vertex $v'$ 
in Int~$\ell_O$.
Then $(e^\dagger_p,
\ell_O[v^\dagger_{p-1},v'],e')$ 
and 
$(e',\ell_O[v',v^*_{q-1}],e^*_q)$ 
are staples with label pair $(m,k)$. 
Hence 
by Lemma~\ref{lemStaple}, 
there exist two outside edges for 
$(E,{\rm Int}~\ell_O)$ 
of label $m$; 
an outside edge for 
$(E,{\rm Int}~\ell_O[v^\dagger_{p-1},v'])$ 
middle at a vertex in 
Int~$\ell_O[v^\dagger_{p-1},v']$, and 
an outside edge for 
$(E,{\rm Int}~\ell_O[v',v^*_{q-1}])$ 
middle at a vertex in 
Int~$\ell_O[v',v^*_{q-1}]$.
This contradicts Statement (c).
Thus Statement (d) holds.
\end{proof}

\begin{lemma}
\LABEL{lemNoPrincipalQuad}
There does not exist a principal quad in 
a minimal chart.
\end{lemma}

\begin{proof}
Suppose that there exists a principal quad 
$(P^\dagger,P^*,\ell_I,\ell_O)$ 
with a label pair $(k,m)$ 
in a minimal chart $\Gamma$.

Let $(v^\dagger_0,v^\dagger_1,\cdots,v^\dagger_p)$
and 
$(e^\dagger_1,e^\dagger_2,\cdots,e^\dagger_p)$ 
are a vertex sequence and an edge sequence 
of the upward principal path $P^\dagger$ 
of label $k$ respectively,  
here $p\geq 2$.

Let $(v^*_0,v^*_1,\cdots,v^*_q)$ and  
$(e^*_1,e^*_2,\cdots,e^*_q)$  
be a vertex sequence and an edge sequence 
of the upward principal path $P^*$ 
of label $k$ respectively, 
here $q\geq 2$.

The associated disk $E$ 
for the principal quad
is the disk bounded by 
$P^\dagger[v^\dagger_0,v^\dagger_{p-1}]\cup 
\ell_O\cup P^*[v^*_0,v^*_{q-1}]\cup \ell_I$. 
Since $(e^\dagger_p,\ell_O,e^*_q)$ 
is a staple with label pair $(m,k)$, 
there exists an outside edge $\tau'$ for 
$(E,{\rm Int}~\ell_O)$ 
of label $m$ 
middle at a vertex $v'$ in Int~$\ell_O$ 
by Lemma~\ref{lemStaple}. 
Without loss of generality we can assume that
\begin{enumerate}
\item[(1)] 
$v^\dagger_0,v^*_0,v'$ are situated 
counterclockwise in this order 
on $\partial E$.
\end{enumerate} 
Let $e'$ be the inside edge for 
$(E,{\rm Int}~\ell_O)$ 
of label $k$ 
middle at $v'$. 
Then $e'$ is inward at $v'$. 
Let $P'$ be a directed path of label $k$  
leading to $e'$ 
downward maximal with respect to $E$ 
with a vertex sequence 
$(v'_0,v'_1,\cdots, v'_r)$ 
and 
an edge sequence 
$(e'_1,e'_2,\cdots, e'_r)$, 
here $v'_r=v'$ and $e'_r=e'$. 
By Lemma~\ref{lemPrincipalQuad}(b)(ii), 
we have $v'_0\in\ell_I$.

{\bf Claim}. 
$P'\cap\ell_I=v'_0$.\\
For, if $v'_i\in {\rm Int~}\ell_I$ 
for some integer $i~(0<i<r)$, 
then 
there exists an outside edge for 
$(E,{\rm Int}~\ell_I)$ 
of label $k$ 
inward at $v'_i$. 
Since $e'_i$ is inward at $v'_i$, 
the edge $e'_{i+1}$ is middle at $v'_i$. 
Hence $P'[v'_{i},v'_r]$ is 
a dichromatic M\&M directed path. 
This contradicts Lemma~\ref{LemNoM-M}. 
Since $P'\cap(P^\dagger\cup P^*)=\emptyset$ 
by Lemma~\ref{lemPrincipalQuad}(b)(i), 
Claim follows 
from Lemma~\ref{lemPrincipalQuad}(b)(ii).

Thus the path $P'$ splits 
the disk $E$ into disks.
There exist two disks; 
a disk intersecting 
both of $P^\dagger$ and $P'$,
say $\Delta^\dagger$,
and
a disk intersecting 
both of $P'$ and $P^*$, 
say $\Delta^*$
(see Fig.~\ref{fig22}(a)).
By Lemma~\ref{lemPrincipalQuad}(b)(ii),
there exists 
an inside edge $\widehat e$ for 
$(E,{\rm Int}~\ell_I)$ 
of label $k$ outward at $v'_0$ 
different from $e'_1$.
By Claim, 
one of 
the two disks $\Delta^\dagger, \Delta^*$
containing the edge $\widehat e$, 
say $\widehat \Delta$. 
Without loss of generality 
we can assume that 
$\widehat \Delta = \Delta^\dagger$.

\begin{figure}
\begin{center}
\includegraphics{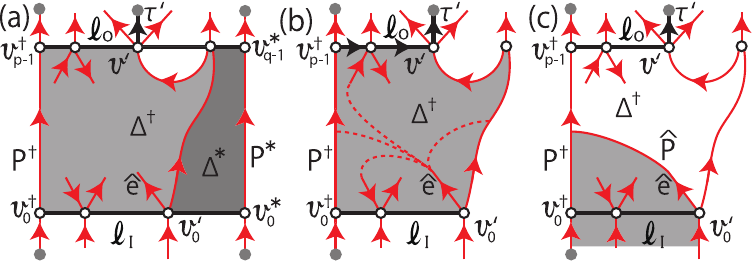}
\end{center}
\caption{\LABEL{fig22}
(c) The gray area is a half spindle.}
\end{figure}

Let $\widehat P$ be 
an upward-right-selective 
(upward-left-selective 
if $\widehat \Delta=\Delta^*$)
directed path
of label $k$
starting from $\widehat e$ 
dominated by $\widehat \Delta$ 
with a vertex sequence 
$(\widehat v_0,\widehat v_1,\widehat v_2,
\cdots,\widehat v_s)$,  
an edge sequence 
$(\widehat e_1,\widehat e_2,
\cdots,\widehat e_s)$, and 
$\widehat v_0=v'_0$ 
such that $\widehat P$ is 'maximal' 
with respect to $\widehat\Delta$ 
in the sense that 
the path $\widehat P$ is not 
contained in 
another upward-right-selective 
(upward-left-selective 
if $\widehat \Delta=\Delta^*$)
directed path starting from $\widehat e$
dominated by $\widehat\Delta$. 

Then we can show 
$\widehat v_s\in\partial\Delta^\dagger$ 
by the same way as the one of 
Lemma~\ref{lemMaxAtBoundaryUp}(a). 
There are four cases 
(see Fig.~\ref{fig22}(b)):
\begin{enumerate}
\item[] Case 1. 
$\widehat v_s\in 
P^\dagger\cap \Delta^\dagger$.
\item[] Case 2. 
$\widehat v_s\in P'\cap \Delta^\dagger$.
\item[] Case 3. 
$\widehat v_s\in$ 
Int $\ell_I[v^\dagger_0,v'_0]$.
\item[] Case 4. 
$\widehat v_s\in$~Int~$\ell_O[v^\dagger_{p-1},v'_r]$.
\end{enumerate}

{\bf Case 1}.
Let $v^\dagger_i=\widehat v_s$. 
Considering Condition (ii)(d) 
of a principal quad, 
we have $0< i < p$. 
Since the path $P^\dagger$ is upward principal, 
the path is upward-left-selective 
by Remark~\ref{remarkPrincipal}(2), 
and so is $P^\dagger[v^\dagger_0,v^\dagger_i]$. 
Since $\widehat P$ is upward-right-selective, 
there exists a half spindle containing 
$\widehat P\cup
P^\dagger[v^\dagger_0,v^\dagger_i]$
(see Fig.~\ref{fig22}(c)). 
This contradicts Lemma~\ref{lemHalfSpindle}.

{\bf Case 2}. 
Let $v'_j=\widehat v_s$, 
here $0<j<r$.
Since $\widehat e_s$ 
is inward at $v'_j$ and 
since $e'_j$ is inward at $v'_j$, 
the edge $e'_{j+1}$ is middle at $v'_j$.
Hence the path $P'[v'_j,v'_r]$ 
is a dichromatic M$\&$M directed path.
This contradicts Lemma~\ref{LemNoM-M}.

{\bf Case 3}. 
Since there does not exist 
any outside (terminal) edge 
of label $m$ for $\Delta^\dagger$
intersecting Int $\ell_I$ 
by Condition (ii)(b) 
of the definition of a principal quad, 
there exists only one outside edge 
for $\Delta^\dagger$ 
at $\widehat v_s$.
Now the edge is of label $k$ 
and inward at $\widehat v_s$.
Since $\widehat e_s$ is inward at $\widehat v_s$, 
there exists an inside edge $\overline e$ 
for $\Delta^\dagger$ 
of label $k$ outward at $\widehat v_s$ 
by considering Condition (ii)(b) 
of the definition of a principal quad. 
Thus $\widehat P\cup \overline e$ is 
an upward-right-selective directed path 
dominated by $\Delta^\dagger$.
This contradicts the fact that
$\widehat P$ is maximal 
with respect to $\Delta^\dagger$.

{\bf Case 4}. 
If there exist two outside edges 
for $\Delta^\dagger$ of label $k$ 
at $\widehat v_s\in$~Int~$\ell_O$, 
then there exists outside edge 
for $\Delta^\dagger$ 
of label $m$ middle at $\widehat v_s$ 
different from $\tau'$. 
This contradicts Lemma~\ref{lemPrincipalQuad}(c).
Thus there exists 
exactly one outside edge $\overline e'$  
for $\Delta^\dagger$ of label $k$ 
at $\widehat v_s$. 
Thus there exists an inside edge $\overline e''$ 
for $\Delta^\dagger$ of label $k$ 
at $\widehat v_s$ 
different from $\widehat e_s$.
Let $\overline e'''$ be the edge of label $m$ 
situated between $\overline e''$ 
and $\widehat e_s$ around $\widehat v_s$.  
If $\overline e''$ is inward at $\widehat v_s$, 
then $\overline e'$ is 
middle at $\widehat v_s$. 
This contradicts Lemma~\ref{lemPrincipalQuad}(d).
Thus $\overline e''$ is outward at $\widehat v_s$. 
Since there does not exist 
any outside edge for $\Delta^\dagger$ 
of label $k$ nor label $m$
middle at a vertex in
Int~$\ell_O[v^\dagger_{p-1},v']$ 
by Lemma~\ref{lemPrincipalQuad}(c),(d), 
the path $\ell_O[v^\dagger_{p-1},v']$ 
is oriented from $v^\dagger_{p-1}$ to $v'$ 
(see Fig. \ref{fig22}(b)). 
Hence the edge $\overline e''$ is outward 
at $\widehat v_s$ and situated 
between $\overline e'''$ 
and $\ell_O[\widehat v_s,v']$ 
around $\widehat v_s$. 
And $\widehat e_s,\overline e'',\overline e'$ 
are counterclockwise 
situated around $\widehat v_s$  
in this order by Statement (1).
Thus $\widehat P\cup \overline e''$ is 
an upward-right-selective directed path 
dominated by $\Delta^\dagger$. 
This contradicts the fact that 
$\widehat P$ is 
maximal with respect $\Delta^\dagger$. 
Since for each of four cases 
we got a contradiction, 
there does not exist any principal quad 
in a minimal chart. 
\end{proof}


\section{Proof of Theorem~1.4}
\LABEL{s:Th4}

Throughout this section we assume that
(see Fig.~\ref{fig23})
\begin{enumerate}
\item[$\bullet$] 
$\Gamma$ is a  minimal chart, 
and $m$ is a label of $\Gamma$, 
\item[$\bullet$] 
$(\Gamma\cap D,D)$ is 
a Type-I elementary IO-tangle of label $m$
with a boundary IO-arc pair $(L_I,L_O)$,
\item[$\bullet$] 
$k$ is a label of $\Gamma$ with $|m-k|=1$
and $\Gamma\cap D\subset\Gamma_m\cup\Gamma_k$,
\item[$\bullet$] 
$e_I,e_O$ are the edges of label $m$ 
intersecting $\partial D$,
\item[$\bullet$]
$\Delta$ is a disk 
containing 
all the white vertices in $D$
with $\partial\Delta\subset 
\Gamma_m\cap D$, 
\item[$\bullet$]
$\Delta_I$ 
(resp. $\Delta_O$) 
is the closure of 
a connected component of 
$D-(\Delta\cup e_I\cup e_O)$ with
$L_I\subset \Delta_I$ 
(resp. $L_O\subset \Delta_O$), and
\item[$\bullet$]
$J_I$ (resp. $J_O$) is 
the closure of 
a connected component of 
$\partial \Delta-(e_I\cup e_O)$ with
$J_I\subset \Delta_I$ 
(resp. $J_O\subset \Delta_O$).
\end{enumerate}


\begin{figure}
\begin{center}
\includegraphics{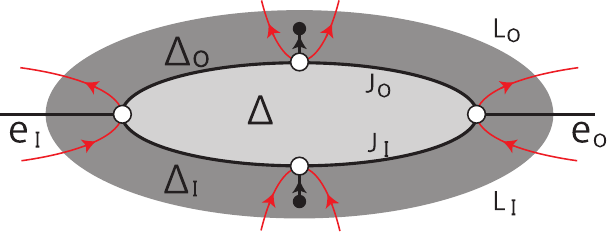}
\end{center}
\caption{ \LABEL{fig23} 
Thick arcs are of label $m$. Thin arcs are of label $k$.}
\end{figure}


Since $\Gamma_m\cap D$ contains a cycle, 
the tangle is simple by 
Lemma~\ref{TwoColorGateTangle}.
Thus we have the following 
by the definition of a simple IO-tangle.
\begin{enumerate}
\item[(I)] 
any edge dominated by $\Delta_I$ is 
inward at a vertex in $J_I$.
\item[(II)] 
any edge dominated by $\Delta_O$ is 
outward at a vertex in $J_O$.
\item[(III)] 
any outside edge  
for $(\Delta,{\rm Int}~J_I)$ of label $m$ 
is a terminal edge 
inward at a vertex in Int $J_I$.
\item[(IV)] 
any outside edge 
for $(\Delta,{\rm Int}~J_O)$ 
of label $m$ 
is a terminal edge 
outward at a vertex in Int $J_O$.
\end{enumerate}

\begin{lemma}\LABEL{Lemma1}
Let $\tau$ be a terminal edge of label $m$ 
in $\Delta_I$
inward at a white vertex $v_0\in J_I$, and 
$P$ a directed path 
of label $k$ 
upward maximal 
with respect to  $\Delta$ 
with a vertex sequence 
$(v_0,v_1,v_2,\cdots,v_p)$. 
Then
\begin{enumerate}
\item[{\rm (a)}] 
$v_p=P\cap J_O$, and
\item[{\rm (b)}] 
the disk $\Delta_O$ dominates 
an outside edge $e^\dagger$ 
for $\Delta$ 
of label $k$ 
outward and middle at $v_p$.
\end{enumerate}
\end{lemma}

\begin{proof}
Let 
$(e_1,e_2,\cdots,e_p)$ be 
an edge sequence 
of the path $P$.
Since the terminal edge $\tau$ is 
middle at $v_0$, so is the edge $e_1$. 

Suppose that  $v_p\not\in J_O$. 
Since $P$ is 
upward maximal with respect to $\Delta$, 
we have 
$v_p\in\partial \Delta=J_I\cup J_O$ 
by Lemma~\ref{lemMaxAtBoundaryUp}(a). 
Hence $v_p\in J_I-J_O=~{\rm Int}~J_I$. 

Since there exists an outside edge 
for $(\Delta,{\rm Int}~J_I)$ of label $k$ 
inward at $v_p$ 
by (I), 
and %
since $e_p$ is of label $k$ and 
inward at $v_p$, 
there exists 
an edge $e'$ 
of label $k$ 
outward at $v_p$. 
Thus $e'$ is an inside edge for 
$(\Delta,{\rm Int}~J_I)$ by (I). 
Now $P\cup e'$ is 
a directed path dominated by $\Delta$ 
by Lemma~\ref{LemOneWayCycle}. 
This contradicts the fact that 
the path $P$ is upward maximal 
with respect to $\Delta$.
Thus $v_p\in J_O$. 

Suppose $v_t\in P\cap J_O$
for some integer $0<t\le p$.
There are two cases.\\
Case 1. $v_t\in \partial J_O$.\\
Case 2. $v_t\in$~Int~$J_O$.

{\bf Case 1.} 
There exist 
two edges 
$e',e''$ 
of label $k$  
at $v_t$ 
such that
$e'$ is dominated by $\Delta_O$ and 
outward at $v_t$, and that
$e''$ is dominated by $\Delta_I$ and
inward at $v_t$. 
Since $e_t$ is inward at $v_t$, 
the edge $e'$ is middle at $v_t$. 
Since $e_t$ is only the inside edge for $\Delta$ 
at $v_t$,
the path $P[v_0,v_t]$ is a directed path 
upward maximal with respect to $\Delta$. 
Thus $t=p$.
Hence Statement (a) and (b) hold for this case.

{\bf Case 2.}
Let $e$ be the edge 
of label $m$ at $v_t$ with 
$e\not\subset J_O$. 
We claim that $e$ is an inside edge for $\Delta$. 
For, if $e$ is an outside edge for $\Delta$, 
then $e$ is dominated by $\Delta_O$. 
Hence $e$ is a terminal edge outward at $v_t$ 
by (IV).
Thus $e_t$ is middle at $v_t$. 
Since $e_1$ is middle at $v_0$, 
the path $P[v_0,v_t]$ is 
a dichromatic M$\&$M directed path. 
This contradicts Lemma~\ref{LemNoM-M}.
Hence $e\subset \Delta$.

Now $e\subset \Delta$ implies that 
there exist two edges $e',e''$ 
of label $k$ 
containing $v_t$ 
different from $e_t$ such that 
$e'$ is an outside edge for $\Delta$, and 
$e''$ is an inside edge for $\Delta$.  
If $e''$ is outward at $v_t$, 
then $e_t$ is middle at $v_t$. 
Again $P[v_0,v_t]$ is 
a dichromatic M\&M directed path. 
This contradicts Lemma~\ref{LemNoM-M}.
Hence $e''$ is inward at $v_t$. 
Since $e_t$ is inward at $v_t$, 
the edge $e'$ is outward and middle at $v_t$. 
Further $P[v_0,v_t]$ is a directed path 
upward maximal with respect to $\Delta$. 
Thus $t=p$.
Hence Statement (a) and (b) hold 
for this case, too.
\end{proof}

The path $P$ in Lemma~\ref{Lemma1} 
is an upward principal path 
by Lemma~\ref{lemPrincipalPath}(a), 
because $e_1$ is middle at $v_0$. 
We call $P$  
an {\it upward principal path} 
for $\tau$, 
and $e^\dagger$ a {\it corresponding edge} 
for $\tau$ with respect to $\Delta$.

By the similar way as 
the one of Lemma~\ref{Lemma1}
we can show the following lemma.

\begin{lemma}\LABEL{Lemma2}
Let $\tau$ be a terminal edge of label $m$ 
in $\Delta_O$ 
outward at a white vertex $v^\dagger\in J_O$, and 
$P$ a directed path 
of label $k$ 
downward maximal 
with respect to  $\Delta$ 
with a vertex sequence 
$(v_0,v_1,v_2,\cdots,v_p)$, 
and $v_p=v^\dagger$. 
Then
\begin{enumerate}
\item[{\rm (a)}] 
$v_0= P\cap J_I$, and
\item[{\rm (b)}] 
the disk $\Delta_I$ dominates 
an outside edge $e^*$ 
of label $k$ for $\Delta$
inward and middle at $v_0$. {\hfill {$\square$}\vspace{1.5em}}
\end{enumerate}
\end{lemma}

The path $P$ in Lemma~\ref{Lemma2} is 
a downward principal path 
by Lemma~\ref{lemPrincipalPath}(b). 
The path $P$ is called 
a {\it downward principal path}\index{downward principal} 
for $\tau$, 
and $e^*$ a {\it corresponding edge}\index{corresponding edge} 
for $\tau$ with respect to $\Delta$.


Let $P$ be a path of $\Gamma$, and 
$\tau^\dagger,\tau^*$ 
terminal edges with 
$v^\dagger=\tau^\dagger\cap P,
v^*=\tau^*\cap P$ white vertices.
If 
Int~$P[v^\dagger,v^*]$ 
does not intersect any terminal edge, 
then $\tau^\dagger,\tau^*$ are said\index{adjacent} 
to be {\it adjacent} with respect to $P$.

\begin{lemma}\LABEL{Lemma3}
\begin{enumerate}
\item[{\rm (a)}] If there exist  
two terminal edges $\tau^\dagger,\tau^*$ 
of label $m$ in $\Delta_I$, 
then two upward principal paths 
for the two terminal edges do not intersect 
each other.
\item[{\rm (b)}] If there exist  
two terminal edges $\tau',\tau''$ 
of label $m$ in $\Delta_O$, 
then two downward principal paths 
for the two terminal edges do not intersect 
each other.
\end{enumerate}
\end{lemma}

\begin{proof}
{\bf Statement (a)}. 
Since each upward principal path 
for a terminal edge 
splits the disk $\Delta$ 
by Lemma~\ref{Lemma1}, 
it is sufficient to show 
for the case that 
$\tau^\dagger$ and $\tau^*$ are adjacent 
with respect to $J_I$.

Let $v^\dagger=\tau^\dagger\cap J_I, 
v^*=\tau^*\cap J_I$, 
and 
$P^\dagger,P^*$ be 
the upward principal paths 
for $\tau^\dagger,\tau^*$.
Since $J_I[v^\dagger,v^*]$ satisfies 
Condition (i) and (ii) in 
Lemma~\ref{lemTwoPrincipal}(a),
Statement (a) follows from 
Lemma~\ref{lemTwoPrincipal}(a).

Similarly Statement (b) follows from
Lemma~\ref{lemTwoPrincipal}(b).
\end{proof}

\begin{lemma}\LABEL{Lemma4}
Each of the disks 
$\Delta_I,\Delta_O$ contains 
exactly one terminal edge of label $m$.
\end{lemma}

\begin{proof} 
Suppose that 
$\Delta_I$ contains 
at least two terminal edges of label $m$.
Let $\tau^\dagger,\tau^*$ be 
adjacent terminal edges with respect to $J_I$. 
Let $P^\dagger,P^*$ be 
the upward principal paths 
of label $k$
for $\tau^\dagger,\tau^*$ 
respectively, and
$e^\dagger,e^*$ 
the corresponding edges of label $k$ 
for $\tau^\dagger,\tau^*$ 
respectively. 
Then $P^\dagger\cap P^*=\emptyset$ 
by Lemma~\ref{Lemma3}(a).

Let $(v^\dagger_0,v^\dagger_1,\cdots,v^\dagger_p),
(v^*_0,v^*_1,\cdots,v^*_q)$ be 
vertex sequences of $P^\dagger,P^*$ 
respectively. 
Then $\tau^\dagger\cap J_I=v^\dagger_0,
\tau^*\cap J_I=v^*_0,
e^\dagger\cap J_O=v^\dagger_p,
e^*\cap J_O=v^*_q$.
Let 
\begin{enumerate}
\item[] 
$s=\max\{i~|~v^\dagger_i
\in J_I[v^\dagger_0,v^*_0]\}$
\item[] 
$t=\max\{i~|~v^*_i
\in J_I[v^\dagger_0,v^*_0]\}$
\end{enumerate}
Since the corresponding 
edges $e^\dagger, e^*$ are middle at 
$v^\dagger_p,v^*_q$ respectively 
by Lemma~\ref{Lemma1}, 
the paths 
$P^\dagger[v^\dagger_s,v^\dagger_p]
\cup e^\dagger, 
P^*[v^*_t,v^*_q]\cup e^*$ are upward principal. 
Then 
$(P^\dagger[v^\dagger_s,v^\dagger_p]
\cup e^\dagger,
P^*[v^*_t,v^*_q]\cup e^*,
J_I[v^\dagger_s,v^*_t],J_O[v^\dagger_p,v^*_q])$ 
is a principal quad of label $(k,m)$. 
This contradicts Lemma~\ref{lemNoPrincipalQuad}.
Thus 
the disk $\Delta_I$ 
contains at most one terminal edge.

Similarly we can show that 
the disk $\Delta_O$ 
contains at most one terminal edge. 
Let $v_I=e_I\cap\partial\Delta, 
v_O=e_O\cap\partial\Delta$. 
Neither $e_I$ nor $e_O$ is 
middle at $v_I$ nor $v_O$ respectively. 
Thus  
Lemma~\ref{LemTwoColorTangle}, (III), and (IV) 
assure that 
there exist at least two terminal edges 
of label $m$
in $\Delta_I\cup\Delta_O$. 
Therefore each of the disks 
$\Delta_I,\Delta_O$ 
contains exactly one terminal edge of label $m$. 
Thus Lemma~\ref{Lemma4} holds. 
\end{proof}

Throughout this section 
further we assume that 
(see Fig.~\ref{fig24}(a))
\begin{enumerate}
\item[$\bullet$] 
$v_I=e_I\cap\partial\Delta,~
v_O=e_O\cap\partial\Delta$.
\item[$\bullet$]
$\tau_I$ is 
the terminal edge of label $m$ 
in $\Delta_I$, and
$\tau_O$ is 
the terminal edge of label $m$ 
in $\Delta_O$. 
\item[$\bullet$]
$w_I=\tau_I\cap\partial\Delta$, and
$w_O=\tau_O\cap\partial\Delta$.
\item[$\bullet$] 
$P_I$ is an upward principal path 
of label $k$ for $\tau_I$.
\item[$\bullet$] 
$P_O$ is a downward principal path 
of label $k$ for $\tau_O$.
\end{enumerate}

Then we have the following by 
Lemma~\ref{Lemma4}:
\begin{enumerate}
\item[(V)] 
Neither $\Delta_I$ nor $\Delta_O$ 
dominates an outside edge 
for $\Delta$ 
of label $m$ 
intersecting 
$\partial\Delta-\{v_I,w_I,v_O,w_O\}$.
\end{enumerate}

\begin{figure}
\begin{center}
\includegraphics{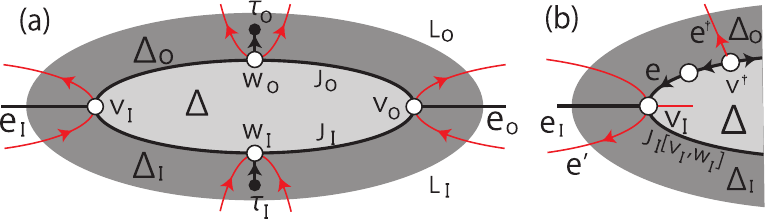}
\end{center}
\caption{ \LABEL{fig24} 
Thick arcs are of label $m$.}
\end{figure}

\begin{lemma}\LABEL{Lemma5}
There does not exist 
any outside edge for $\Delta$ 
of label $k$ 
middle at a vertex in 
$\partial\Delta-\{v_I,v_O,w_I,w_O\}$.
\end{lemma}

\begin{proof} 
Suppose that 
there exists an outside edge $e^\dagger$ 
for $(\Delta,{\rm Int}~J_O[v_I,w_O])$ 
of label $k$ 
middle at a vertex $v^\dagger$ in 
Int~$J_O[v_I,w_O]$. 
Without loss of generality 
we can assume that 
there does not exist 
any outside edge for 
$(\Delta,{\rm Int}~J_O[v_I,v^\dagger])$ 
of label $k$ 
middle at a vertex 
in Int~$J_O[v_I,v^\dagger]$. 
Hence the path $J_O[v_I,v^\dagger]$  
is a directed path, 
because  
there does not exist 
any outside edge for 
$(\Delta,{\rm Int}~J_O[v_I,w_O])$ 
of label $m$ by (V). 
Since $e^\dagger$ is outward at $v^\dagger$, 
the directed path 
$J_O[v_I,v^\dagger]$
is oriented from $v^\dagger$ to $v_I$. 
Thus an edge $e$ 
in $J_O[v_I,v^\dagger]$ with $v_I\in e$
is inward at $v_I$. 
Hence by Condition (iii) 
of the definition of a chart,
the edge $e'$ 
situated 
between $e_I$ and $J_I[v_I,w_I]$ 
around $v_I$ 
is outward at $v_I$ 
(see Fig.~\ref{fig24}(b)). 
But the edge $e'$ is dominated by 
$\Delta_I$. 
Thus $e'$ must be 
inward at $v_I$. 
This is a contradiction.
Thus there does not exist 
an outside edge 
for $(\Delta,{\rm Int}~J_O[v_I,w_O])$ 
of label $k$ 
middle at a vertex in 
Int~$J_O[v_I,w_O]$. 

Similarly we can show that 
there does not exist 
an outside edge 
of label $k$ for $\Delta$ 
middle at a vertex in 
Int~$J_O[w_O,v_O]$,  
Int~$J_I[v_O,w_I]$ nor 
Int~$J_I[w_I,v_I]$.
Thus Lemma~\ref{Lemma5} holds. 
\end{proof}

\begin{lemma}\LABEL{Lemma6}
One of $e_I,e_O$ is an I-edge 
for the disk $D$ 
and 
the other is an O-edge 
for the disk $D$.
\end{lemma}

\begin{proof}
Suppose that the both of $e_I,e_O$ 
are I-edges.
Namely $e_I$ is inward at $v_I$, and
$e_O$ is inward at $v_O$. 
By Lemma~\ref{Lemma5}, 
there does not exist 
an outside edge of label $k$ 
for $\Delta$ middle at 
a vertex in 
$\partial\Delta-\{v_I,v_O,w_I,w_O\}$. 
Since $e_I,e_O$ are inward at 
$v_I,v_O$ respectively,
the disk $\Delta_O$ does not dominate 
an edge of label $k$ middle at $v_I,v_O$. 
Further the sibling edges of $\tau_O$ are not 
middle at $w_O$. 
Hence $\Delta_O$ does not dominate 
an edge of label $k$ 
middle at a vertex in $J_O$. 
This contradicts Lemma~\ref{Lemma1} 
by considering 
the upward principal path $P_I$. 
Thus one of $e_I$ and $e_O$ is 
an O-edge. 
Similarly we can show that
one of $e_I$ and $e_O$ is an I-edge.
Hence Lemma~\ref{Lemma6} holds. 
\end{proof}

Therefore throughout this section 
furthermore 
by Lemma~\ref{Lemma4}, Lemma~\ref{Lemma5},
Lemma~\ref{Lemma6}, 
we can assume that 
(see Fig.~\ref{fig25})
\begin{enumerate}
\item[$\bullet$] 
$e_I$ is inward at $v_I$, and 
$e_O$ is outward at $v_O$.
\item[$\bullet$] 
$J_I[w_I,v_I]$ is a directed path of label $m$ 
oriented from $w_I$ to $v_I$.
\item[$\bullet$] 
$J_I[w_I,v_O]$ is 
a directed path of label $m$ 
oriented from $w_I$ to $v_O$.
\item[$\bullet$] 
$J_O[v_I,w_O]$ is a directed path of label $m$ 
oriented from $v_I$ to $w_O$.
\item[$\bullet$] 
$J_O[v_O,w_O]$ is a directed path of label $m$ 
oriented from $v_O$ to $w_O$.
\end{enumerate}

\begin{figure}
\begin{center}
\includegraphics{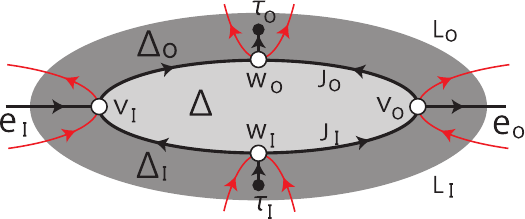}
\end{center}
\caption{ \LABEL{fig25} Thick arcs are of label $m$.}
\end{figure}

\begin{lemma}\LABEL{Lemma7}
Let $(v^*_0,v^*_1,v^*_2,\dots,v^*_p)$ 
be a vertex sequence of 
the upward principal path $P_I$ for $\tau_I$, 
and 
$(v^\dagger_0,v^\dagger_1,
v^\dagger_2,\dots,v^\dagger_q)$  
a vertex sequence of 
the downward principal path $P_O$ 
for $\tau_O$.
Then we have the following:
\begin{enumerate}
\item[{\rm (a)}]
$v^*_p=v_O,~v^\dagger_0=v_I$.
\item[{\rm (b)}] 
$P_I\cap$~{\rm Int}~$J_O=\emptyset,~
P_O\cap$~{\rm Int}~$J_I=\emptyset$. 
\end{enumerate}
\end{lemma}

\begin{proof} 
Let 
$(e^*_1,e^*_2,\dots,e^*_p),
(e^\dagger_1,e^\dagger_2,\dots,
e^\dagger_q)$
be edge sequences of $P_I, P_O$ respectively.

{\bf Statement (a)}.
By Lemma~\ref{Lemma1} 
and Lemma~\ref{Lemma5}, we have 
$v^*_p=v_I,~v^*_p=w_O,$ or $v^*_p=v_O$.
If $v^*_p=w_O$, 
then the path $P_I$ is 
a dichromatic M$\&$M directed path. 
This contradicts Lemma~\ref{LemNoM-M}. 
Since $e^*_p$ is inward at $v^*_p$, 
we have $v^*_p=v_O$.
Similarly we can show that 
$v^\dagger_0=v_I$.
Hence Statement (a) holds.

{\bf Statement (b)}. 
By Lemma~\ref{Lemma1}(a), 
$v^*_p=P_I\cap J_O$. 
Hence $v^*_p=v_O\in \partial J_O$ implies  
$P_I\cap$~Int~$J_O=\emptyset$.
Similarly we can show that 
$P_O\cap$~Int~$J_I=\emptyset$. 
Hence Statement (b) holds.
\end{proof}

Let $\Gamma$ be a chart in a disk $D^2$, 
and $m,k$ be labels of $\Gamma$ 
with $|m-k|=1$.
A disk $E$ with 
$\Gamma\cap E\subset\Gamma_m\cup\Gamma_k$
is {\it bigonal}\index{bigonal} 
if $\partial E$ is a union of 
two paths $P,L$ such that 
(see Fig.~\ref{fig26}) 
\begin{enumerate}
\item[(i)] 
there exist two white vertices $v^*,v^{**}$ 
on $\partial E$
such that
\begin{enumerate} 
\item[(a)] $P$ is a directed path of label $k$ 
oriented from $v^*$ to $v^{**}$, and 
\item[(b)] 
$L$ is a directed path of label $m$ 
oriented from $v^*$ to $v^{**}$,
\end{enumerate}
\item[(ii)] 
the disk $E$ does not dominate 
any edge of label $m$ 
at $v^*$ nor $v^{**}$,
\item[(iii)]
there does not exist any outside edge for 
$(E,{\rm Int}~L)$ of label $m$,

\item[(iv)]
the path $P$ is upward principal 
or downward principal,
and 
if $P$ is upward principal 
(resp. downward principal), then
any outside edge for $(E,{\rm Int}~L)$ 
of label $k$  
is inward (resp. outward) at 
a vertex in Int~$L$.
\end{enumerate}


\begin{figure}
\begin{center}
\includegraphics{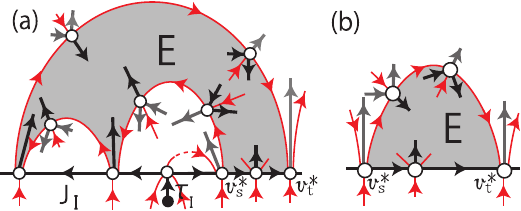}
\end{center}
\caption{ \LABEL{fig26} Bigonal disks.}
\end{figure}


\begin{lemma}
\LABEL{lemBigonal}
Any bigonal disk 
in a minimal chart  
is a bigon.
\end{lemma}

\begin{proof}
Let $E$ be a bigonal disk 
in a minimal chart $\Gamma$.
We use all the notations 
in the definition of a bigonal disk.

We only show the case that 
$P$ is an upward principal directed path 
of label $k$.
Then we have that
\begin{enumerate}
\item[(1)] 
if an edge $e$ of label $k$ 
intersects Int~$P$ but $e\not\subset P$, 
then the edge $e$ is inward at 
a vertex of Int~$P$.
\end{enumerate}
Thus 
by considering Statement (1),
Condition (i),(iii) and (iv) for a bigonal disk,
among the edges of label $m$ 
dominated by $E$, for each vertex $v$ in 
$\partial E-(P\cap L)$ 
(see Fig.~\ref{fig26})
\begin{enumerate}
\item[(2)] 
\begin{enumerate}
\item[(a)] there exists
exactly one edge of label $m$
outward at $v$, or
\item[(b)] there exist
exactly two edges of label $m$; 
one is inward at $v$ and 
the other outward at $v$.
\end{enumerate}
\end{enumerate}
Suppose that
$\partial E-(P\cap L)$ 
contains a white vertex $v'$. 
Then
the disk $E$ dominates 
an edge $e'$ of label $m$ 
outward at $v'$.
Let $P'$ be 
a directed path of label $m$ 
starting from $e'$
upward maximal 
with respect to $E$ with 
a vertex sequence 
$(v'_0,v'_1,v'_2,\cdots,v'_t)$ 
here $v'_0=v'$.
By Lemma~\ref{LemOneWayCycle}, 
we have $v'_t\neq v'_0$. 
Further 
$v'_t\in\partial E$
by Lemma~\ref{lemMaxAtBoundaryUp}.
Thus Condition (ii) of 
a bigonal disk implies that 
$v'_t$ is a white vertex of Statement (2)(b). 
Since there exists an edge of label $m$ 
outward at $v'_t$, 
Lemma~\ref{LemOneWayCycle} assures that 
$P'$ is not upward maximal 
with respect to $E$.
This is a contradiction. 
Thus the white vertices of $\partial E$ 
are $\{v^*,v^{**}\}=P\cap L$. 
Namely $E$ is a bigon 
by Condition (ii) of 
a bigonal disk. 
Hence Lemma~\ref{lemBigonal} holds.
\end{proof}

For the sets $X,Y$, define
$X\Delta Y=(X-Y)\cup (Y-X)$.

\begin{lemma}\LABEL{Lemma8}~\\
{\rm (a)} 
$P_I\cap {\rm Int}~J_I[v_I,w_I]=\emptyset$ and
$P_O\cap {\rm Int}~J_O[w_O,v_O]=\emptyset$.\\
{\rm (b)} 
Neither 
$P_I\Delta J_I[w_I,v_O]$ nor 
$P_O\Delta J_O[v_I,w_O]$ 
contains a white vertex.\\
{\rm (c)}
Neither {\rm Int}~$J_I[v_I,w_I]$ nor 
{\rm Int}~$J_O[w_O,v_O]$ 
contains a white vertex.
\end{lemma}

\begin{proof} 
Let $(v^*_0,v^*_1,\cdots,v^*_p)$ be 
a vertex sequence of $P_I$. 
Then $v^*_0=w_I,v^*_p=v_O$ 
by Lemma~\ref{Lemma7}(a).

{\bf Statement (a).}
We only show that 
$P_I\cap {\rm Int}~J_I[v_I,w_I]=\emptyset$. 
Suppose that there exists a vertex $v$ in 
$P_I\cap$ Int $J_I[v_I,w_I]$. 
Then $v=v^*_j$ for some integer $j~(0<j<p)$ 
by Lemma~\ref{Lemma7}(a). 
let
\begin{enumerate}
\item[] 
$s=\max\{i~|~i<j,~v^*_i\in
J_I[w_I,v_O]\}$ and\\ 
$t=\min\{i~|~j<i,~v^*_i\in   
J_I[w_I,v_O]\}$.
\end{enumerate}
Then the disk $E$ bounded by 
$P_I[v^*_s,v^*_t]\cup  
J_I[v^*_s,v^*_t]$ 
is a bigonal disk in $\Delta$
(see Fig.~\ref{fig26}(a)).  
Also $v\in E$. 
On the other hand, 
the disk $E$ is a bigon with 
by Lemma~\ref{lemBigonal}.
Since $v^*_s,v^*_t\in J_I[w_I,v_O]$ are 
the only vertices of the bigon $E$, 
namely $v\not\in E$. 
This is a contradiction. 
Hence $P_I\cap 
{\rm Int}~J_I[v_I,w_I]=\emptyset$.

{\bf Statement (b)}. 
Now $P_I\cap J_I\subset J_I[w_I,v_O]$.
Let $L$ be 
the closure of 
the connected component of 
$P_I\Delta J_I[w_I,v_O]$. 
Then $\partial L$ consists of
two vertices $v^*_s,v^*_t$ 
for some integers $0\le s<t\le p$.
Let $E$ be the disk bounded by 
$P_I[v^*_s,v^*_t]\cup J_I[v^*_s,v^*_t]$ 
(see Fig.~\ref{fig26}(b)). 
Then $E$ is a bigonal disk in $\Delta$. 
Hence $E$ is a bigon by Lemma~\ref{lemBigonal}.
Thus neither Int~$P_I[v^*_s,v^*_t]$ 
nor Int~$J_I[v^*_s,v^*_t]$ 
contains a white vertex. 
Hence Statement (b) holds.

{\bf Statement (c).} 
Suppose that 
there exists a white vertex $v$ in
Int~$J_I[w_I,v_I]$. 
Then $\Delta$ contains 
an edge $e$ of label $k$
outward at $v$. 
Let $P'$ be a directed path of label $k$ 
starting from $e$ 
upward maximal 
with respect to $\Delta$ 
with 
a vertex sequence 
$(v'_0,v'_1,v'_2,\cdots,v'_t)$. 
Then $v'_t\in \partial\Delta$ 
by Lemma~\ref{lemMaxAtBoundaryUp}(a). 
Since $\{v_I,v_O,w_I,w_O\}\subset  P_I\cup  P_O$
by Lemma~\ref{Lemma7}(a),
we have $v'_t\not\in \{v_I,v_O,w_I,w_O\}$. 
Thus by (V) 
there exists exactly 
one outside edge $e^\dagger$ 
for $\Delta$ at $v'_t$, 
which is of label $k$. 
Further $e^\dagger$ is
not middle at $v'_t$ by Lemma~\ref{Lemma5}. 
Hence one of the two inside edges of label $k$ 
at $v'_t$ is outward at $v'_t$, 
say $e^*$. 
Then $P'\cup e^*$ is 
a directed path 
of label $k$
starting from $e$ 
dominated by $\Delta$. 
This contradicts the fact that 
$P'$ is upward maximal 
with respect to $\Delta$.
Hence there does not exist 
any white vertex in Int~$J_I[w_I,v_I]$. 
Similarly we can show that 
there does not exist 
any white vertex in Int~$J_O[v_O,w_O]$. 
Hence Statement (c) holds. 
\end{proof}

\begin{figure}
\begin{center}
\includegraphics{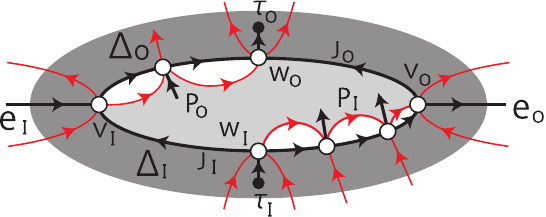}
\end{center}
\caption{ \LABEL{fig27} 
The disk $E$ is light gray. 
Thick arcs are of label $m$.}
\end{figure}

{\bf Proof of Theorem~\ref{StIITheorem4}.}
If there exists a white vertex $v$ in 
Int~$\Delta$, 
there exists a directed path of label $k$ 
starting from an edge  
outward at $v$ 
upward maximal with respect $\Delta$.
We get the same contradiction 
as the one of Lemma~\ref{Lemma8}(c). 
Thus  
there does not exist 
any white vertex in Int~$\Delta$. 
(Now the tangle is like 
the one in Fig.~\ref{fig27}.) 
Let $s$ be 
the number of white vertices 
in Int~$J_I[w_I,v_O]$, and 
$t$
the number of white vertices 
in Int~$J_O[v_I,w_O]$. 
Since there does not exist 
any white vertex in Int~$\Delta$, 
we have $s=t$.  
Thus $(\Gamma\cap D,D)$ is
a Type-$I_{s+1}$ elementary IO-tangle. 
This proves Theorem~\ref{StIITheorem4}. {\hfill {$\square$}\vspace{1.5em}}



\section{Indices}
\LABEL{s:Indices}

We define indices of a simple IO-tangle and 
a net-tangle.

Let $\Gamma$ be a minimal chart, and 
$m$ a label of the chart. 
If 
a terminal edge is inward
at its black vertex,
then the edge 
is called an {\it I-terminal edge},\index{I-terminal edge}
otherwise 
the edge is called 
an {\it O-terminal edge}.\index{O-terminal edge} 
Let $(\Gamma\cap D,D)$ be 
a simple IO-tangle 
of label $m$ 
with a boundary IO-arc pair $(L_I,L_O)$. 
Let $k$ be a label of $\Gamma$ with 
$|m-k|=1$ and 
$\Gamma\cap D\subset\Gamma_m\cup\Gamma_k$. 
Let $\sigma_1,\sigma_2,\cdots,\sigma_s$
be all the O-terminal edges 
of label $m$ in $D$
and 
$\widetilde\sigma_1,\widetilde\sigma_2,
\cdots,\widetilde\sigma_t$
all the I-terminal edges of 
label $m$ in $D$.
For each terminal edge $\tau$ in $D$,
let $e^*,e^{**}$ be the sibling edges 
of the terminal edge.
The union $e^*\cup e^{**}$ splits
the disk $D$ into two disks.
Let $\Delta(\tau)$ be
the one of the two disks containing
the terminal edge.
Then we can show 
$\Delta(\sigma_i)
\cap\partial D\subset L_I,
\Gamma\cap({\rm Int}~\Delta(\sigma_i)
-\sigma_i)
=\emptyset~~(i=1,2,\cdots,s)$ and 
$\Delta(\widetilde\sigma_j)
\cap\partial D\subset L_O,
\Gamma\cap({\rm Int}~\Delta(\widetilde\sigma_j)
-\widetilde\sigma_j)
=\emptyset~~(j=1,2,\cdots,t)$ 
(cf. \cite[Lemma~7.6]{StI}).
Set\\
\ \ \ $D^\dagger=
Cl(D-\cup_{i=1}^s\Delta(\sigma_i)
-\cup_{j=1}^t\Delta(\widetilde\sigma_j))$ 
(see Fig.~\ref{fig28}(b)).\\
Let
$A_0,A_1,\cdots,A_s$ be the connected components of
$\partial D^\dagger\cap L_I$
situated counterclockwise on $\partial D^\dagger$
in this order, and
$B_0,B_1,\cdots,B_t$ the connected components of
$\partial D^\dagger\cap L_O$
situated clockwise on $\partial D^\dagger$
in this order.
Let $X_I$ be the union of 
all the I-edges 
for $D$ of label $k$, and
$X_O$ the union of 
all the O-edges 
for $D$ of label $k$. 
For each $i=0,1,\cdots,s$ and 
$j=0,1,\cdots,t$, 
let\\
\ \ \ \ $a_i=|A_i\cap X_I|$, and
$b_j=|B_j\cap X_O|$.\\
Then $(s+1)$-tuple\index{${\bf a}(\Gamma,D)$}
${\bf a}(\Gamma,D)=
(a_0,a_1,\cdots,a_s)$ and\index{${\bf b}(\Gamma,D)$}
$(t+1)$-tuple ${\bf b}(\Gamma,D)=
(b_0,b_1,\cdots,b_t)$
are called 
the {\it I-index} and the {\it O-index} 
of the simple IO-tangle respectively.\index{Index$(\Gamma,D)$}
The pair 
$({\bf a}(\Gamma,D),{\bf b}(\Gamma,D))$ 
is called the {\it index} of the 
simple IO-tangle, 
and denoted by Index$(\Gamma,D)$. 
For the simple IO-tangle shown 
in Fig.~\ref{fig28}(a),
the I-index is $(4,3,1)$ and 
the O-index is $(2,5)$.

\begin{figure}
\begin{center}
\includegraphics{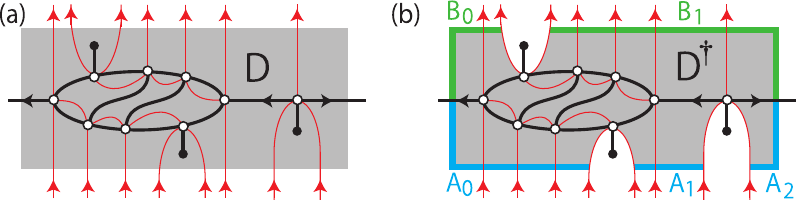}
\end{center}
\caption{ \LABEL{fig28}
(a) The thick edges are of label $m$.
(b) The thick edges are of label $m$, and
the thin edges are of label $k$.}
\end{figure}

Let $\Gamma$ be a chart.
A tangle $(\Gamma\cap D,D)$ 
is called a {\it net-tangle}\index{net-tangle} 
provided that
\begin{enumerate}
\item[(i)]
the disk $D$ contains no crossing, hoop, nor free edge
but a white vertex,
\item[(ii)]
there exist two labels $\alpha,\beta$ 
with $\alpha<\beta$ and
$\Gamma\cap D
\subset
\cup_{i=\alpha}^{\beta}\Gamma_{i}$,
and
\item[(iii)]
there exist two arcs 
$L_\alpha,L_\beta$ on $\partial D$ 
with $L_\alpha\cap L_\beta$ two points 
such that
\begin{enumerate}
\item[(a)] 
$\Gamma\cap\partial D=
(\Gamma_\alpha\cap {\rm Int}~L_\alpha)
\cup 
(\Gamma_\beta\cap {\rm Int}~L_\beta)$,
\item[(b)]
all the edges intersecting $L_\alpha$ 
are I-edges of label $\alpha$ or\\
all the edges intersecting $L_\alpha$ 
are O-edges of label $\alpha$,
and
\item[(c)]
all the edges intersecting $L_\beta$ 
are O-edges of label $\beta$ or\\
all the edges intersecting $L_\beta$ 
are I-edges of label $\beta$.
\end{enumerate}
\end{enumerate}
The pairs $(\alpha,\beta)$ and 
$(L_\alpha,L_\beta)$ are called 
a {\it label pair}\index{label pair} 
and 
a {\it boundary arc pair} of the net-tangle\index{boundary arc pair}
respectively.
If all the edges of labe $\alpha$ 
intersecting $L_\alpha$
are I-edges (resp. O-edges), 
and
if all the edges of label $\beta$ 
intersecting $L_\beta$
are 
O-edges (resp. I-edges),
then the net-tangle is said to be 
{\it upward} (resp. {\it downward})\index{upward net-tangle}\index{downward net-tangle} 
(see Fig.~\ref{fig29}).
An upward or downward 
net-tangle with 
a label pair $(\alpha,\alpha+1)$\index{N-tangle} 
is called an {\it N-tangle}.

\begin{figure}
\begin{center}
\includegraphics{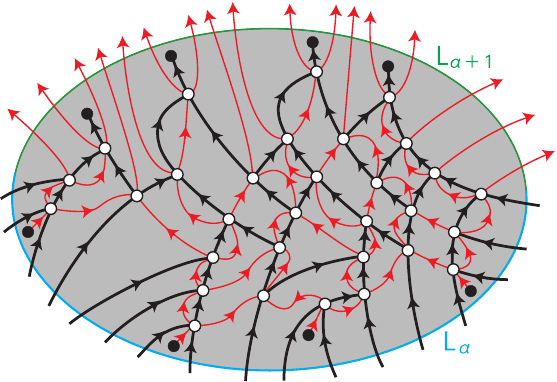}
\end{center}
\caption{ \LABEL{fig29}
An example of an upward net-tangle with 
a label pair $(\alpha,\alpha+1)$.
The thick edges are of label $\alpha$, 
and
the thin edges are of label $\alpha+1$.
}
\end{figure}

To define an index of an N-tangle, 
we need the following theorem.

\begin{theorem}
{\rm (\cite[Theorem 1.1]{StI})}
\LABEL{StITheorem1}
Let $\Gamma$ be a minimal chart, and 
$(\Gamma\cap D,D)$ a net-tangle
with a label pair $(\alpha,\alpha+1)$. 
Then we have the following:
\begin{enumerate}
\item[{\rm (a)}] 
The tangle is an N-tangle.
\item[{\rm (b)}] 
The number of the edges of label $\alpha$ 
intersecting $\partial D$
is equal to
the number of the edges of label $\alpha+1$ 
intersecting $\partial D$.
\item[{\rm (c)}] 
There exists a terminal edge in $D$.
\item[{\rm (d)}]
The number of terminal edges of 
label $\alpha$ in $D$ 
is equal to 
the number of terminal edges of 
label $\alpha+1$ in $D$.
\item[{\rm (e)}]
If the tangle is upward $($resp. downward$)$, 
then 
all the terminal edges of label $\alpha$ in $D$ 
are I-terminal $($resp. O-terminal$)$ edges
and 
all the terminal edges of label $\alpha+1$ in $D$ 
are O-terminal $($resp. I-terminal$)$ edges
\end{enumerate}
\end{theorem}

Now we can define an index of the N-tangle
as follows. 

Let $\Gamma$ be a minimal chart.
Let $(\Gamma\cap D,D)$ be an N-tangle 
with 
a label pair $(\alpha,\alpha+1)$ and a boundary arc pair 
$(L_{\alpha},L_{\alpha+1})$.
Let $s_I$ be the label of the I-edges, and 
 $s_O$ the label of the O-edges.
Further, we denote by $L_I$ (resp. $L_O$)
the one of the arcs 
$L_{\alpha},L_{\alpha+1}$
which intersects I-edges (resp. O-edges).
Considering 
Theorem~\ref{StITheorem1}(d),
let\\ 
$k=$the number of O-terminal edges 
in $D$\\
~~=the number of I-terminal edges 
in $D$.\\
Let $\sigma_1,\sigma_2,\cdots,\sigma_k$
be all the O-terminal edges in $D$
and 
$\widetilde\sigma_1,\widetilde\sigma_2,
\cdots,\widetilde\sigma_k$
all the I-terminal edges in $D$.

For each terminal edge $\tau$ in $D$, 
we define $\Delta(\tau)$ 
in the same way as the one for an IO-tangle.

\begin{lemma}{\rm (\cite[Lemma~7.6]{StI})}
\label{lemDelta}
Let $\Gamma$ be a minimal chart, and 
$(\Gamma\cap D,D)$ an N-tangle 
with a label pair $(\alpha,\alpha+1)$ 
and a boundary arc pair 
$(L_\alpha,L_{\alpha+1})$.
Let $\tau$ be a terminal edge in $D$. 
Then we have the following.
\begin{enumerate}
\item[{\rm (a)}] 
If the label of $\tau$ is $\alpha+1$, 
then $\Delta(\tau)\cap\partial D\subset L_\alpha$
otherwise 
$\Delta(\tau)\cap\partial D\subset L_{\alpha+1}$.
\item[{\rm (b)}] 
$\Gamma\cap ({\rm Int}~\Delta(\tau)-\tau)
=\emptyset$.
\end{enumerate}
\end{lemma}
Set\\
\ \ \ $D^\dagger=
Cl(D-\cup_{i=1}^k\Delta(\sigma_i)
-\cup_{j=1}^k\Delta(\widetilde\sigma_j))$ 
(see Fig.~\ref{fig30}(b)).\\
Let
$A_0,A_1,\cdots,A_k$ be 
the connected components of
$\partial D^\dagger\cap L_I$
situated counterclockwise 
on $\partial D^\dagger$
in this order, and
$B_0,B_1,\cdots,B_k$ be 
the connected components of
$\partial D^\dagger\cap L_O$
situated clockwise 
on $\partial D^\dagger$
in this order.
Let $X_I$ be the union of 
all the I-edges for $D$, and
$X_O$ the union 
of all the O-edges for $D$.
For each $i=0,1,\cdots,k$, let\\
\ \ \ \ $a_i=|A_i\cap X_I|$, and
$b_i=|B_i\cap X_O|$.\\
Then $(k+1)$-tuples
${\bf a}(\Gamma,D)=(a_0,a_1,\cdots,a_k)$ and
${\bf b}(\Gamma,D)=(b_0,b_1,\cdots,b_k)$
are called 
the I-{\it index} and O-{\it index}\index{${\bf a}(\Gamma,D)$}\index{${\bf b}(\Gamma,D)$}
of the N-tangle respectively.\index{Index$(\Gamma,D)$}
The pair 
$({\bf a}(\Gamma,D),{\bf b}(\Gamma,D))$
is called the {\it index} of the N-tangle, 
and denoted by Index$(\Gamma,D)$. 
According to Theorem~\ref{StITheorem1}(b),
we have $\sum_{i=0}^ka_i=\sum_{i=0}^kb_i$.

For the N-tangle as the one shown 
in Fig.~\ref{fig30}(a),
the I-index is $(3,3,3,1)$ and 
the O-index is $(2,3,4,1)$.

Let ${\bf x}=(x_1,x_2,\cdots,x_k)$ 
be an I-index or O-index of a tangle, 
then\index{$||\bf x||$}
the sum of components, $\sum_{i=1}^kx_i$, 
is called {\it index sum}, and denoted by $||\bf x||$. 

\begin{figure}
\begin{center}
\includegraphics{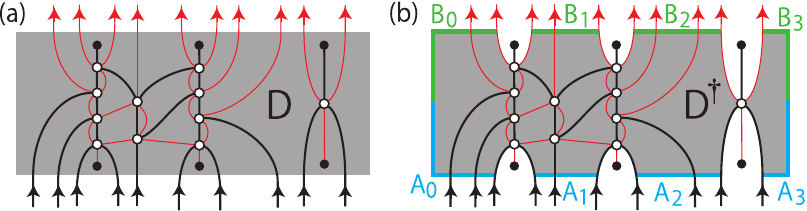}
\end{center}
\caption{\LABEL{fig30}
The thick edges are of label $s_I$, and
the thin edges are of label $s_O$.}
\end{figure}

\section{Normal forms}
\LABEL{s:NormalForm}

To define a normal form 
for 2-crossing minimal charts, 
we need 
one more theorem 
proved in \cite{StI}.

There exists a special C-move 
called a C-I-M2 move 
(see Fig.~\ref{fig05} in Section~\ref{s:Prel}).
Let $\Gamma$ and $\widetilde\Gamma$ be 
minimal charts,
and $(\Gamma\cap D,D)$ a net-tangle with 
a label pair $(\alpha,\beta)$.
Then the chart $\widetilde\Gamma$ is said to be\index{M2-related to $\Gamma$} 
M2-{\it related to $\Gamma$ with respect to $D$}
provided that
\begin{enumerate}
\item[(i)]
$\Gamma\cap D^c=\widetilde\Gamma\cap D^c$, and
\item[(ii)]
the chart $\widetilde\Gamma$ is obtained from
the chart $\Gamma$ by a finite sequence of 
C-I-M2 moves in $D$ each of which modifies 
two edges of label $m$ with $\alpha<m<\beta$.
\end{enumerate}

\begin{theorem}
{\rm (\cite[Theorem 1.2]{StI})}
\LABEL{StITheorem2}
Let $\Gamma$ be a minimal chart, and 
$(\Gamma\cap D,D)$ a net-tangle
with a label pair $(\alpha,\beta)$. 
Then 
there exists a minimal chart $\widetilde\Gamma$
M2-related to $\Gamma$ with respect to $D$ such that
there exist N-tangles 
$(\widetilde\Gamma\cap D_{\alpha},D_{\alpha}),
(\widetilde\Gamma\cap D_{\alpha+1},D_{\alpha+1}),
\cdots,
(\widetilde\Gamma\cap D_{\beta-1},D_{\beta-1})$
equipped with
\begin{enumerate}
\item[{\rm (a)}] for each 
$i=\alpha,\alpha+1,\cdots,\beta-1$,
the tangle 
$(\widetilde\Gamma\cap D_{i},D_{i})$ 
is an N-tangle with the label pair $(i,i+1)$,
\item[{\rm (b)}] $D=\cup_{i=\alpha}^{\beta-1}D_i$, 
\item[{\rm (c)}] 
for each $i=\alpha,\alpha+1,\cdots,\beta-2$,
the intersection $D_i\cap D_{i+1}$ 
is a proper arc of $D$,
\item[{\rm (d)}]
all the N-tangles are 
upward or downward
simultaneously.
\end{enumerate}
\end{theorem}

\begin{remark} \LABEL{remWhite2}
{\rm 
Let $\Gamma$ be a minimal chart, 
$m,k$ labels of $\Gamma$, and 
$(\Gamma\cap D,D)$ 
a nontrivial IO-tangle of label $m$ with 
$\Gamma\cap D\subset\Gamma_m\cup\Gamma_k$.
Considering Fig.~\ref{fig14},
\begin{enumerate}
\item[(i)] 
the IO-tangle $(\Gamma\cap D,D)$ 
contains at most one white vertex\\
if and only if 
there exist at most two I-edges for $D$ 
of label $k$ and 
at most two O-edges for $D$ 
of label $k$. 
\end{enumerate}
Therefore 
\begin{enumerate}
\item[(ii)]
the IO-tangle $(\Gamma\cap D,D)$ 
contains at least two white vertices\\
if and only if 
there exist at least three I-edges for $D$ 
of label $k$ or 
at least three O-edges for $D$ 
of label $k$. 
\end{enumerate}
}\end{remark}

Throughout this section,
$\Gamma$ is a $2$-crossing minimal chart 
in a disk $D^2$ different from the chart 
in Fig.~\ref{fig31}(a).
Set
$\alpha=\alpha(\Gamma)$, and 
$\beta=\beta(\Gamma)$.
Then by Lemma~\ref{StIILemma1},
we can assume that 
\begin{enumerate}
\item[$\bullet$]
there exist two cycles 
$C_\alpha,C_\beta$ 
with 
$C_\alpha\subset\Gamma_\alpha,C_\beta
\subset\Gamma_\beta$ such that 
$C_\alpha\cap C_\beta$ consists of 
the two crossings,
\item[$\bullet$]
there exists an annulus $A$ 
with $A\cap \partial D^2=\emptyset$
such that $A$ contains all the white vertices of $\Gamma$ 
but does not intersect hoops nor free edges,
\item[$\bullet$]
each connected component of $Cl(D^2-A)$ contains a crossing,
\item[$\bullet$]
 $\Gamma\cap \partial A=(C_\alpha\cup C_\beta)\cap \partial A$, 
$\Gamma\cap \partial A$ consists of eight points,
\item[$\bullet$] 
$\Gamma_\alpha\cap A$ 
consists of 
two connected components 
$X_1,X_3$ 
separated by $C_\beta$, 
\item[$\bullet$]
$\Gamma_\beta\cap A$ 
consists of 
two connected components 
$X_2,X_4$ 
separated by $C_\alpha$. 
\end{enumerate}
\begin{figure}
\begin{center}
\includegraphics{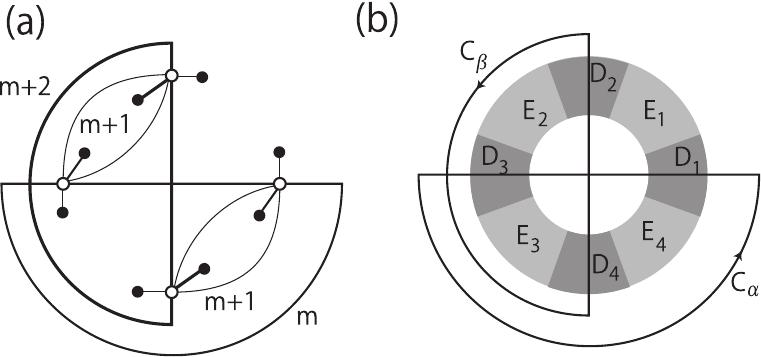}
\end{center}
\caption{ \LABEL{fig31}
(a) the thin edges are of label $m+1$.
(b) $A=D_1\cup E_1\cup D_2\cup E_2\cup
D_3\cup E_3\cup D_4\cup E_4$.}
\end{figure}
For each $i=1,2,3,4$, 
let $D_i$ be a regular neighbourhood 
of the SC-closure $SC(X_i)$ in the annulus $A$.
By Boundary Condition Lemma 
(Lemma~\ref{BoundaryConditionLemma}),
\begin{enumerate}
\item[$\bullet$] $\Gamma\cap D_i\subset\Gamma_\alpha\cup\Gamma_{\alpha+1}$ for $i=1,3$,
\item[$\bullet$] $\Gamma\cap D_i\subset\Gamma_\beta\cup\Gamma_{\beta-1}$ for $i=2,4$.
\end{enumerate}
By Theorem~\ref{StIITheorem2} and Remark~\ref{remWhite1},
\begin{enumerate}
\item[$\bullet$]
the tangle $(\Gamma\cap D_i,D_i)$ is an IO-tangle for
$i=1,2,3,4$.
\end{enumerate}
Without loss of generality
we can assume that
(see Fig.~\ref{fig31}(b))
\begin{enumerate}
\item[$\bullet$]
$Cl(A-(D_1\cup D_2\cup D_3\cup D_4))$ consists of
four disks, say
$E_1,E_2,E_3,E_4$,
\item[$\bullet$]
$D_1,E_1,D_2,E_2,
D_3,E_3,D_4,E_4$
are situated on the annulus $A$ 
in this order,
\item[$\bullet$]
an O-edge for $D_1$ is an I-edge for $E_1$.
\end{enumerate}

\begin{lemma}
\LABEL{FoundamentalTwo}
For each fundamental tangle 
$(\Gamma\cap D_i,D_i)$ $(i=1,2,3,4)$,
the disk $D_i$ contains 
at least two white vertices.
\end{lemma}

\begin{proof}
We need two claims.

{\bf Claim 1}. 
Each of $\Gamma_\alpha$ and $\Gamma_\beta$
contains at least two white vertices.

For, by Lemma~\ref{Lemma3-2},
$\Gamma_\alpha$ contains 
at least one white vertex $v$.
If $\Gamma_\alpha$ contains 
exactly one white vertex,
then Main$(\Gamma_\alpha)$ 
consists of a loop $C_\alpha$ and 
a terminal edge 
not middle at $v$. 
Thus we can eliminate 
the white vertex 
by a C-III move. 
This contradicts the fact that
the chart is minimal. 
Thus $\Gamma_\alpha$ contains 
at least two white vertices.
Similarly 
we can show that $\Gamma_\beta$ 
contains at least two white vertices.
Hence Claim 1 holds.

{\bf Claim 2}.
One of $D_1,D_2,D_3,D_4$ 
contains at least two white vertices.

{\it Proof of Claim 2}. 
Suppose that
each of $D_1,D_2,D_3,D_4$ 
contains at most one white vertex.
Since each of $\Gamma_\alpha$ and $\Gamma_\beta$
contains at least two white vertices,
each of $D_1,D_2,D_3,D_4$  
contains exactly one white vertex.
Thus for each $i=1,2,3,4$ 
the tangle $(\Gamma\cap D_i,D_i)$ is 
one of the four tangles as shown in
Fig.~\ref{fig14}
by Remark~\ref{remWhite1}.
Without loss of generality
we can assume that
\begin{enumerate}
\item[] $\Gamma\cap(D_1\cap E_1)$ 
consists of at most one point.
\end{enumerate}
There are two cases: Case 1. $\beta-\alpha>2$ or 
Case 2. $\beta-\alpha=2$.

{\bf Case 1}. 
Suppose $\beta-\alpha>2$.
If $E_1$ contains at least one white vertex,
then $(\Gamma\cap E_1,E_1)$ is 
an NS-tangle of label $\beta-1$.
This contradicts Lemma~\ref{LemNS-Tangle}.
Thus $E_1$ does not contain any white vertex.
Since $\beta-\alpha>2$, 
we have $\Gamma\cap E_1=\emptyset$. 
Hence $(\Gamma\cap D_1,D_1)$
is a tangle as shown in Fig.~\ref{fig14}(c),(d) 
and $E_1\cap D_1=E_1\cap D_2=\emptyset$. 
Thus $(\Gamma\cap D_2,D_2)$
is a tangle as shown in Fig.~\ref{fig14}(c),(d) 
and
$D_2\cap E_2$ consists of two points. 
Hence $D_3\cap E_2\neq\emptyset$, 
and 
$D_3\cap E_3$ consists of at most one point.
Thus again $E_3$ does not contain any white vertex. So on...
finally $\Gamma$ is the one like 
in Fig.~\ref{fig32}(a). 
Since $\beta-\alpha>2$,
by C-II moves, a C-I-M2 move and a C-III move,
we obtain the chart as shown in 
Fig.~\ref{fig32}(e). 
The resulting chart is not a minimal chart. 
This contradicts the fact that 
the chart $\Gamma$ is not minimal.

\begin{figure}
\begin{center}
\includegraphics{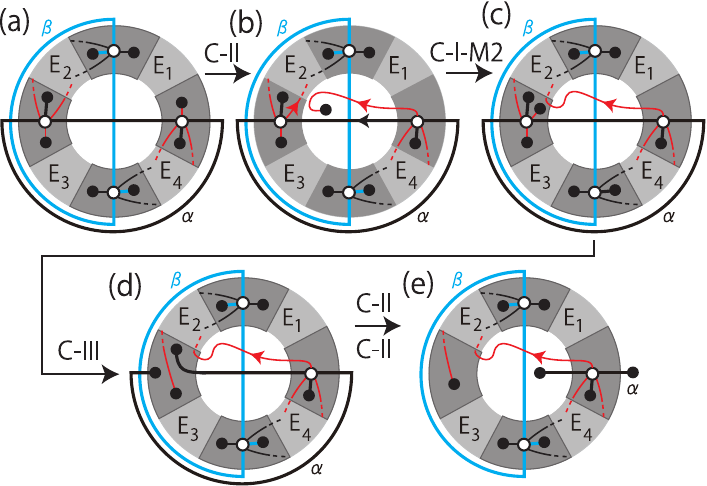}
\end{center}
\caption{ \LABEL{fig32} If $|\alpha-\beta|>2$, then any 2-crossing chart with four white vertices  is not minimal.}
\end{figure}


{\bf Case 2}. 
Suppose $\beta-\alpha=2$. 
By the same way we can show that 
$E_1$ does not contain any white vertex. 
If $\Gamma\cap(D_1\cap E_1)$ 
consists of one point, 
then a part of $\Gamma\cap A$ is 
the one as shown in Fig.~\ref{fig-33}(a). 
Looking at the vertex in $D_3$, 
we will find 
incorrect orientation 
of edges of label $\alpha$ 
around the vertex (see Fig.~\ref{fig-33}(b)). 
Thus $\Gamma\cap(D_1\cap E_1)=\emptyset$.
Hence we can show that the chart $\Gamma$ is 
the one shown in 
Fig.~\ref{fig31}(a).
But our chart $\Gamma$ is assumed to be 
different from 
the one shown in 
Fig.~\ref{fig31}(a).

Therefore one of $D_1,D_2,D_3,D_4$ 
contains at least two white vertices.
Thus Claim 2 holds.


Without loss of generality
we can assume that 
$D_1$ contains at least two white vertices.
Then $(\Gamma\cap D_1,D_1)$ is 
a simple IO-tangle
by Theorem~\ref{StIITheorem2}. 
Hence 
there exist
at least three O-edges 
of label $\alpha+1$ 
for $D_1$ by 
Remark~\ref{remWhite2}(ii). 
Namely there exist 
at least three I-edges 
for $E_1$ of label $\alpha+1$.
Since $(\Gamma\cap D_1,D_1)$ is 
an IO-tangle of label $\alpha$
and since $(\Gamma\cap D_2,D_2)$ is 
an IO-tangle of label $\beta$,
the tangle $(\Gamma\cap E_1,E_1)$ is 
a net-tangle with 
a label pair $(\alpha+1,\beta-1)$
by Boundary Condition Lemma
(Lemma~\ref{BoundaryConditionLemma}).
By Theorem~\ref{StITheorem1}
and Theorem~\ref{StITheorem2}, 
there are 
at least three O-edges 
of label $\beta-1$ 
for $E_1$ 
(here possibly $\alpha+1=\beta-1$,  
in this case, $E_1$ contains 
just parallel arcs). 
Namely there are 
at least three I-edges 
of label $\beta-1$ 
for $D_2$. 
Hence by Remark~\ref{remWhite2}(ii), 
the disk $D_2$ contains 
at least two white vertices.

By the similar way, 
we can show that 
each of $D_1,D_2,D_3,D_4$ 
contains at least two white vertices.
Thus Lemma~\ref{FoundamentalTwo} holds.
\end{proof}

\begin{figure}
\begin{center}
\includegraphics{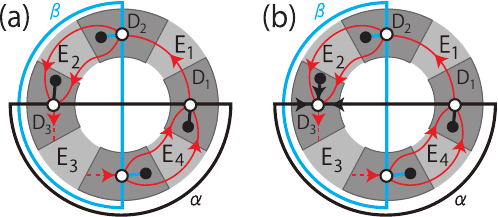}
\end{center}
\caption{ \LABEL{fig-33} An incorrect chart.}
\end{figure}

Now we have
\begin{enumerate}
\item[$\bullet$] 
for each of $E_1,E_3$, 
there are 
at least three I-edges 
of label $\alpha+1$ 
and 
at least three O-edges 
of label $\beta-1$.
\item[$\bullet$] 
for each of $E_2,E_4$, 
there are 
at least three I-edges 
of label $\beta-1$ 
and 
at least three O-edges 
of label $\alpha+1$.
\end{enumerate}

Let $e_1^*$ be 
an O-edge for $D_1$ of label $\alpha+1$.
Let $C_0$ be a simple closed curve in Int~$A$
containing the edge $e_1^*$ 
and intersecting 
each of the eight disks 
$D_1,E_1,D_2,E_2,D_3,$ $E_3,D_4,E_4$ 
by a proper arc.
The oriented edge $e_1^*$
induces the orientation of 
the simple closed curve $C_0$.
Let $\ell$ be a simple arc 
connecting a point in $\partial A$
and a point in a brim of $D^2$ 
with 
$\ell \cap A=E_1\cap D_2$.

If the simple closed curve $C_0$
is oriented clockwise 
(see Fig.~\ref{fig-34}(a)),
then apply the chart $\Gamma$ 
by DH-tricks along $\ell$ 
(see Fig.~\ref{fig-34}(b) and (c)), 
we can assume

\begin{enumerate}
\item[$\bullet$]
the simple closed curve $C_0$
is oriented counterclockwise 
(see Fig.~\ref{fig-35}).
\end{enumerate}

\begin{figure}[thb]
\begin{center}
\includegraphics{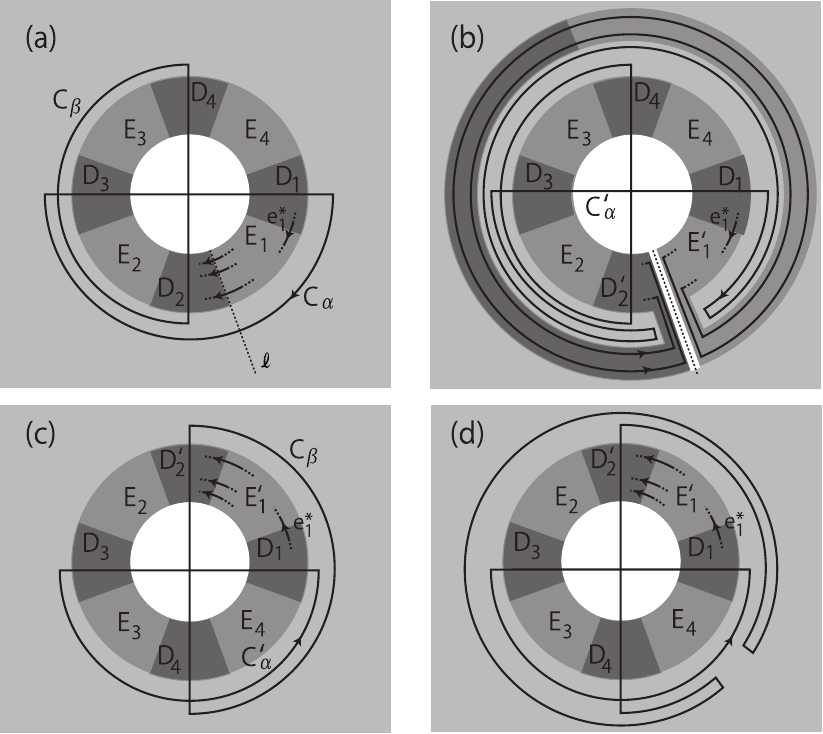}
\end{center}
\caption{ \LABEL{fig-34}Normalizing 
$C_\alpha,C_\beta$, and the annulus $A$.}
\end{figure}

Let $U$ be the connected component of $D^2-A$
containing $\partial D^2$.
Now consider the cycles 
$C_\alpha,C_\beta$ 
as non-oriented simple closed curves. 
The oriented edge $e_\alpha$ of $\Gamma_\alpha$
containing $C_\alpha\cap U$ 
induces 
the orientation of 
the simple closed curve $C_\alpha$.
Similarly 
the oriented edge $e_\beta$ of $\Gamma_\beta$
containing $C_\beta\cap U$
induces 
the orientation of 
the simple closed curve $C_\beta$.

If necessary 
we apply the chart $\Gamma$ by
a DH-trick for  
the edge $e_\alpha$,
we can assume that 
\begin{enumerate}
\item[$\bullet$] the simple closed curve 
$C_\alpha$ is 
oriented counterclockwise.
\end{enumerate}
If necessary 
we apply the chart $\Gamma$ by
a DH-trick for 
the edge $e_\beta$ 
(see Fig.~\ref{fig-34}(c) and (d))
and if necessary 
we renumber $E_1,D_1,E_2,D_2,E_3,D_3,E_4,D_4$, 
we can assume that
(see Fig.~\ref{fig-35})
\begin{enumerate}
\item[$\bullet$] 
$E_1$ does not intersect any of disks bounded by 
$C_\alpha$ nor $C_\beta$.
\end{enumerate}
Define\\
$
\delta=\left\{
\begin{array}{ll}
1&~~\text{if the simple closed curve 
$C_\beta$ 
is oriented counterclockwise,}\\
2&~~{\rm otherwise.}
\end{array}
\right.
$\vspace{2mm}

\begin{figure}[hbt]
\begin{center}
\includegraphics{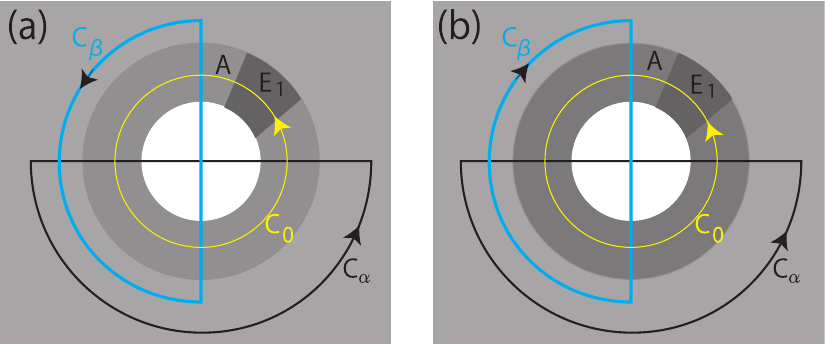}
\end{center}
\caption{ \LABEL{fig-35}
(a) $\delta=1$ (b) $\delta=2$.}
\end{figure}

{\bf Case 1}: $\beta-\alpha<2$.
Then the chart is a ribbon chart.

{\bf Case 2}: $\beta-\alpha=2$.
Then we can assume that
\begin{enumerate}
\item[(i)] 
$D_1\cup D_2\cup D_3\cup D_4=A$,
\item[(ii)] for each $i=1,2,3,4$,\\
$||{\bf a}(\Gamma,D_{i+1})||=
||{\bf b}(\Gamma,D_{i})||$,
here $D_5=D_1$.
\end{enumerate}\index{normal form $F(\Gamma)$}
We define the {\it  normal form} for the chart $\Gamma$
by\\
$F(\Gamma)=
((n,\beta-\alpha,\alpha,\delta),
(||{\bf b}(\Gamma,D_1)||,
||{\bf b}(\Gamma,D_2)||,
||{\bf b}(\Gamma,D_3)||,
||{\bf b}(\Gamma,D_4)||);\\
{\bf b}(\Gamma,D_1),
{\bf a}(\Gamma,D_2),
{\bf b}(\Gamma,D_2),
{\bf a}(\Gamma,D_3),
{\bf b}(\Gamma,D_3),
{\bf a}(\Gamma,D_4),
{\bf b}(\Gamma,D_4),
{\bf a}(\Gamma,D_1)
).$

{\bf Case 3}: $\beta-\alpha\ge 3$.
By Theorem~\ref{StITheorem2},
for each $i=1,2,3,4$ and
$j=\alpha+1,\alpha+2,\cdots,\beta-2$,
there exists 
an N-tangle 
$(\Gamma\cap D_i(j),D_i(j))$
with a label pair $(j,j+1)$
such that  
\begin{enumerate}
\item[(i)]
$A=(\cup_{i=1}^4 D_i)\cup(
\cup_{i=1}^4\cup_{j=\alpha+1}^{\beta-2} D_i(j))$,
\item[(ii)] 
$E_i=\cup_{j=\alpha+1}^{\beta-2} D_i(j)$
for each $i=1,2,3,4$,
\item[(iii)] by Theorem~\ref{StITheorem1} 
for each $i=1,2,3,4$,\\
$
||{\bf a}(\Gamma,D_{i+1})||=
||{\bf b}(\Gamma,D_i)||=
||{\bf a}(\Gamma,D_{i}(\alpha+1))||=
||{\bf b}(\Gamma,D_{i}(\alpha+1))||\\=
||{\bf a}(\Gamma,D_{i}(\alpha+2))||=
||{\bf b}(\Gamma,D_{i}(\alpha+2))||=
\cdots=
||{\bf a}(\Gamma,D_{i}(\beta-2))||\\=
||{\bf b}(\Gamma,D_{i}(\beta-2))||$,
here $D_5=D_1$.  
\end{enumerate}
We define the {\it  normal form} for the chart $\Gamma$\index{normal form $F(\Gamma)$}
by $F(\Gamma)=$\\ 
$((n,\beta-\alpha,\alpha,\delta),
(||{\bf b}(\Gamma,D_1)||,
||{\bf b}(\Gamma,D_2)||,
||{\bf b}(\Gamma,D_3)||,
||{\bf b}(\Gamma,D_4)||);\\
{\bf b}(\Gamma,D_1),
$Index$(\Gamma,D_1(\alpha+1)),
$Index$(\Gamma,D_1(\alpha+2)),
\cdots,
$Index$(\Gamma,D_1(\beta-2)),
{\bf a}(\Gamma,D_2);\\
{\bf b}(\Gamma,D_2),
$Index$(\Gamma,D_2(\beta-2)),
$Index$(\Gamma,D_2(\beta-3)),
\cdots,
$Index$(\Gamma,D_2(\alpha+1)),
{\bf a}(\Gamma,D_3);\\
{\bf b}(\Gamma,D_3),
$Index$(\Gamma,D_3(\alpha+1)),
$Index$(\Gamma,D_3(\alpha+2)),
\cdots,
$Index$(\Gamma,D_3(\beta-2)),
{\bf a}(\Gamma,D_4);\\
{\bf b}(\Gamma,D_4),
$Index$(\Gamma,D_4(\beta-2)),
$Index$(\Gamma,D_4(\beta-3)),
\cdots,
$Index$(\Gamma,D_4(\alpha+1)),
{\bf a}(\Gamma,D_1)
).$

For example, 
the normal form
for the 5-chart in Fig.~\ref{fig-36}
is\\ 
$((5,3,1,2),(8,9,8,7);\\
(1,3,4),(3,3,2),(2,3,3),(6,2);\\
(1,3,5),(2,3,4),(4,3,2),(2,3,3,1);\\
(2,4,2),(1,3,2,2),(1,2,4,1),(4,3,1);\\
(2,5),(2,3,2),(2,3,2),(5,2))$.

\begin{figure}[h]
\begin{center}
\includegraphics{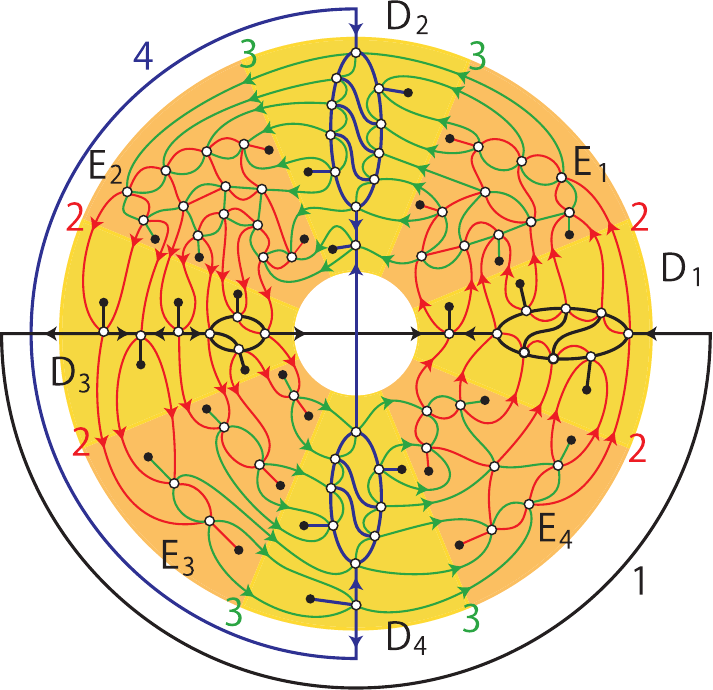}
\end{center}
\caption{ \LABEL{fig-36} A 5-chart with two crossings.}
\end{figure}

We have the following:
\begin{enumerate}
\item[(1)]
Let
$\Gamma$ and $\Gamma^*$ be 
2-crossing minimal charts.
If $F(\Gamma)=F(\Gamma^*)$ and 
$\Gamma- {\rm Main}(\Gamma)=
\Gamma^*- {\rm Main}(\Gamma^*)$, 
then  
the two charts are C-move equivalent.
\item[(2)]
There does not exist 
a 2-crossing $w$-minimal chart $\Gamma$ with 
$\alpha(\Gamma)-\beta(\Gamma)\ge 3$.
We leave the proof for this fact 
as an easy exercise.
\item[(3)]
There does not exist
a 2-crossing minimal chart
representing a surface braid whose closure 
is a 2-knot \cite[Theorem 1.2]{TwoCrossingII}.
\item[(4)]
There exists a 2-crossing minimal 4-chart (see Figure~$1$ in \cite{4-chartShpere} and Fig.~\ref{fig31}(a)).
\end{enumerate}

But we do not know 
if there exists a 2-crossing 
$c$-minimal chart with
$\alpha(\Gamma)-\beta(\Gamma)\ge 3$.

\vspace{5mm}

\begin{minipage}{65mm}
{Teruo NAGASE
\\
{\small Tokai University \\
4-1-1 Kitakaname, Hiratuka \\
Kanagawa, 259-1292 Japan\\
\\
nagase@keyaki.cc.u-tokai.ac.jp
}}
\end{minipage}
\begin{minipage}{65mm}
{Akiko SHIMA 
\\
{\small Department of Mathematics, 
\\
Tokai University
\\
4-1-1 Kitakaname, Hiratuka \\
Kanagawa, 259-1292 Japan\\
shima@keyaki.cc.u-tokai.ac.jp
}}
\end{minipage}

\vspace{1cm}

{\bf List of notations}\vspace{7mm}\\
{\small $
\begin{array}{ll}
\text{$\Gamma_m$} & p 2 \\
\text{$\alpha(\Gamma)$, $\beta(\Gamma)$} & p 2 \\
\text{${\rm Main}(\Gamma)$} & p 11 \\
\text{${\mathcal{W}}_O^{{\rm Mid}}(C,m)$} & p 12 \\
\text{$SC(X)$ } & p 15 \\
\text{${\rm Main}(\Gamma_m)$} & p 16 \\
\end{array}
~~\begin{array}{||ll}
\text{ $P[v_i,v_j]$} & p 24 \\
\text{ $\Delta(\tau)$} & p 44, p46\\
\text{ ${\bf a}(\Gamma,D)$, {\bf b}$(\Gamma,D)$} & p 45, p 47 \\
\text{ {\rm Index}$(\Gamma,D)$} & p 45, p 47 \\
\text{ $||\bf x||$} & p 47\\
\text{ $F(\Gamma)$} &  p53\\
\end{array}
$}

\newpage

\vspace{1cm}

{\bf List of terminologies}\vspace{7mm}\\
{\small $
\begin{array}{ll}
\text{$2$-crossing chart} & p 4 \\
\text{adjacent} & p 38 \\
\text{admissible tangle} & p 18 \\

\text{bigon} & p 8\\
\text{bigonal disk} & p 42\\
\text{boundary IO-arc pair} & p 5 \\
\text{brim} & p 11 \\

\text{C-move equivalent} & p 1 \\
\text{consecutive triplet} & p 11 \\
\text{corresponding edge} & p 38 \\
\text{cut-edge} & p 15 \\
\text{cycle} & p 12 \\

\text{decomposition} & p 6 \\
\text{DH-trick (double hoops trick)} & p 13 \\
\text{dichromatic path} & p 24 \\
\text{directed cycle} & p 25 \\
\text{directed path} & p 24 \\
\text{dominate} & p 27\\
\text{downward-left-selective} & p 25 \\
\text{downward maximal} & p 28\\
\text{downward net-tangle} & p 45 \\
\text{downward principal} & p 29, p 38\\
\text{downward-right-selective} & p 25\\

\text{edge sequence} & p 24 \\
\text{elementary IO-tangle} & p 6 \\

\text{free edge} & p 2\\
\text{fundamental tangle} & p 5 \\

\text{half spindle} & p 26 \\
\text{hoop} & p 2 \\

\text{I-edge} & p 26 \\
\text{inside edge} &  p 29 \\
\text{internal edge} & p 15 \\
\text{inward} & p 6, p 25 \\
\text{IO-tangle} & p 4 \\
\text{I-terminal edge} & p 44 \\

\text{leading to $e$} & p 25 \\
\text{locally inward} & p 4 \\
\text{locally left-side} & p 25 \\
\text{locally outward} & p 4 \\
\text{locally right-side} & p 25 \\
\end{array}
$}
~~{\small $
\begin{array}{||ll}
\text{loop} & p 25 \\

\text{M\&M path} & p 24  \\
\text{M2-related to $\Gamma$} & p 47 \\
\text{middle arc} & p 8 \\
\text{middle at $v$} & p 8 \\
\text{minimal chart} & p 9 \\

\text{net-tangle} & p 45 \\
\text{normal form $F(\Gamma)$} & p 53 \\
\text{NS-tangle} & p 12 \\
\text{N-tangle} & p 45 \\

\text{O-edge} & p 26 \\
\text{O-terminal edge} & p 44 \\
\text{outside edge} & p 12, p 29 \\
\text{outward} & p 6, p 25 \\
\text{oval nest} & p 9\\

\text{path of label $m$} & p 24 \\
\text{principal quad} & p 32\\
\text{proper arc} & p 5 \\

\text{ring} & p 10 \\
\text{SC-closure $SC(X)$} & p 15 \\
\text{sibling edge} & p 21 \\
\text{simple hoop} & p 2 \\
\text{simple IO-tangle} & p 5 \\
\text{situated between $X_1$ and $X_3$} & p 29 \\
\text{staple} & p 31\\
\text{starting from $e$} & p 25 \\

\text{tangle} & p 2 \\
\text{terminal edge} & p 4  \\
\text{trivial tangle} & p 6 \\
\text{Type-I} & p 5 \\
\text{Type-I$_p$} & p 6 \\
\text{Type-II} & p 5 \\
\text{Type-II$_t$} & p 19 \\

\text{upward-left-selective} & p 25\\
\text{upward maximal} & p27 \\
\text{upward net-tangle} & p 45 \\
\text{upward principal} & p 29, p38\\
\text{upward-right-selective} & p 25 \\

\text{vertex sequence} & p 24 \\
 & \\
\end{array}
$}


\end{document}